%% file: main.tex
\documentclass[twoside]{article}

\usepackage{aistats2019}
\usepackage[algo2e]{algorithm2e} 
\usepackage{algorithm}
\usepackage[utf8]{inputenc} 
\usepackage[T1]{fontenc}    
\usepackage{url}            
\usepackage{booktabs}       
\usepackage{amsfonts}       
\usepackage{nicefrac}       
\usepackage{amsmath}
\usepackage{amssymb}
\usepackage{times}
\usepackage{color}
\usepackage{amsmath,amssymb} 
\usepackage[]{algorithm2e}
\usepackage{bbm}
\usepackage{amsthm}
\usepackage{graphicx}
\usepackage{cite}

\usepackage{subcaption} 
\newtheorem{lemma}{Lemma}
\newtheorem{theorem}{Theorem}

\DeclareMathOperator*{\argmin}{arg\,min}

\SetCommentSty{mycommfont}

\usepackage[colorinlistoftodos,bordercolor=orange,backgroundcolor=orange!20,linecolor=orange,textsize=scriptsize]{todonotes}

\newcommand{\E}[1]{\mathbb{E}\left[#1\right] } 
\newcommand{\eqdef}{\overset{\text{def}}{=}} 
\newcommand{\cL}{\mathcal{L}}
\newcommand{\cG}{\mathcal{G}}
\newcommand{\mA}{{\bf A}}

\newcommand{\mD}{{\bf D}}
\newcommand{\mE}{{\bf E}}
\newcommand{\mF}{{\bf F}}

\newcommand{\mI}{{\bf I}}
\newcommand{\mJ}{{\bf J}}

\newcommand{\mR}{{\bf R}}
\newcommand{\mS}{{\bf S}}

\newcommand{\mW}{{\bf W}}

\usepackage{hyperref}       

\begin{document}

%

%


\twocolumn[

\aistatstitle{Improving SAGA via a Probabilistic Interpolation with Gradient Descent}

\aistatsauthor{ Adel Bibi \footnotemark[1] \And Alibek Sailanbayev\footnotemark[1] \And  Bernard Ghanem\footnotemark[1] \And Robert Mansel Gower\footnotemark[2] \And Peter Richt\'{a}rik\footnotemark[1]\footnotemark[3]}

\aistatsaddress{  } 
]



\begin{abstract}
\input{./sections/abstract}
\end{abstract}

\input{./sections/introduction}

\input{./sections/saga_v_gd}

\input{./sections/saga_v_minibatch}

\input{./sections/theoretical_analysis}


\input{./sections/experiments}

\input{./sections/conclusion}

\newpage

{\small
\bibliographystyle{ieeetr}
\bibliography{mybib.bib}{}
}

\input{./sections/appendix}

\end{document}

%% file: sections/abstract.tex
We develop and analyze a new algorithm for empirical risk minimization, which is the key paradigm for training supervised machine learning models. Our method---SAGD---is based on a probabilistic interpolation of SAGA and gradient descent (GD). In particular, in each iteration we take a gradient step with probability $q$ and a SAGA step with probability $1-q$. We show that, surprisingly,  the total expected complexity of the method (which is obtained by multiplying the number of iterations  by the expected number of gradients computed in each iteration) is minimized for a non-trivial probability $q$. For example, for a well conditioned problem the choice $q=1/(n-1)^2$, where $n$ is the number of data samples, gives a method with an overall complexity which is better than both the complexity of GD and SAGA.  We further generalize the results to a probabilistic interpolation of SAGA and minibatch SAGA, which allows us to compute both the optimal probability and the optimal minibatch size. While the theoretical improvement may not be large, the practical improvement is robustly present across all synthetic and real data we tested for, and can be substantial. Our theoretical results suggest that for this optimal minibatch size our method achieves linear speedup in minibatch size, which is of key practical importance as minibatch implementations are used to train machine learning models in practice. Moreover, empirical evidence suggest that a linear speedup in minibatch size can be attained with a parallel implementation.

%% file: sections/introduction.tex
\section{Introduction}
The success of modern machine learning models and applications has been 
pushing forward the development of algorithms for tackling large scale optimization problems.
In particular, the development of stochastic methods for minimizing objective functions with a finite sum structure, such as the empirical risk minimization (ERM) problem \cite{shai_book}.
In the ERM problem each element of the finite sum structure corresponds to the loss of the model evaluated at a particular data point:
\begin{equation}
    \label{original_problem}
    \min_{x \in \mathbb{R}^d} f(x) \eqdef \frac{1}{n} \sum_{i=1}^n f_i(x).
\end{equation}

In the big data regime, $n$ is typically very large. This finite sum structure can be leveraged to design efficient optimization algorithms. In this paper, we focus on the case where each $f_i$ is $L_i$--smooth and  $f$ is $\mu$--strongly convex, e.g., the regularized ridge regression $f_i(x) = \frac{1}{2}(a_i^\top x - y_i)^2 + \frac{\lambda}{2}\|x\|_2^2$, regularized logistic loss $f_i(x) = \frac{1}{2}\log(1 + \exp(-y_i a_i^\top x)) + \frac{\lambda}{2}\|x\|_2^2$ and regularized conditional random fields~\cite{SchmidtBADCS15}.

{\bf Variance-reduced methods.} In the last five years, {\em variance-reduced} stochastic gradient-type methods were devised for solving \eqref{original_problem}, improving upon the running time of the immensely popular stochastic gradient descent (SGD) method \cite{RobinsMonro1951, Nemirovski-Juditsky-Lan-Shapiro-2009}. The first such method was SAG  \cite{schmidt2017minimizing}, which operates by keeping track of the average of the latest stochastic gradients associated with each data point. The method uses a biased estimate of the gradient, and while its analysis is notoriously very technical, it enjoys a linear rate of convergence on problem \eqref{original_problem} as opposed to a sublinear rate of SGD.
This was soon followed by SDCA \cite{shalev2013stochastic}, which requires an explicit strongly convex regularizer to be present in \eqref{original_problem}, and operates by applying the coordinate descent method \cite{UCDC} to the dual problem. The final complexity rate is the same as that of SAG. One of the most popular variance reduced methods is 
SVRG \cite{johnson2013accelerating} (see also S2GD \cite{S2GD}), which relies on a few computations of the full gradient of $f$ in an outer loop, and using this gradient to correct the standard SGD estimator in an inner loop, and leads to an estimator with diminishing variance. In a similar spirit to both SAG and SVRG, SAGA \cite{defazio2014saga} is yet another incremental gradient-based method with the same fast linear convergence rate. SAGA modifies SAG so as to make the gradient estimates unbiased, which leads to a much more streamlined analysis. Similar to SAG, SAGA maintains a table of the latest stochastic gradients computed for each data point, which is then used via an averaging and relaxation step to produce the unbiased gradient estimator. 

It is well known that all of the above methods can be improved using several orthogonal tricks,  such as mini-batching \cite{pegasos2, PCDM, PCD_complexity, mS2GD, ASDCA, dmSDCA, csiba2016importance, hofmann2015variance}, importance sampling \cite{UCDC, richtarik2016optimal, QUARTZ, ALPHA, csiba2016importance, allen2016even},  acceleration \cite{APPROX, Hydra2,  ASDCA, lin2014accelerated, ALPHA, allen2017katyusha, allen2016even}, and higher-order information \cite{Schmidt2009,
Schmidt2011a,
Erdogdu2015nips, GowerGold2016,
moritz2016linearly, GRB-HessianSVRG,
acceleratedBFGSrules}. These approaches lead to further speedup and practical benefits. However, we will not elaborate further on these contributions.

{\bf Contributions.} 
Recently, Gower, Richt\'{a}rik and Bach~\cite{2018arXiv180502632G} proposed a novel family of variance reduced methods, called JacSketch, based on iterative Jacobian sketching and controlled stochastic optimization reformulations of \eqref{original_problem}. 
Our work is inspired by their results. In particular, we look more deeply into some exotic special cases of JacSketch not considered in \cite{2018arXiv180502632G} and uncover an interesting phenomenon, which we capture theoretically, and also observe in practice.

In particular,  we propose and analyze a new algorithm that  interpolates between SAGA and GD in a probabilistic sense. That is, we choose a probability $q\in[0,1]$, and then take a gradient step with  probability $q$ and a SAGA step with probability $(1-q)$  . We prove that by choosing $q$ dependent on the condition number of the problem~\eqref{original_problem} the resulting method always has a better \emph{total} complexity than either SAGA or GD. Indeed let $\kappa := \left.4\overline{L}\right/\mu$  be this condition number where $\overline{L} : = \tfrac{1}{n}\sum_{i=1}^n L_i.$ If problem~\eqref{original_problem} is \emph{well conditioned} as captured by the inequality $ \kappa \leq n-1$, then the simple choice of $q = \left.1\right/(n-1)^2$ leads to a method with better total complexity than both SAGA and GD. On the other hand, if problem~\eqref{original_problem} is \emph{badly conditioned} with $ \kappa \geq n-1$, then choosing $q = \left. \mu \right/4n \overline{L}$ gives a method that also enjoys a strictly better total complexity than both SAGA and GD. These results are captured by Theorem~\ref{theo:t-n}.

Motivated by this result, we then modify our method by replacing the GD step with a minibatch SAGA step, with minibatch size $\tau$. Since  minibatch SAGA with $\tau=n$ reduces to GD, this is a generalization of our previous method. However, we now have an extra free parameter, $\tau$, and can optimize for it as well. We compute the optimal $\left(q^*,\tau^*\right)$  pair minimizing the total complexity of the algorithm using the more relaxed approximation $L_i = L_{\max} := \max_{i=1,\ldots, n} L_i$ of the smoothness constants.
Our results show  that  the optimal minibatch size often is non-trivial, that is, it is neither $1$ nor $n$. The resulting method is  better than SAGA and minibatch SAGA with any minibatch size, since the choice $(q,\tau) = (0,1)$ reduces to SAGA and the choice $(q,\tau) = (1,\tau)$ reduces to minibatch SAGA.  We show that computing the optimal $(q^*,\tau^*)$ pair can be performed cheaply. 

Lastly, we conduct extensive experiments on both synthetic and real datasets to demonstrate that our method is better in practice than both SAGA and GD. We also demonstrate experimentally that the theory matches the experiments in predicting the optimal ($\tau^*$)  given $q^*$. The experiments also suggest that any parallel implementation of our method achieves  {\em linear speedup} in minibatch size $\tau$. SAGD is the first variance-reduced method for \eqref{original_problem} with this property. Linear or superlinear speedup in minibatch size was previously only demonstrated for the SDNA method \cite{SDNA}, which applies to the dual of a regularized version of \eqref{original_problem}. However, the super-linear speedup of SDNA is for iteration complexity, and does not take the cost of each iteration into consideration.

%% file: sections/saga_v_gd.tex
\section{Probabilistic interpolation of SAGA and gradient descent}
Loosely inspired by SVRG \cite{johnson2013accelerating}, which computes a gradient step in an outer loop and then performs variance reduced SGD steps in an inner loop, we propose the following modification of SAGA:
\begin{equation} 
\label{eq:sagd update}
\begin{aligned}
x^{k+1} = \begin{cases}
\text{Gradient step at } x^k &\text{ with prob} \ q,\\
\text{SAGA step at } x^k     & \text{ with prob} \ (1-q),
\end{cases}
\end{aligned}
\end{equation}
where $q\in[0,1]$ and the stepsizes are $\alpha$. Every so often (i.e., with ``small'' probability $q$),  the current table (matrix) of gradient estimates $\mJ^k \in \mathbb{R}^{d \times n}$ is reset to consist of the true gradients $\mJ^k \leftarrow \begin{bmatrix}f^\prime_1(x^k), f^\prime_2(x^k), \dots, f^\prime_n(x^k)\end{bmatrix}$, where $f^\prime_i(.) \in \mathbb{R}^d$, and with probability $(1-q)$ one column of the table is updated such that is $\mJ^k_{:i} = f^\prime_i(x^k)$. 
However, we observe two major differences between the proposed algorithm and SVRG: (i) the number of inner and outer loop iterations are not fixed but are randomly controlled by $q$, and (ii) with this new method, all the samples that are used to build the unbiased estimator of the gradient are used to update the full gradient estimate. This is in contrast with SVRG, where only one of the samples of the gradient step update in the inner loop is later used to compute the total gradient in the outer loop.

Note that a somewhat similar albeit fundamentally different strategy is employed in the CSGS method \cite{SCSG} of Lei and Jordan, which modifies SVRG so that the outer loop estimates the gradient via a minibatch, as opposed to the full batch. From this perspective, we work with SAGA instead of SVRG, and include a probabilistic interpolation with minibatch SAGA in order to reduce variance more aggressively every so often, instead of doing this in an outer loop, as in CSGS. 

To analyze this interpolation scheme, we consider the total complexity $\Omega^{q,n} = \Omega_k^{q,n} \times \Omega_c^{q,n} $ where $\Omega_k^{q,n}$ denotes the iteration count and $\Omega_c^{q,n}$ is the average complexity per iteration. Given that computing the $f_i'$ gradient is the dominating cost, it is clear that under such a sampling $\Omega_c^{q,n} = q(n-1) + 1$. That is, on average $q(n-1) + 1$ gradients are computed per iteration. In what follows we showcase that in fact there exists a non trivial simple choice for $q$ that is better than both SAGA and GD.
    \begin{theorem}
\label{theo:t-n}
For the update \eqref{eq:sagd update}, consider the Lyapunov function
\begin{align*}
&\Psi^k \eqdef \|x^k - x^* \|_2^2 + \frac{\alpha}{2 L_{\max} (q(n-1)+1)}\|\mJ^k - \nabla \mF(x^*) \|^2_{F},
\end{align*}
where $\mJ^k$ is the table of gradients and $\nabla F(x^*)=\begin{bmatrix} f^\prime_1(x^*), f^\prime_2(x^*), \dots, f^\prime_n(x^*)\end{bmatrix}$.
If the problem is well conditioned with $\frac{4 \overline{L}}{\mu} \leq n-1$, then by choosing $q = \frac{1}{(n-1)^2}$ and using the stepsize
\begin{equation}\label{eq:alphawell}
    \alpha= \frac{1}{4(1-\frac{2}{n})L_{\max}+ \mu(n-1)},
\end{equation}
the expected total complexity to achieve $\epsilon$ accuracy, that is $\mathbb{E}[\Psi^k] \leq \epsilon \Psi^0$, is given by
\begin{eqnarray}
    \Omega^{\frac{1}{(n-1)^2},n} 
    &= &  
    \left(n+ \frac{n-2}{n-1} \frac{4L_{\max}}{\mu}\right)\log \left(\frac{1}{\epsilon}\right).
    \label{eq:total_compleixty_tau-n98js8}
\end{eqnarray}
On the other hand, if the problem is badly conditioned with $\frac{4 \overline{L}}{\mu} \geq n-1$, then by choosing $q= \frac{\mu}{4n \overline{L}}$ and using the stepsize 
\begin{equation}\label{eq:alphabad}
\alpha = \frac{1}{4 \overline{L}} \frac{\left((4\overline{L}+\mu)-\frac{\mu}{n}\right)^2}{\mu n((4\overline{L}+\frac{\mu}{n})-\mu)+4L_{\max}(4\overline{L} - \frac{\mu}{n})},
\end{equation}
the total complexity is
\begin{eqnarray}
    &\Omega^{\frac{\mu}{4n \overline{L}},n} 
    =\left(n + \frac{4L_{\max}}{\mu} \left( 1 - \frac{ 1}{\frac{4\overline{L}}{\mu}+1-\frac{1}{n}} \right)\right)\log \left(\frac{1}{\epsilon}\right).
    \label{eq:total_compleixty_tau-nY(q)}
\end{eqnarray}
    \end{theorem}

The details of the proof are left to the \textit{supplementary material}.
Comparing  the total complexity of SAGA, which is given by $\Omega^{0,n} =(n +\left. 4L_{\max} \right/ \mu)\log(1/\epsilon)$, to the total complexity results ~\eqref{eq:total_compleixty_tau-n98js8} and~\eqref{eq:total_compleixty_tau-nY(q)}, we see that {\em the total complexity of our method is always better than that of SAGA}. Indeed, this is trivially true for~\eqref{eq:total_compleixty_tau-n98js8} and as for~\eqref{eq:total_compleixty_tau-nY(q)} it follows from $\frac{4\overline{L}}{\mu} \ge n-1$ that
\begin{equation} \label{eq:m98s893s}
\frac{4 \overline{L}}{\mu} \geq n-1 \Rightarrow
\frac{4\overline{L}}{\mu}+1-\frac{1}{n}   \geq 1 \Leftrightarrow 0 \leq \frac{ 1}{\frac{4\overline{L}}{\mu}+1-\frac{1}{n}} \leq 1.\end{equation}
Thus the simple choice of the interpolation parameter $q$ given in Theorem~\ref{theo:t-n} exhibits strictly better performance than both SAGA and GD.
Though we note that both complexity results~\eqref{eq:total_compleixty_tau-n98js8} and~\eqref{eq:total_compleixty_tau-nY(q)} converge to the complexity of SAGA  asymptotically as $n\rightarrow \infty.$ Indeed this is again trivial for~\eqref{eq:total_compleixty_tau-n98js8}. As for~\eqref{eq:total_compleixty_tau-nY(q)} we have from $\frac{4\overline{L}}{\mu} \ge n-1$ that
  $$1 - \frac{ 1}{\frac{4\overline{L}}{\mu}+1-\frac{1}{n}} {\geq}  1 - \frac{ 1}{n-1+1-\frac{1}{n}} = 1 - \frac{ 1}{n-\frac{1}{n}} \underset{n\rightarrow \infty}{ \longrightarrow} 1, $$
and consequently $\Omega^{\frac{\mu}{4n \overline{L}},n}$ is asymptotically lower bounded by the complexity of SAGA. Finally since $\Omega^{\frac{\mu}{4n \overline{L}},n} $ is always upper bounded by the complexity of SAGA (which follows from~\eqref{eq:m98s893s}) we have that
\[\lim_{n\rightarrow \infty} \Omega^{\frac{\mu}{4n \overline{L}},n}  \quad = \quad  (n +\left. 4L_{\max} \right/ \mu)\log(1/\epsilon).\]

 It is salient to note that we do not claim that such $q$ is optimal, it is only an example that a non-trivial interpolation exists that improves upon SAGA and GD. Such a result motivates the direction towards finding the optimal $q^*$ that directly minimizes the total complexity $\Omega^{q, n}$ in a more general setting.
 
\vspace{10pt}

%% file: sections/saga_v_minibatch.tex
\begin{figure*}[t]
\begin{subfigure}[h]{0.99\linewidth}
\centering
\includegraphics[width=0.35\textwidth]{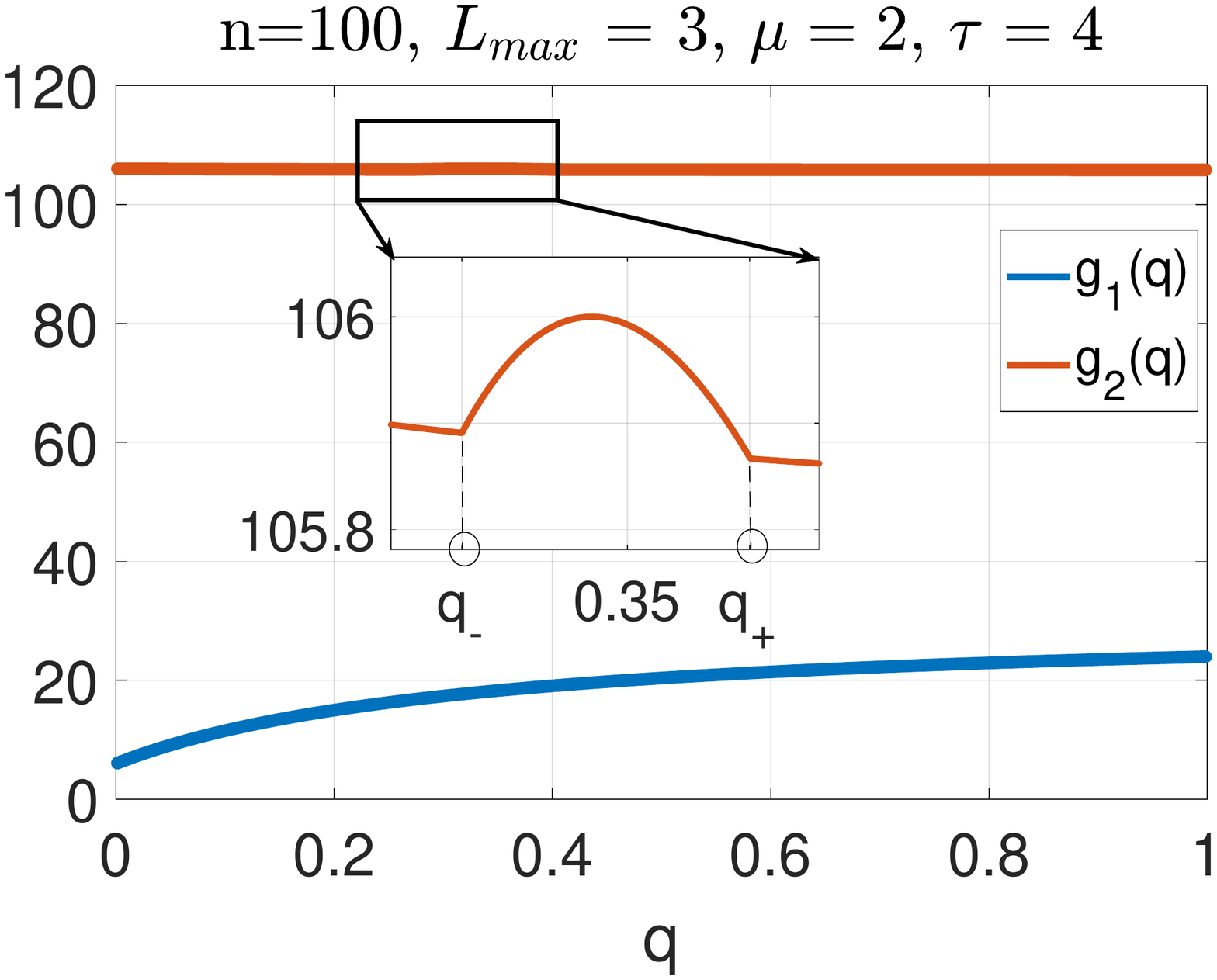}
\includegraphics[width=0.35\textwidth]{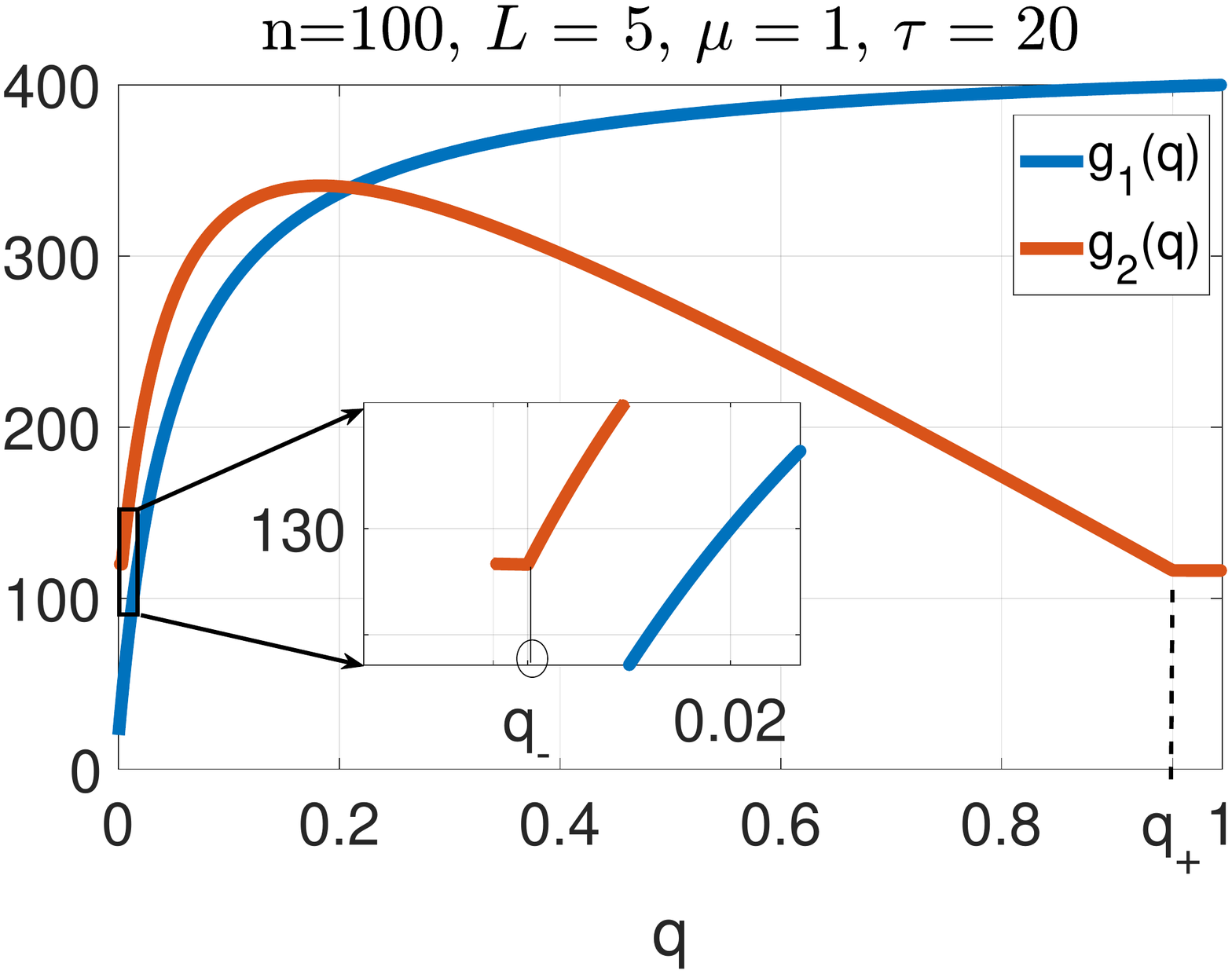}
\end{subfigure}
\vspace{-1.25cm}
\caption{Shows two different examples of the functions $g_1$ and $g_2$. As proven in Lemmas~\ref{g1_increase} and~\ref{g2_behaviour}, $g_1$ is monotonically increasing $\forall q \in [0,1]$ while $g_2$ is concave on $[q_-,q_+]$ and decreasing elsewhere. It is also to be noted that the two functions intersect at atmost one point. }
\vspace{-0.5cm}
\label{plots_of_total_complexities}
\end{figure*}

\begin{algorithm}[t]
\SetAlgoLined
\SetKwInput{KwInit}{Initilization}
\KwInit{$x^0~\in~\mathbb{R}^d, \; q~\in~[0, \, 1],\; \text{batch size }\tau~\in~\mathbb{N},$  $\theta~=~\frac{n}{q(\tau-1)+1}$, stepsize $\alpha \in \mathbb{R}_+$ and randomly initialize table $\mJ~\in~\mathbb{R}^{d \times n}$}
\While{error > $\epsilon$}{
    \textbf{(i)} with probability $\frac{q}{\binom{n}{\tau}}$, sample $\tau$ indices $I_k$\\
    Calculate $f^\prime_j(x^k)$ for $j \in I_k$ \\
    \begin{equation}
      x^{k+1} \leftarrow x^k - \alpha \left[\frac{\theta}{n} \sum_{j \in I_k}\left( f^\prime_j(x^k) - \mJ_{:j}\right) + \frac{1}{n} \sum_{i=1}^{n}  \mJ_{:i} \right]
        \nonumber
    \end{equation} 
    \tcp{minibatch SAGA step}
    Update the table: $\mJ_{:j} =
    \begin{cases}
      \mJ_{:j}, \text{ if } j \not\in I_k \\
      f^\prime_j(x^k), \text{ if } j \in I_k
    \end{cases}
    $ \\
    \noindent\textbf{(ii)} with probability $\frac{1-q}{n}$, select a sample $j$ uniformly at random;\\
    Calculate $f^\prime_j(x^k)$ \\
    $x^{k+1} \leftarrow x^k - \alpha \left[\frac{\theta}{n} \left(f^\prime_j(x^k) - \mJ_{:j}\right) + \frac{1}{n} \sum_{i=1}^n \mJ_{:i}\right]$  \tcp{SAGA step} 
    Update the table: 
        $\mJ_{:i} =
    \begin{cases}
      \mJ_{:i}, \text{ if } i \not= j \\
      f^\prime_i(x^k), \text{ if } i=j
    \end{cases}$
} 
\caption{SAGD: Interpolation between SAGA and minibatch SAGA}
\label{SAGD_algo}
\end{algorithm}

\vspace{-15pt}

\section{Probabilistic interpolation of SAGA and minibatch SAGA}

Now let us consider the more general method that interpolates between SAGA and minibatch SAGA:
\begin{equation} 
\label{eq:sagd_batch update}
\begin{aligned}
x^{k+1} = \begin{cases}
\text{$\tau$-minibatch SAGA step at } x^k   &  \text{with prob} \ q,\\
\text{SAGA step at } x^k &    \text{with prob} \ 1-q. \\
\end{cases}
\end{aligned}
\end{equation}
Note that when $\tau = n$, the model in \eqref{eq:sagd_batch update} reduces to the interpolation between gradient descent and SAGA described in \eqref{eq:sagd update}. 
In this section, we investigate the values of $(q, \tau)$ that achieve the best total complexity.
We are particularly interested in this question since there are currently no guidelines on how to choose the minibatch size for SAGA type methods. It has been observed in practice that using a small minibatch for SAGA is better in terms of effective passes over the data, as compared to sampling a single data point at each iteration (that is, the standard SAGA method).  But using too large a minibatch can degrade the performance. Thus the need for theoretical guidelines.

\indent The complexity of our algorithm is detailed in the following theorem:

\label{compleixty}
\begin{theorem}(Convergence of SAGD)
\label{main_theorem}
Consider Algorithm~\ref{SAGD_algo} for solving problem \eqref{original_problem}. Assume that the   functions $f_i$ are $L_i$--smooth, and that $f$ is strongly convex with parameter $\mu>0$.  Consider the Lyapunov function:
\begin{align*}
\Psi^k \eqdef \|x^k - x^* \|_2^2 + \frac{\theta \alpha}{2 n L_{\max}}\|\mJ^k - \nabla \mF(x^*) \|^2_{F},
\end{align*}
where
$\theta = \frac{n}{q(\tau-1)+1}$. Let 
\begin{equation} \label{eq:cL1new}
    \begin{aligned}
     \cL_1 = \frac{1}{(q(\tau-1)+1)^2} &\left(\left( q \left(\frac{\tau(n-\tau)}{n-1}-1\right) +1 \right) L_{\max} \right. \\
     &\left.+ n q 
     \frac{\tau(\tau-1)}{n-1}\overline{L}\right),
     \end{aligned}
\end{equation}
which is known as the expected smoothness constant, and let
\begin{align}
\rho = \begin{cases}
\theta^2 \left(\frac{1-q}{n} + q \frac{\tau}{n} \frac{n-\tau}{n-1}\right), & q\theta^2 \le \frac{n}{\tau}\frac{n-1}{\tau-1} \\
\theta^2 \left(\frac{1-q}{n} + q \frac{\tau}{n}\frac{n-\tau}{n-1}\right) &   q\theta^2 \geq \frac{n}{\tau}\frac{n-1}{\tau-1}.
\\
\qquad + n \left(\theta^2 q \frac{\tau}{n}\, \frac{\tau-1}{n-1}-1\right), &
\end{cases} \label{eq:rhointheo}
\end{align}
Choose $0\leq q \leq 1$ and $\tau \in [n]$. If the stepsize is given by
\begin{align}
\alpha = \min \left\{\frac{1}{4 \cL_1}, \frac{n}{4L_{\max}\rho + \mu \theta n}\right\},
\label{alpha_value}
\end{align}
then the expected total complexity in order to achieve $\epsilon$ accuracy, that is $\mathbb{E}[\Psi^k] \leq \epsilon \Psi^0$, is
\begin{align}
    \Omega^{q,\tau} &= \left(q(\tau-1)+1\right) \nonumber \\
    &\max \left\{\frac{4 \cL_1}{\mu}, \theta + \frac{4 \rho L_{\max} }{\mu n} \right\}  \log \left(\frac{1}{\epsilon}\right).
    \label{total_compleixty_SAGD}
\end{align}
\end{theorem}
The proof is left for the \textit{supplementary material}.

%% file: sections/theoretical_analysis.tex
\section{Theoretical analysis}

In this section, we revisit the total complexity result from Theorem~\ref{main_theorem} and derive the optimal ($q^*,\tau^*$)  pair that minimizes the total complexity \eqref{total_compleixty_SAGD}. In order to simplify the expression in \eqref{total_compleixty_SAGD}, we use the lax estimate  $L_i = L_{\max}$ of the individual smoothness constants for $\forall i = 1 \dots n$. First we note that
\begin{lemma}
\label{max_of_prod}
$\frac{4\cL_1}{\mu} \ge 0$, $\theta + \frac{4 \rho L_{\max} }{\mu n} \ge 0$ and $q(\tau-1) + 1 \ge 0~~~ \forall q\in[0,1]$ and $~\forall \tau$.
\end{lemma}
In light of Lemma~\ref{max_of_prod}, the total complexity \eqref{total_compleixty_SAGD} can be rewritten in the form of $\Omega^{q,\tau} = \max\left\{g_1,g_2\right\}\log(1/\epsilon)$ where $g_1 \eqdef \left(\frac{4\cL_1}{\mu}\right)\left(q(\tau - 1) + 1\right)$ and $g_2 \eqdef \left(\theta + \frac{4 \rho L_{\max} }{\mu n}\right)\left(q (\tau - 1) + 1\right)$. To obtain a $\tau$ and $q$ that optimize the complexity, we need to analyze how a change in $\tau$ or $q$ affects the functions $g_1$ and $g_2$.
\begin{lemma}
\label{g1_increase}
The function $g_1$ is monotonically increasing in $q \in [0,1]$.
\end{lemma}
\begin{lemma}
\label{g2_behaviour}
The function $g_2$ is concave in $q$ in the interval $[q_-,q_+]$ and monotonically decreasing in $q \in [0,q_-] \cup [q_+,1]$, where $q_{\pm} = \frac{n\tau + 2(1-n) \pm \sqrt{n\tau}{\sqrt{4(1-n) + n\tau}}}{2(n-1)(\tau-1)}$.
\end{lemma}
Note that $q_{\pm}$ are the solutions to the polynomial $q\theta^2 = \frac{n}{\tau}\frac{n-1}{\tau-1}$ in which $\rho$ changes sign. Moreover, it is easy to show that $q_{\pm}$ are valid probabilities through the following Lemma:
\begin{lemma}
For $\tau \ge 4$, $q_{\pm}$ are respectively increasing and decreasing functions in $\tau$; therefore, 
\begin{align}
\frac{1}{(n-1)^2} \leq q_{-} &\leq \frac{n+1-2\sqrt{n}}{3(n-1)} \\
&\leq \frac{1}{3} \leq \frac{n+1 + 2\sqrt{n}}{3(n-1)} \leq q_{+} \leq 1 .\notag
\end{align}
\end{lemma}
The previous Lemma asserts that $q_-$ and $q_+$ are valid probabilities for all values of $\tau$. At last, we analyze the points of intersections between $g_1$ and $g_2$ in the following Lemma.
\begin{lemma}
\label{q_intersects_and_tau_range}
The functions $g_1$ and $g_2$ intersect at exactly $q_{i1} = \frac{n-1}{(\tau-1)(\tau \frac{4L_{\max}}{\mu}+1 -n)}$ for  $\tau \in \left[\tau_{\text{min}},\tau_{\text{max}}\right]$ and at $q_{i2} = \frac{n - \frac{4L_{\max}}{\mu}}{\frac{4L_{\max}}{\mu}\left(\tau - 1\right)}$ for $\tau \in \left(\tau_{\text{max}},n\right]$ where $\tau_{\text{min}} = \frac{n}{\frac{4L_{\max}}{\mu}} +1 - \frac{\mu}{4L_{\max}}$, $\tau_{\text{max}}={ \min\left(\bar{\tau}_{\text{max}},n\right)\mathbbm{1}_{n > \frac{4L_{\max}}{\mu}} + n\mathbbm{1}_{n \leq \frac{4L_{\max}}{\mu}}}$ and lastly $\bar{\tau}_{\text{max}} = \left(\frac{n(n-1)\mu}{(n - \frac{4L_{\max}}{\mu})4L_{\max}}\right)$. 
\end{lemma}
Note that $\mathbbm{1}_{n > \frac{4L_{\max}}{\mu}}$ is an indicator function such that $\mathbbm{1}_{n > \frac{4L_{\max}}{\mu}} = 1$ if $n > \frac{4L_{\max}}{\mu}$ and zero elsewhere. It is important to note that for a given $\tau$ the two functions intersect at exactly one point. This is due to the fact that $q_{i1}$ and $q_{i2}$ occur in a non-overlapping intervals of $\tau$. Lastly, and as an illustration to the behaviour of the two functions $g_1$ and $g_2$, we plot in Figure~\ref{plots_of_total_complexities} both functions $g_1$ and $g_2$ under two different choices of $n$, $d$, $L_{\max}$, $\mu$ and $\tau$. Note that $g_1$ and $g_2$ exhibit the behaviour described in Lemmas~\ref{g1_increase} and~\ref{g2_behaviour}. Now, we are ready to present the third main result.

\begin{figure*}[t]
\begin{subfigure}[ht]{0.99\linewidth}
\includegraphics[width=0.33\textwidth]{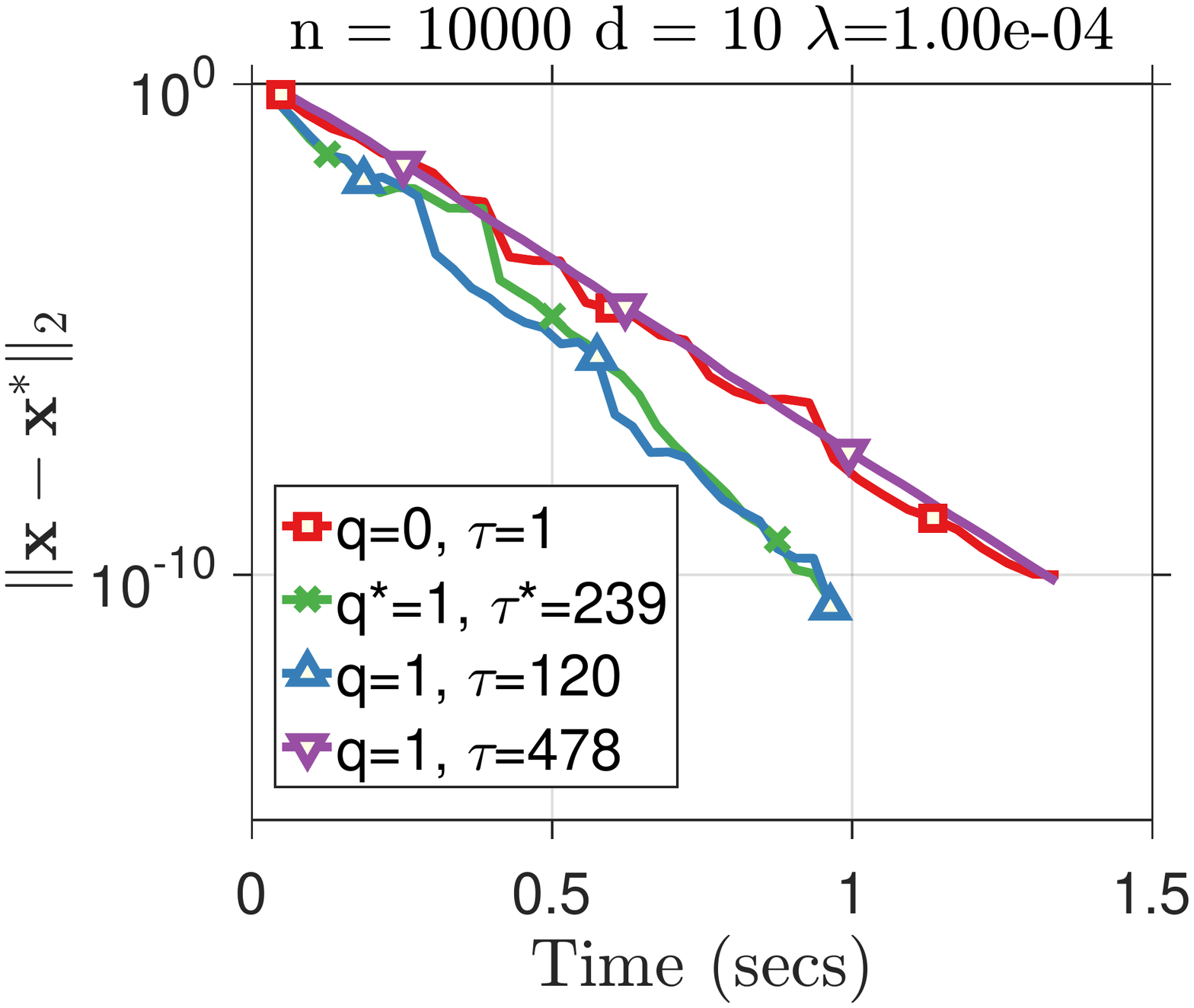}
\includegraphics[width=0.33\textwidth]{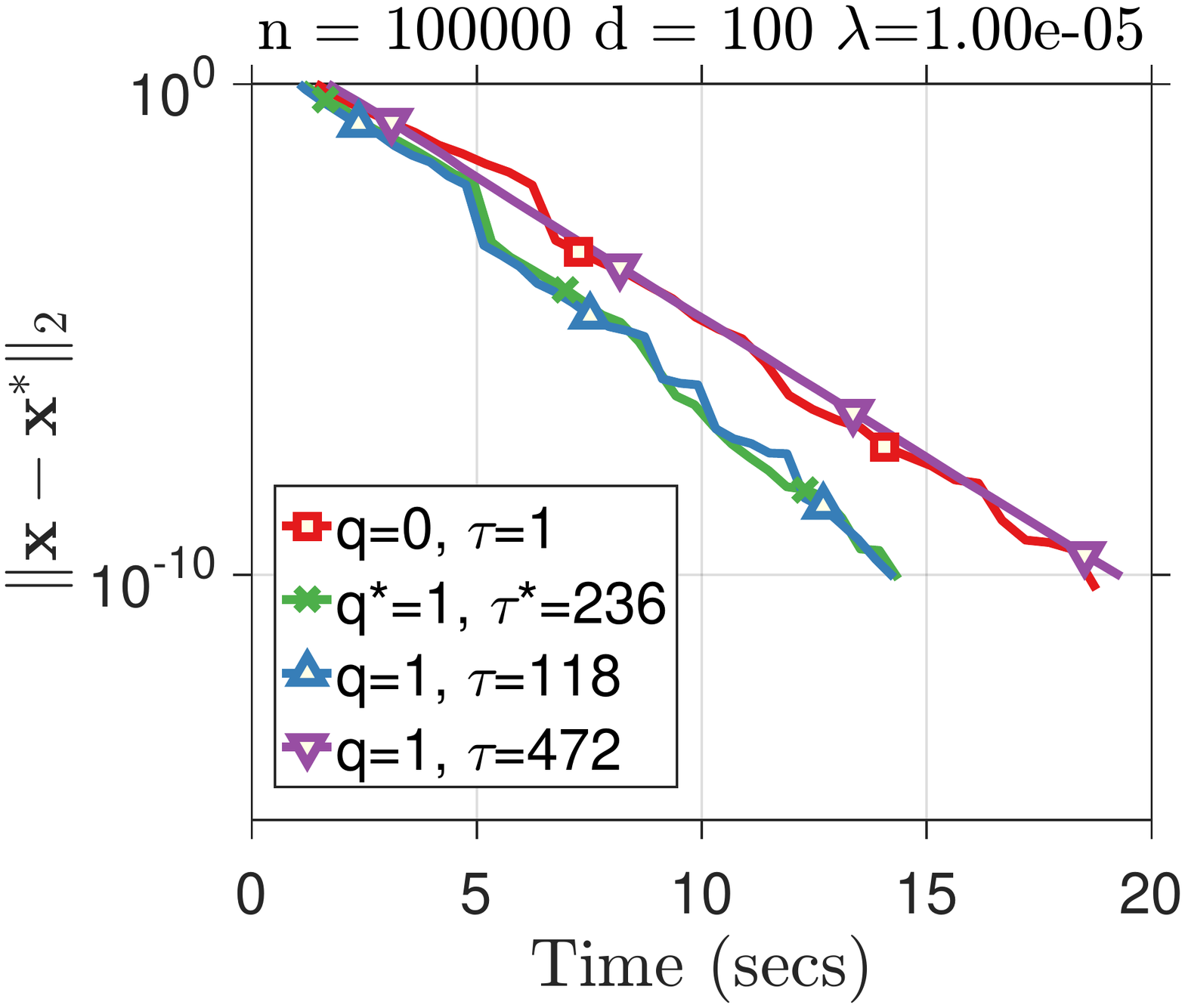}
\includegraphics[width=0.33\textwidth]{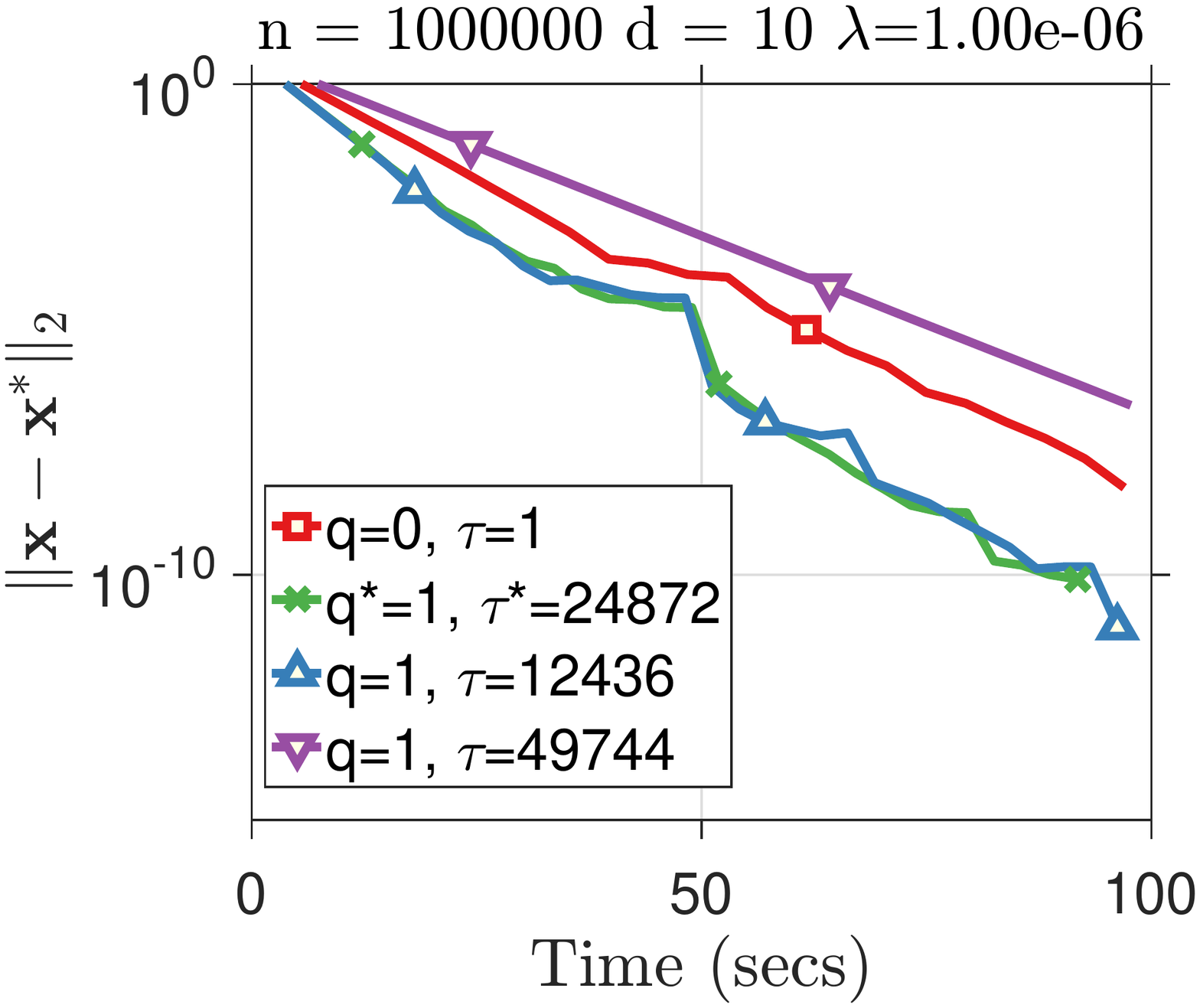}
\vspace{-2.5cm}
\end{subfigure}
\begin{subfigure}[ht]{0.99\linewidth}
\includegraphics[width=0.33\textwidth]{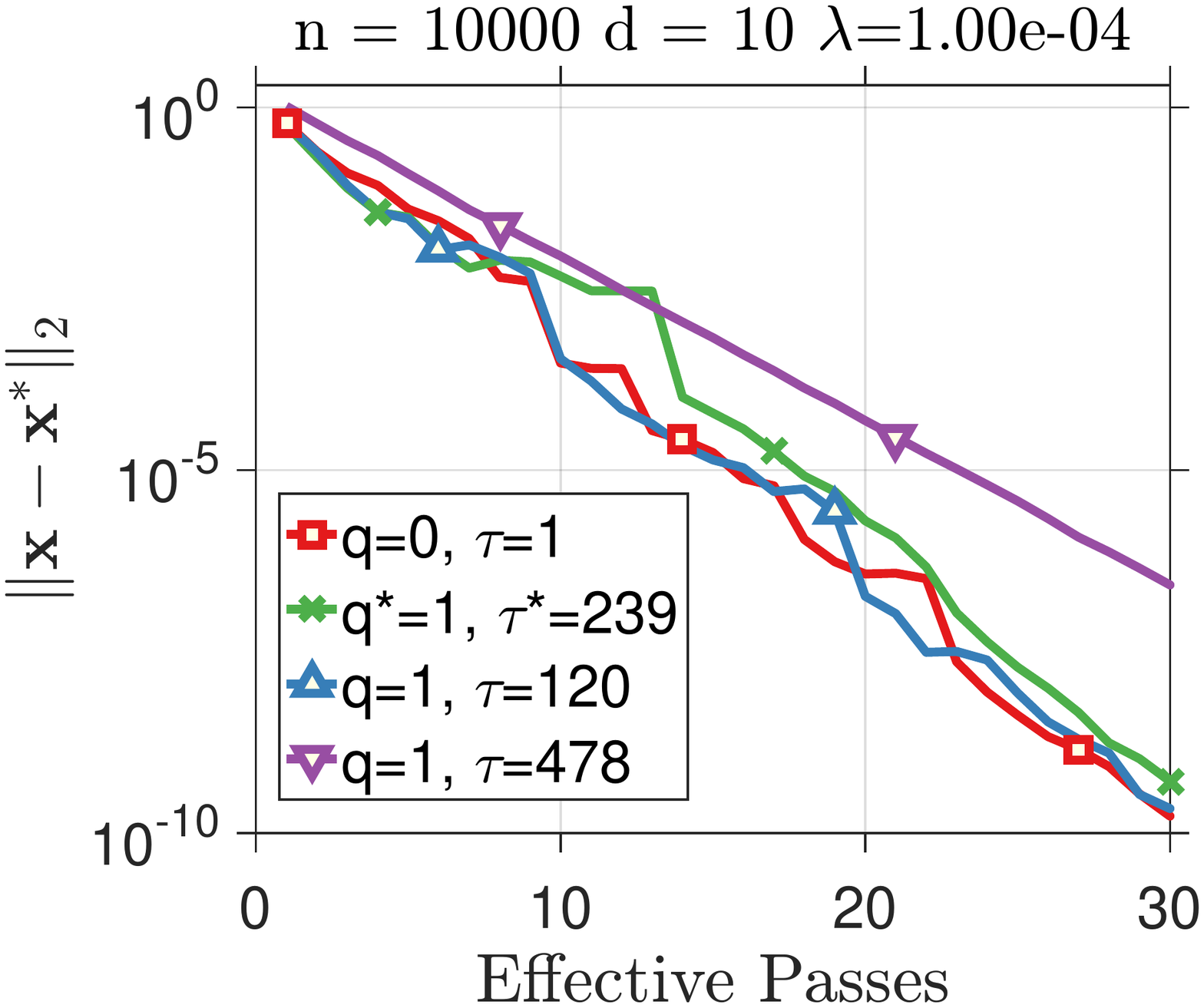}
\includegraphics[width=0.33\textwidth]{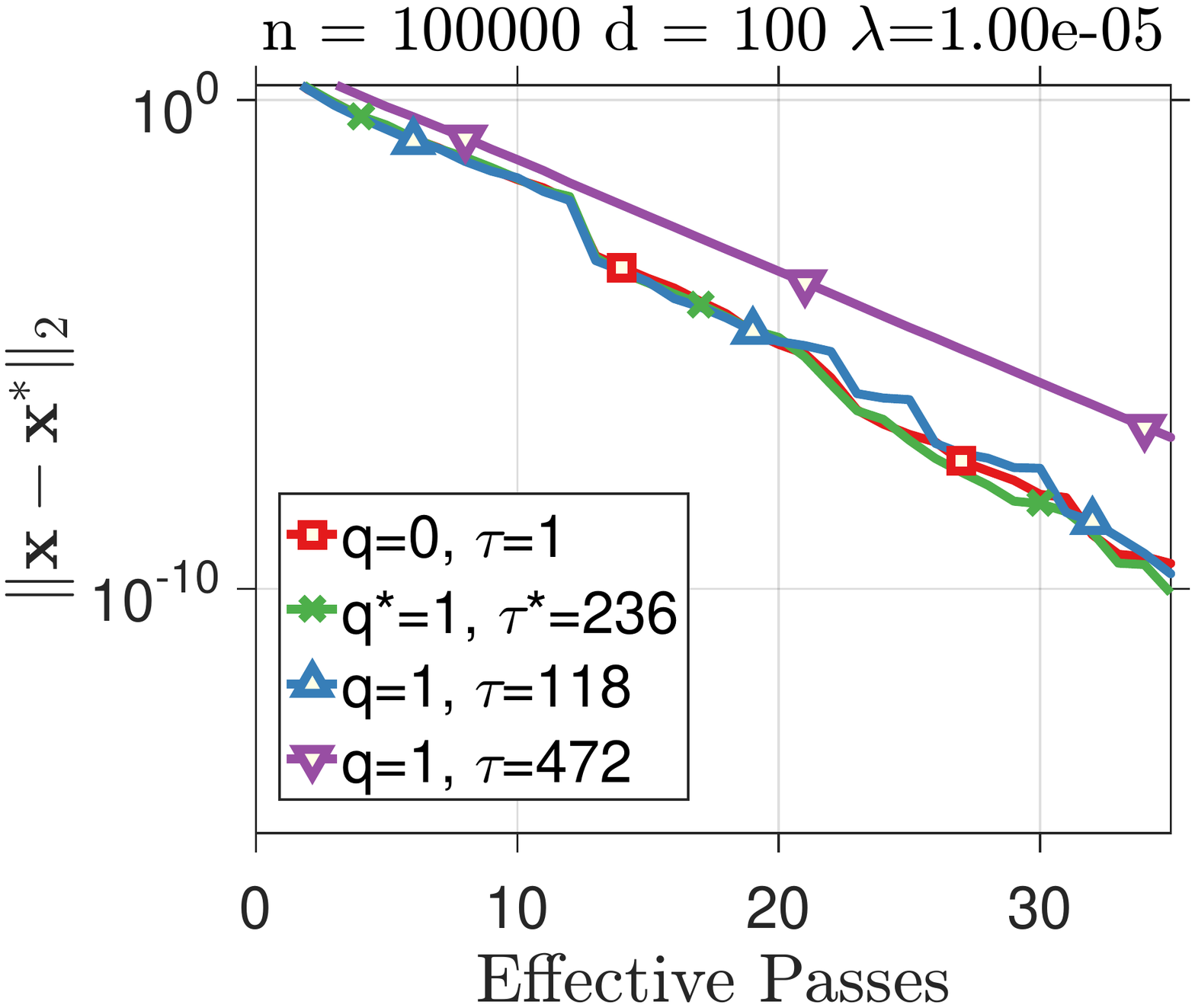}
\includegraphics[width=0.33\textwidth]{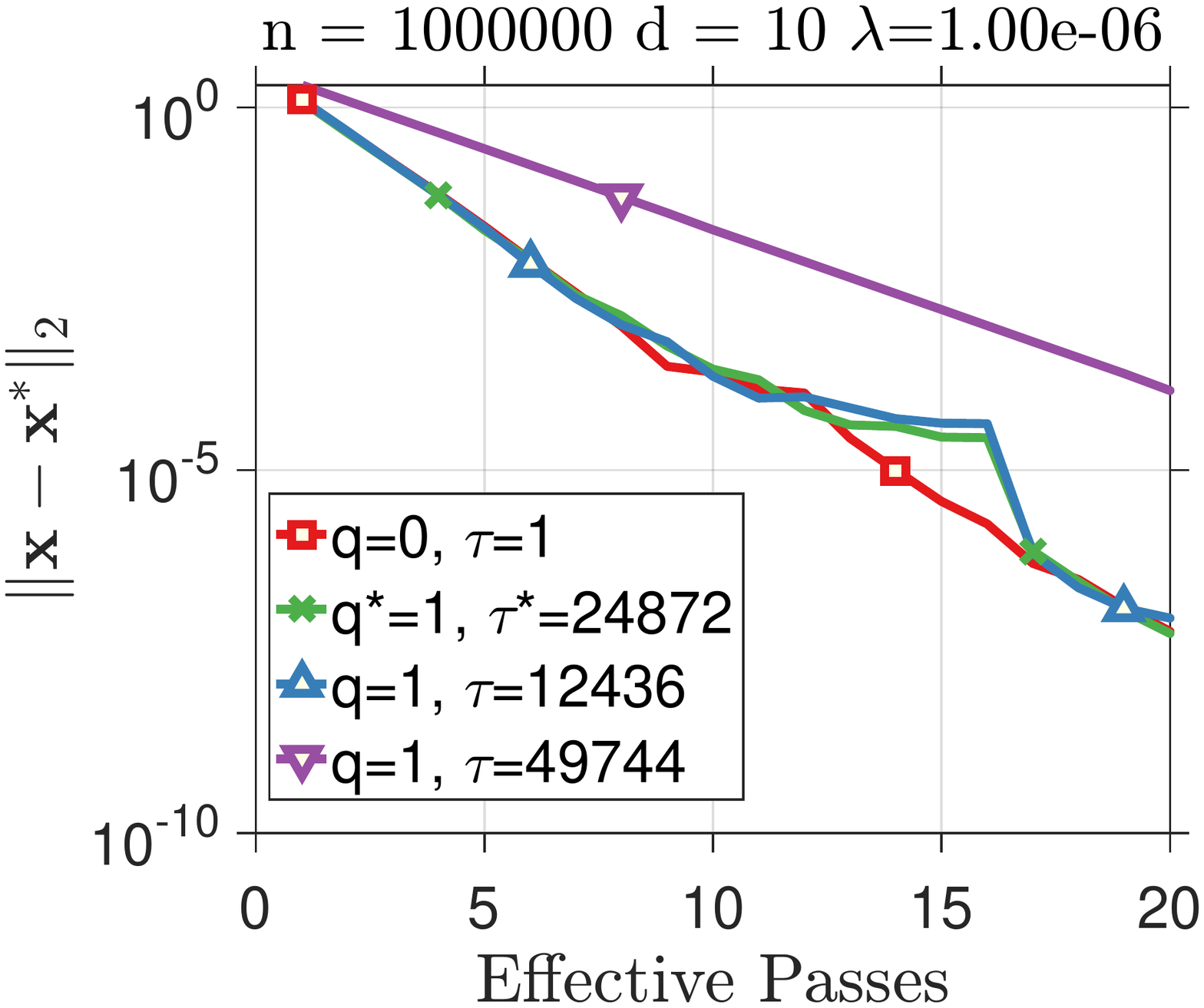}
\vspace{-1.75cm}
\end{subfigure}
\caption{Shows a comparison between SAGA and our method with optimal ($q^*,\tau^*$) pair. We also compare them both against ($q^*,\frac{\tau^*}{2}$) and ($q^*,2\tau^*$). The first row shows a comparison in time while the second row shows the same experiment as compared in the total number of effective passes.}
\label{first_time_epocs_comp}
\vspace{-0.5cm}
\end{figure*}

\begin{theorem}
The optimal $q^*$ of the total complexity in \eqref{total_compleixty_SAGD} is either one of the intersection points ($q_{i1}$, $q_{i2}$) or $q_-$ or $q=1$. That is $q^* \in \left\{1,q_-,q_{i1},q_{i2}\right\}$; therefore,
\begin{align}
    &\tau^* = \argmin_{\tau} \big\{\Omega^{1,\tau}, \Omega^{q_-,\tau}, \Omega^{q_{i1},\tau},\Omega^{q_{i2},\tau}\big\} \label{optimal_tau}\\
    &q^* = \argmin_q \big\{\Omega^{1,\tau^*}, \Omega^{q_-,\tau^*}, \Omega^{q_{i1},\tau^*},\Omega^{q_{i2},\tau^*}\big\}
\end{align}
\label{optimal_q_tau_theorem}
\end{theorem}
\vspace{-15pt}
\begin{proof}The proof is a direct implication of the previous lemmas.\end{proof}
\vspace{-10pt}\textbf{Discussion.} A key consequence of Theorem~\ref{optimal_q_tau_theorem} is that SAGA ($q=0$) is not optimal in the sense of total complexity. Moreover, Theorem~\ref{optimal_q_tau_theorem} provides a consistent way to compute the optimal ($q^*,\tau^*$) that interpolates between SAGA and minibatch SAGA that minimizes the total complexity. 
Such a pair is computed once before the algorithm. Moreover, solving \eqref{optimal_tau} may seem to be expensive since it requires the computation of $4n$ scalar values. However, the computation is in fact much smaller than $4n$. As discussed in Lemma~\ref{q_intersects_and_tau_range}, the intervals of $\tau$ in which the intersection points ($q_{i1},q_{i2}$) occur are in fact non-overlapping. That is, computing $q_{i1}$ and $q_{i2}$ require at most $n$ scalar computations. Moreover, for $q=1$ the total complexity is $\Omega^{1,\tau} = \max\left\{g_1=\frac{4L_{\max}\tau}{\mu},g_2=n + \frac{4L_{\max}}{\mu} \frac{n-\tau}{n-1}\right\}\log(1/\epsilon)$. Note that $g_1$ and $g_2$ are monotonically increasing and decreasing in $\tau$, respectively. Since for $q=\tau=1$, $g_2 > g_1$; therefore, the optimal minibatch size $\tau^*$ occurs at the intersection point $\tau = \text{round}\left(1 + \frac{\mu(n-1)}{4 L_{\max}}\right)$. Therefore, we only need to compute a single scalar for when $q=1$. Moreover, this result is particularly important as it provides a formula for the \textit{optimal} minibatch size of SAGA. Lastly, as for $q_-$, it is only computed for when $\tau \ge 4$. Therefore, the total number of scalar computations required is at most $2(n-1)$ and they are computed once before running the algorithm. Therefore, \eqref{optimal_tau} can be simplified to:

\begin{equation}
\begin{aligned}
    \tau^* = \argmin_{\tau} &\left\{\Omega^{1,\text{round}\left(1 + \frac{\mu}{4 L_{\text{max}}}(n-1)\right)
},\Omega^{q_-,\tau\ge4}, \right. \\
& \left.\Omega^{q_{i1},\tau \in \left[\tau_\text{min},\tau_\text{max}\right]},\Omega^{q_{i2},\tau \in \left(\tau_{\text{max}},n\right]} \right\}.
\end{aligned}
\end{equation}

\begin{figure*}[t]
\begin{subfigure}[ht]{0.99\linewidth}
\includegraphics[width=0.33\textwidth]{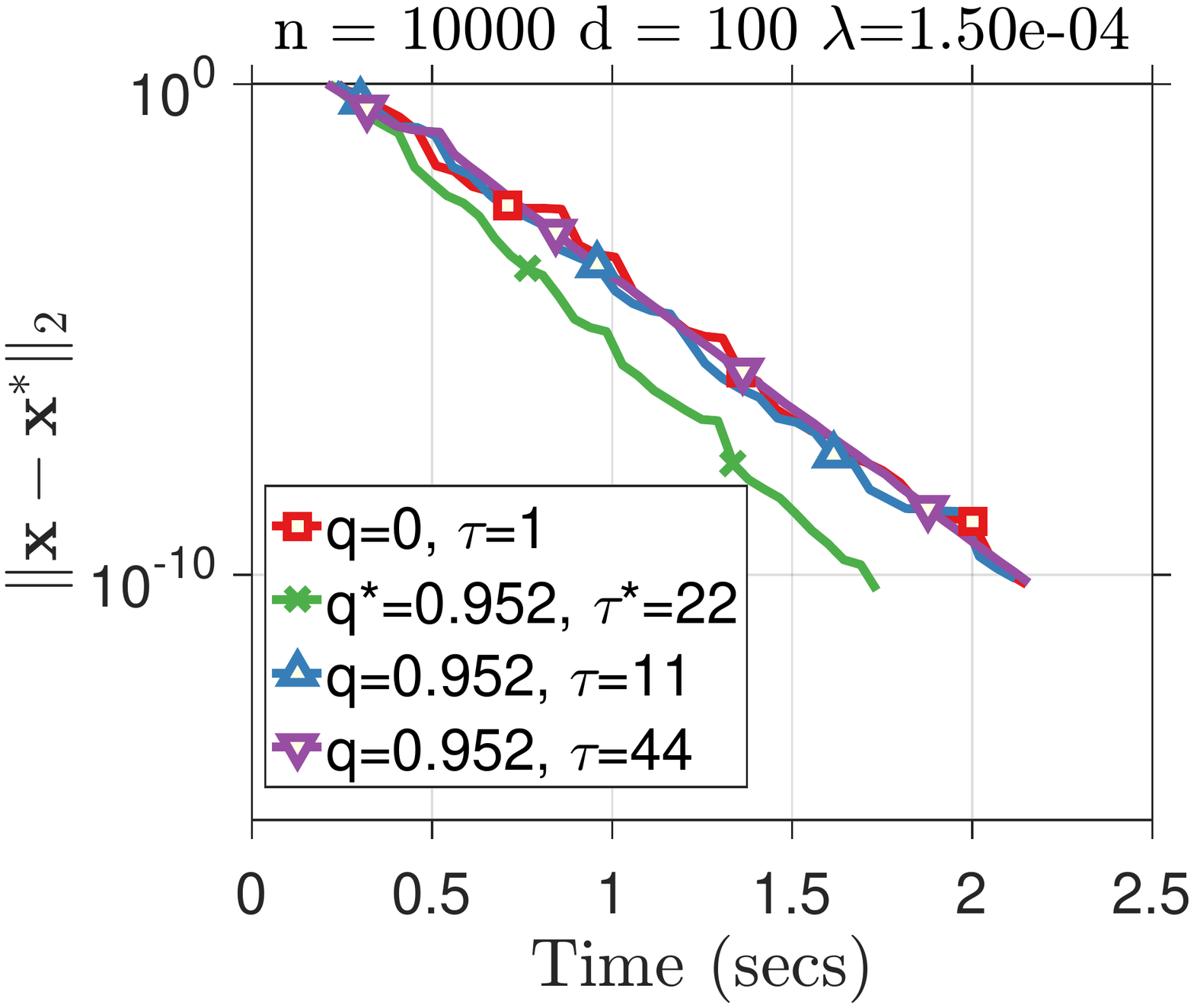}
\includegraphics[width=0.33\textwidth]{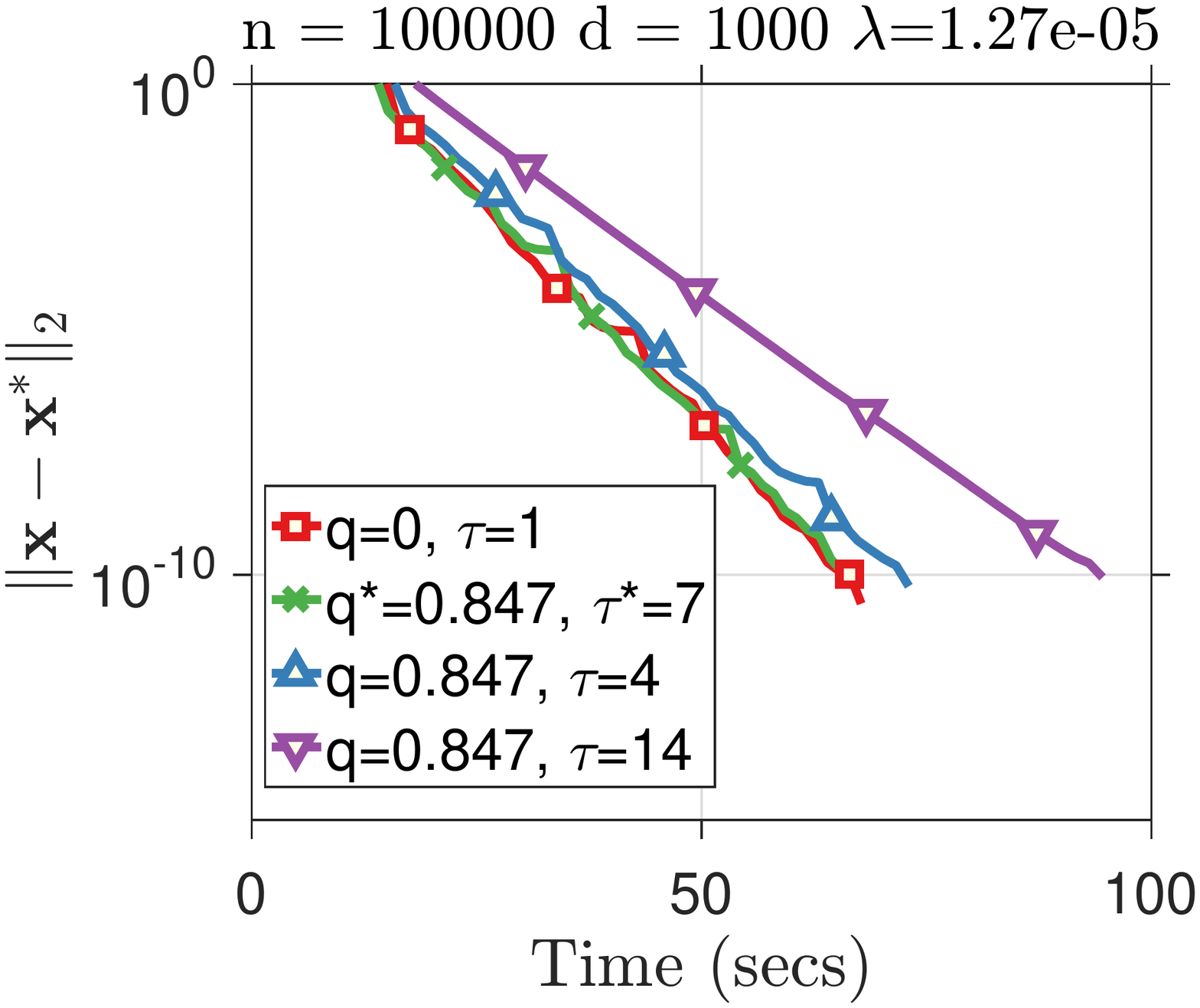}
\includegraphics[width=0.33\textwidth]{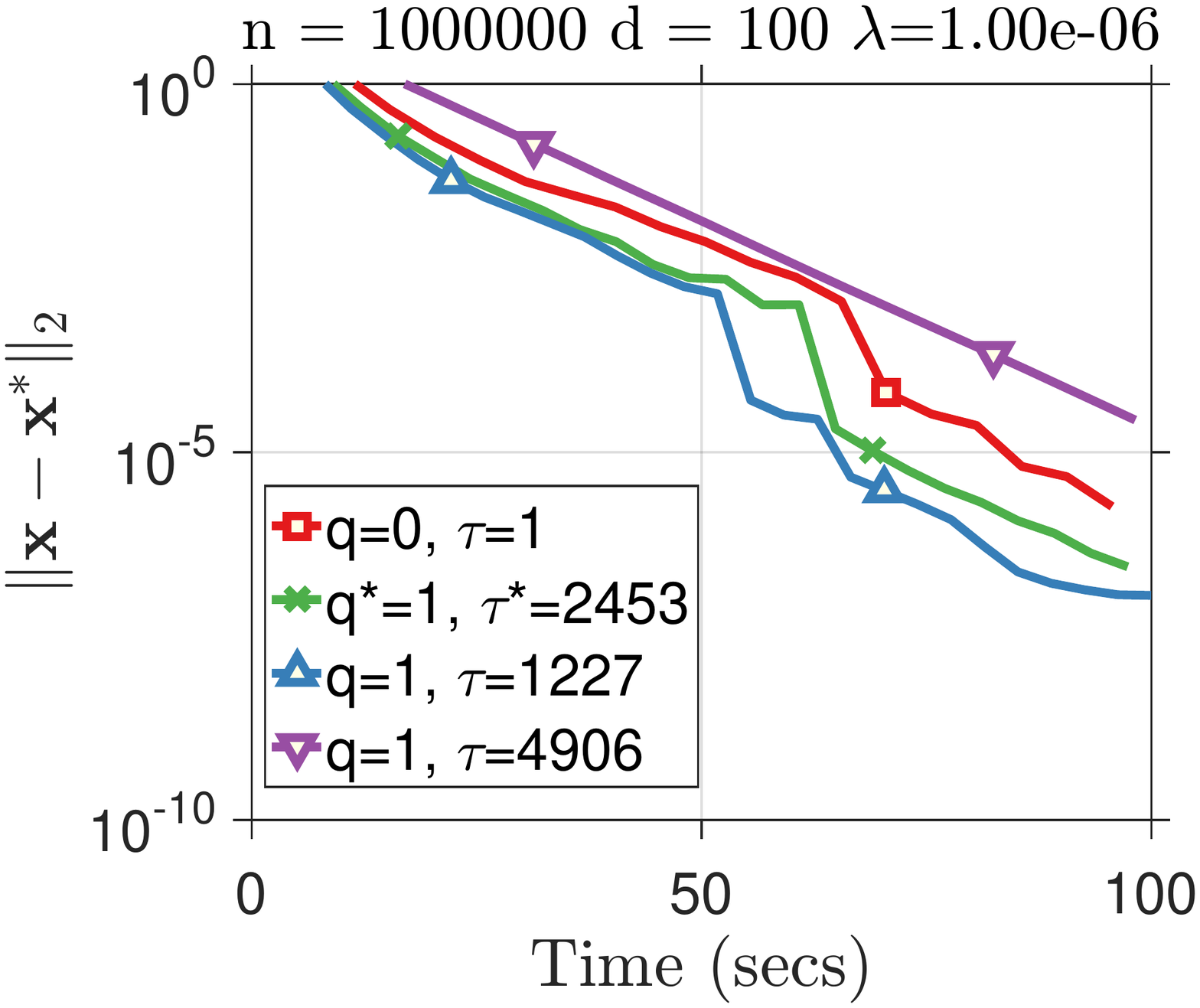}
\vspace{-2.5cm}
\end{subfigure}
\begin{subfigure}[ht]{0.99\linewidth}
\includegraphics[width=0.33\textwidth]{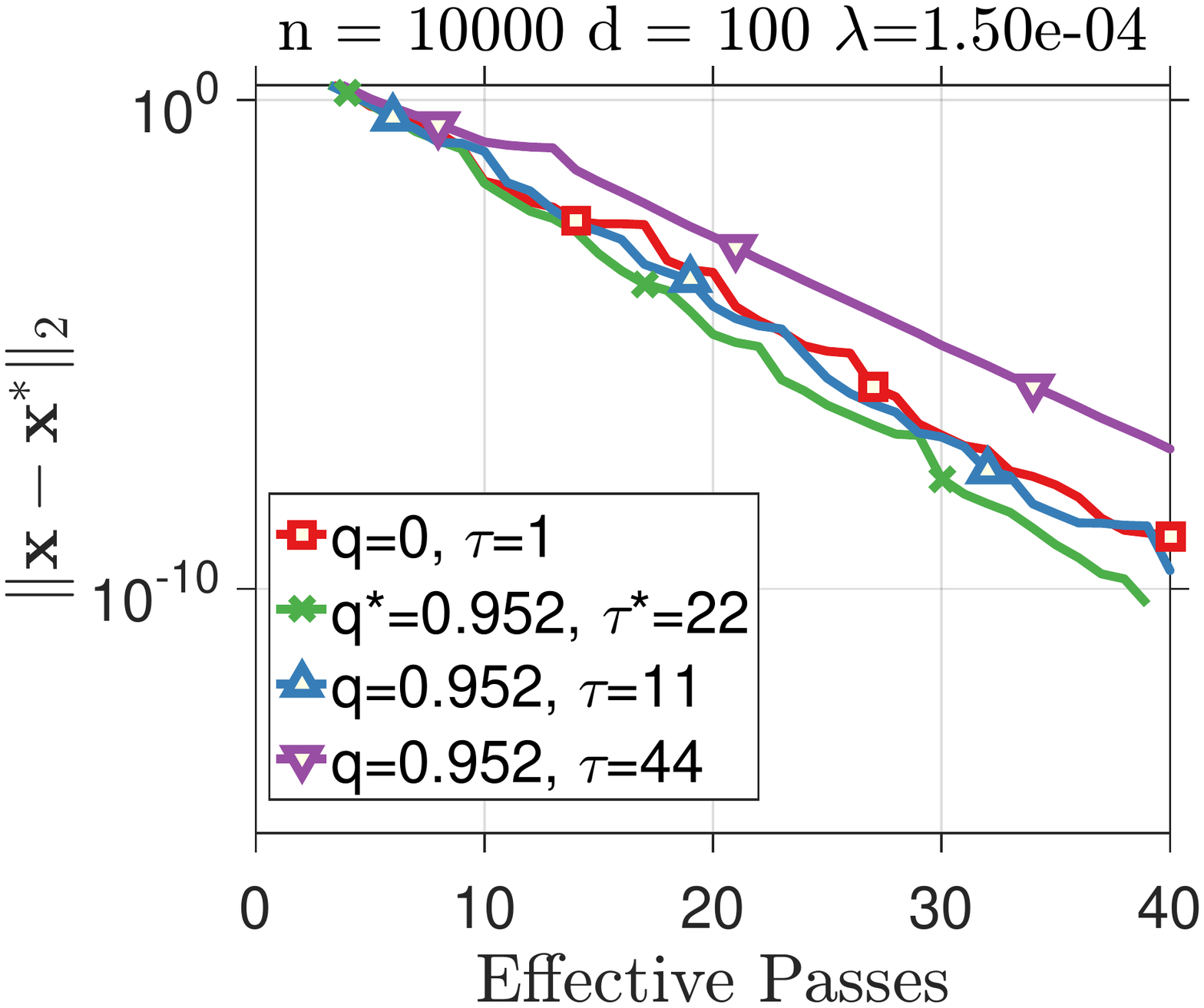}
\includegraphics[width=0.33\textwidth]{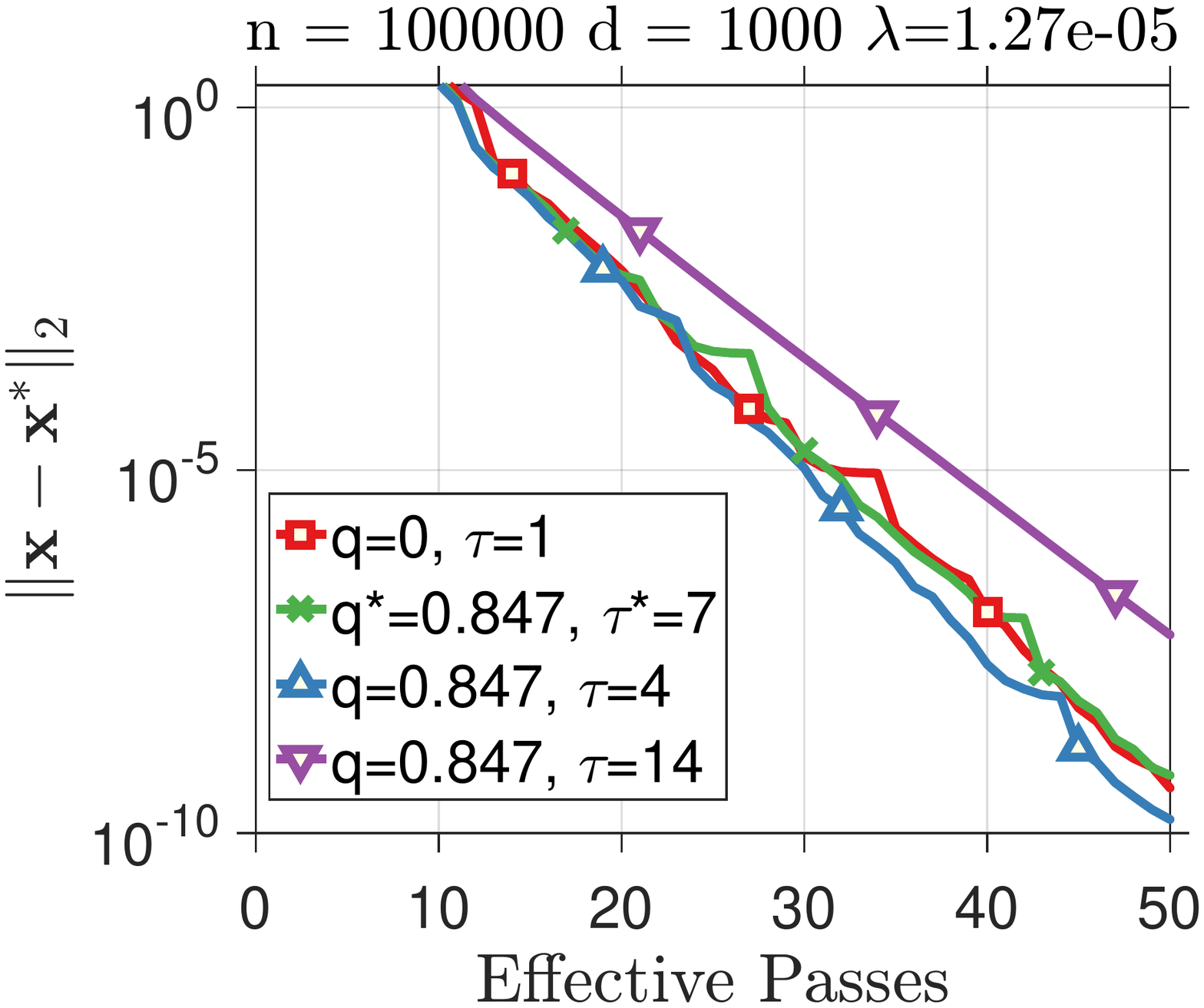}
\includegraphics[width=0.33\textwidth]{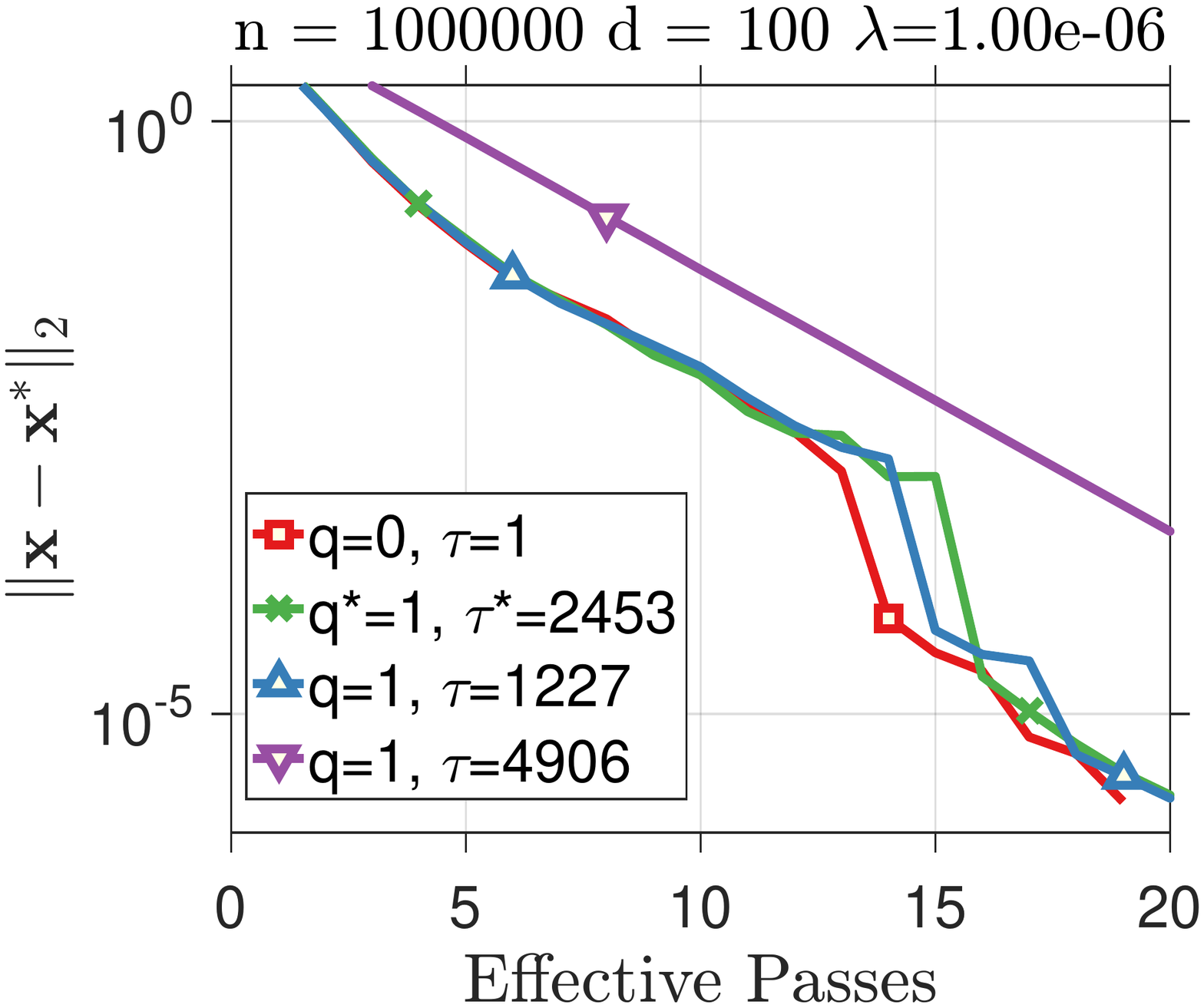}
\vspace{-1.75cm}
\end{subfigure}
\caption{Shows a comparison between SAGA and our method with optimal ($q^*,\tau^*$) pair. We also compare them both against ($q^*,\frac{\tau^*}{2}$) and ($q^*,2\tau^*$). The first row shows a comparison in time while the second row shows the same experiment as compared in the total number of effective passes.}
\label{second_time_epocs_comp}
\vspace{-0.25cm}
\end{figure*}

%% file: sections/experiments.tex
\section{Experiments}
\label{experiments}
In this section, we conduct extensive experiments comparing the performance of SAGD to SAGA as a baseline. The comparisons against GD were omitted as GD was too slow compared to SAGA. The experiments were conducted on several synthetic and real datasets demonstrating that not only SAGD is better in practice but also that SAGD can achieve linear speed up with minibatch size $\tau^*$ for a parallel implementation. We have also conducted experiments to verify how accurately the theory predicts the practical optimal batch size $\tau^*$. In this section, we consider the problem of regularized ridge regression:
\begin{equation}
    \begin{aligned}
    f(x) = \frac{1}{2n} \|\mA x - y \|_2^2 + \frac{\lambda}{2}\|x\|_2^2,
    \label{linear_system}
    \end{aligned}
\end{equation}
where $\mA \in \mathbb{R}^{n \times d}$, $y \in \mathbb{R}^n$ are the data and $\lambda > 0$ is the regularization parameter. However, further experiments on the regularized logistic loss are left for the \textit{supplementary material}. The code was implemented in MATLAB16a and all experiments were run on a dual processor 2.66GHz machine. Any experiments that compare time discarded the time required to read data and the time to perform random sampling. This is due to the assumption that all such data can be hashed and achieved with $\mathcal{O}(1)$ complexity. This will eliminate any dependency on the platform. As for the synthetic data, the rows of $\mA$ and $y$ were sampled from the standard Gaussian distribution $\mathcal{N}(0,1)$ with a regularization parameter set to $\lambda = \frac{1}{n}$. In some few experiments, the data $\mA$ and $y$ were sampled from a uniform distribution with somewhat different choice of $\lambda$. The details and reasons behind such setups for the experiments will be detailed later. In all experiments, the rows of $\mA$ were normalized such that $\|\mA(i,:)\|_2 = 1 ~~\forall i$. Note that for  ridge regression, we have $L_i = \|\mA(i,:)\|_2^2 + \lambda = 1 + \lambda ~~\forall i$ and that $\mu = \frac{1}{n}\lambda_{\text{min}}\left(\mA^\top \mA\right) + \lambda$. In all experiments, we compute a high accuracy solution $x^*$ by solving the linear system in \eqref{linear_system}. The metric used in all experiments is the  error $\|x - x^*\|_2$.

\textbf{Synthetic Experiments.}
In the synthetic experiments, we compare SAGD to SAGA across different choices of $n$ and $d$. As can be noted in Figure~\ref{first_time_epocs_comp}, SAGD with optimal interpolation factor and minibatch size ($q^*$,$\tau^*$) consistently outperforms SAGA in time. The same comparison was also conducted against several other non optimal choices of minibatch size, namely ($2\tau^*$ and $\tau^*/2$). In practice, SAGD with ($q^*$,$\tau^*$) performs very similarly to SAGD with ($q^*$,$\tau^*/2$). However, SAGD with ($q^*$,$\tau^*$) is still favorable as it allows for larger batch sizes suitable for parallelization. Figure~\ref{first_time_epocs_comp} also shows a comparison among SAGA, SAGD with ($q^*$,$\tau^*$), and SAGD with non optimal minibatch sizes ($\tau^*/2$ and $2\tau^*$) on the number of effective data passes. It is clear that SAGD with ($q^*, 2\tau^*$) always requires a larger number of effective data passes to achieve the same accuracy. More interestingly, both SAGA and SAGD with both ($q^*,\tau^*$) and ($q^*,\tau^*/2$) can achieve the same $\epsilon$  accuracy in the same number of effective data passes. This suggests that upon having a parallel implementation of SAGD, one can achieve a linear speedup in minibatch size $\tau^*$ outperforming both competitors. This is the first time linear speedup in minibatch size is obtained for a variance reduced gradient-type method by directly solving the primal empirical risk minimization problem. Similarly, we conduct even further synthetic experiments but on larger dimensional problems. Figure~\ref{second_time_epocs_comp} demonstrates that SAGD with ($q^*, \tau^*$) is either slightly faster or similar to SAGA in time. However, similar to the conclusion drawn in Figure~\ref{first_time_epocs_comp}, SAGD with ($q^*$,$\tau^*$) can have significant speed ups outperforming SAGA in a parallel implementation. As most of the previous experiments had $q^*$ as optimal interpolation factor, we showcase the existence for other non trivial $q^*$. For instance, if we set $\lambda = 3/(2n)$ for the dataset with $n=10000$ and $d=100$, the problem's condition number $\frac{L}{\mu}$ changes that the optimal probabilistic interpolation is non trivially $q^*=0.952$. This further demonstrates the existence of a non-trivial interpolation that can be better than SAGA and any other minibatch version of SAGA. Moreover, we changed the distribution from which $\mA$ and $y$ are sampled to a uniform distribution where we set $\lambda = 401/(\sqrt{n^3})$ for the experiment with $n=100000$ and $d=100$. Such a setup also demonstrated the existence of a non trivial $q^* = 0.847$ as an optimal interpolation factor. In all synthetic experiments, one can observe that SAGD is always at least as good as SAGA with the former enjoying  a linear speedup in $\tau$ with a parallel implementation.

\begin{figure*}[t]
\begin{subfigure}[ht]{0.99\linewidth}
\includegraphics[width=0.33\textwidth]{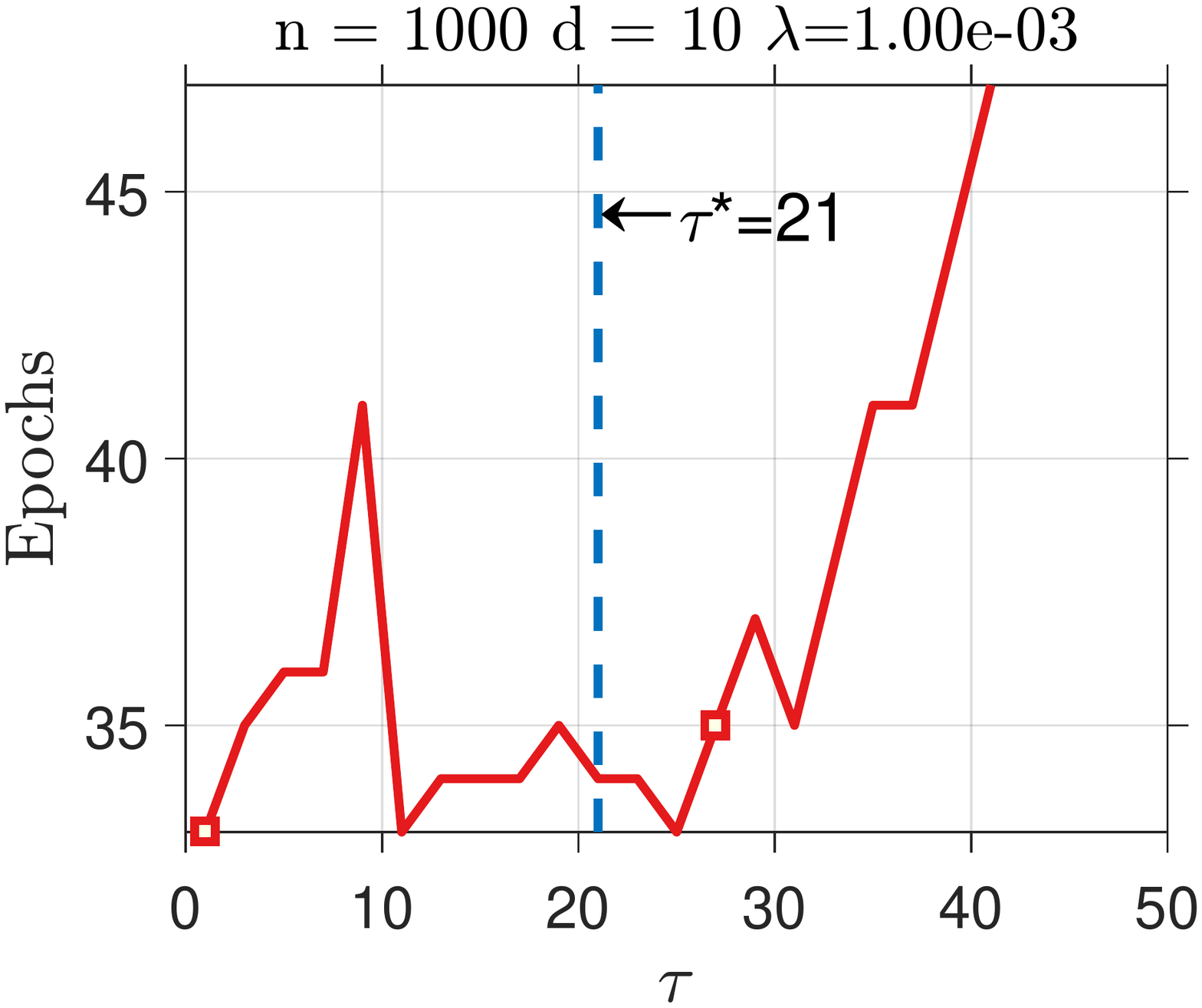}
\includegraphics[width=0.33\textwidth]{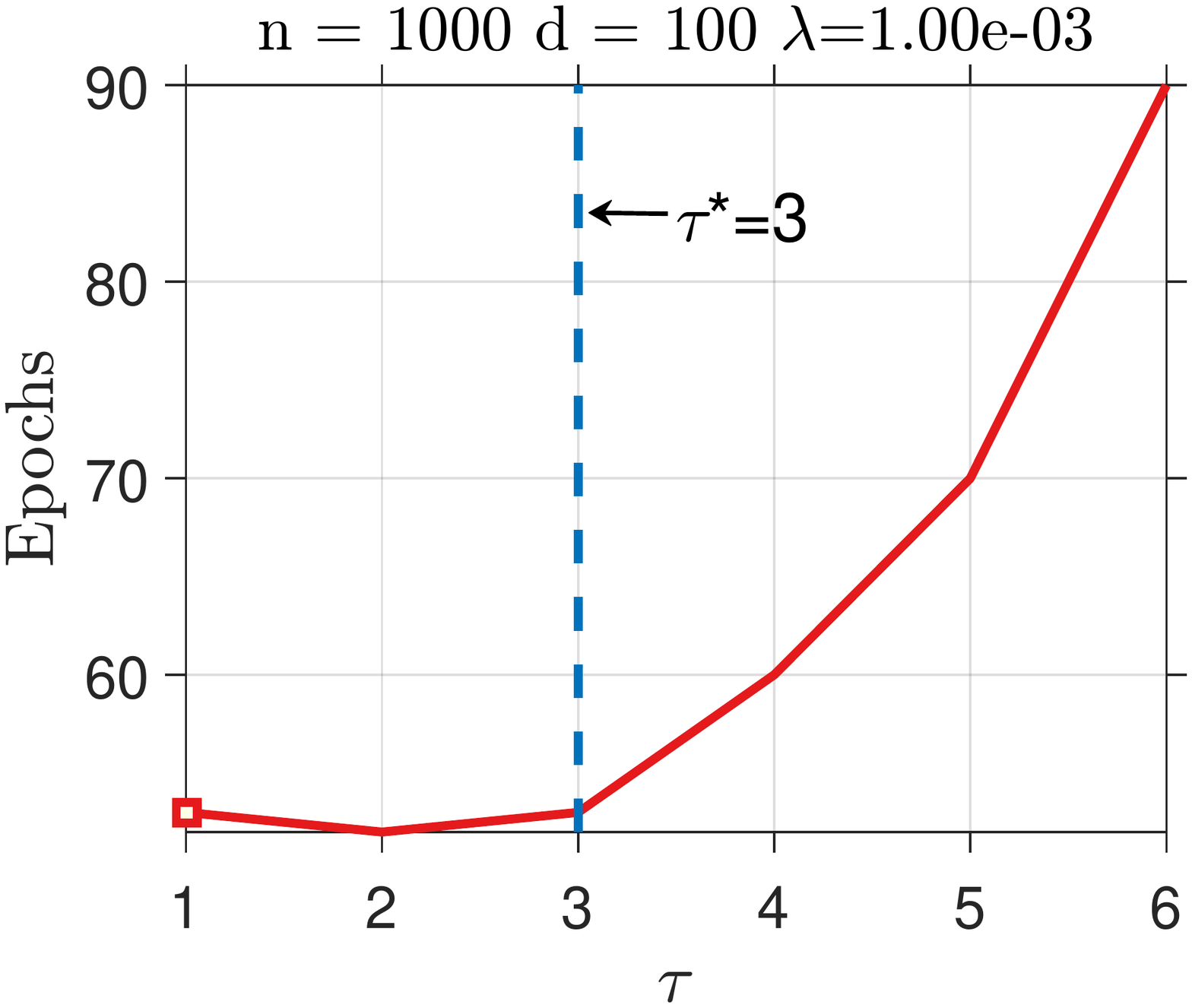}
\includegraphics[width=0.33\textwidth]{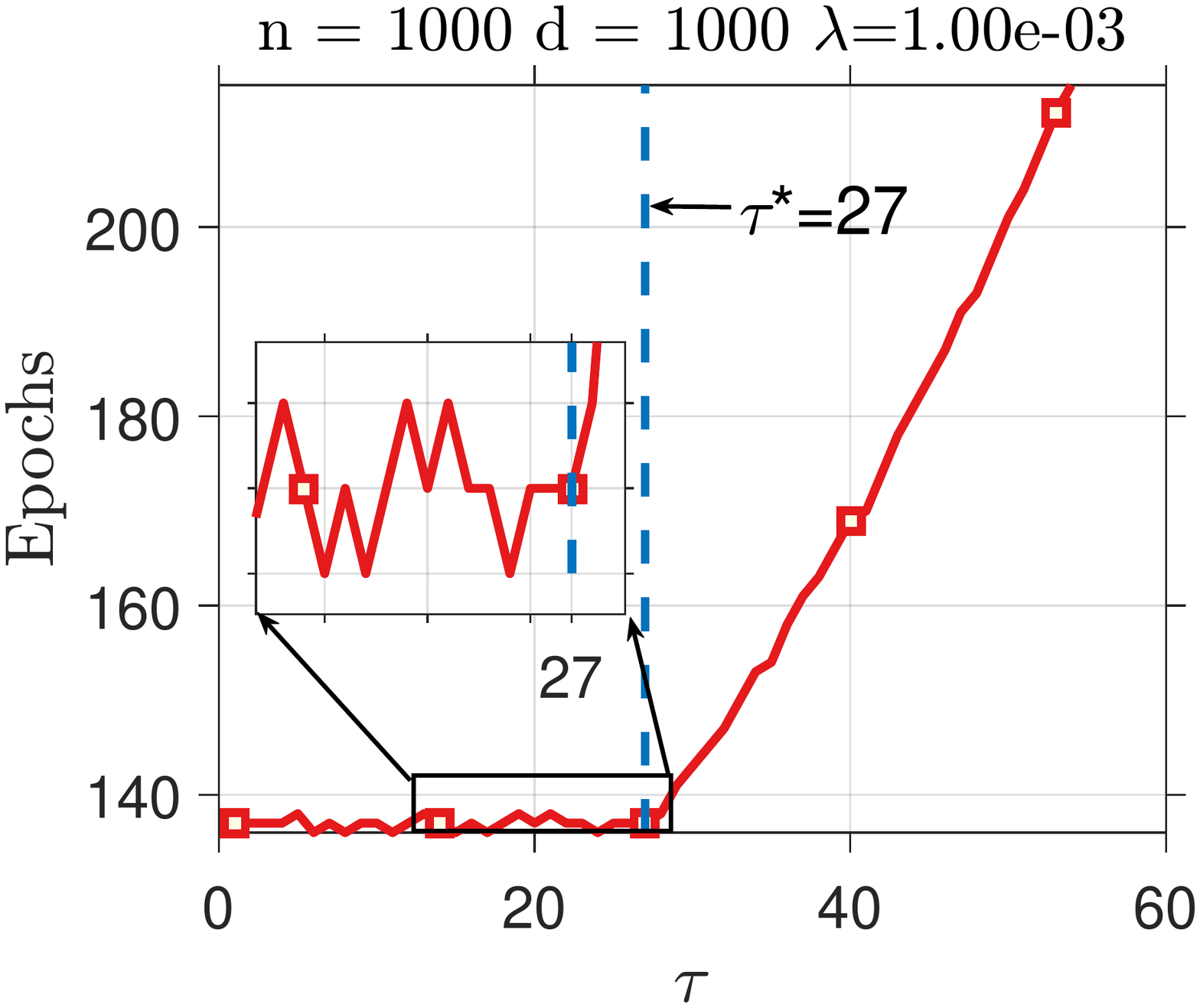}
\vspace{-2.5cm}
\end{subfigure}
\begin{subfigure}[ht]{0.99\linewidth}
\includegraphics[width=0.33\textwidth]{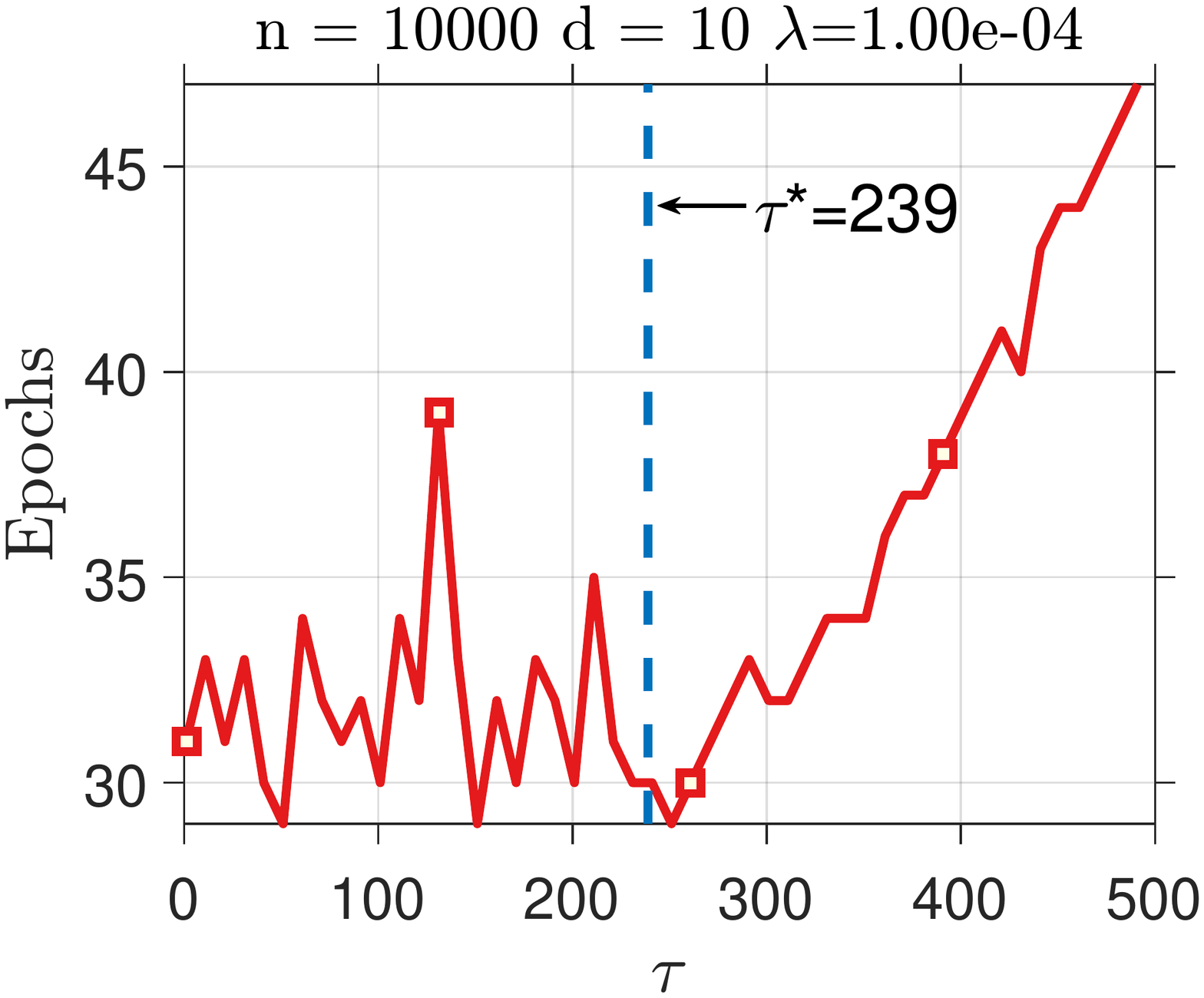}
\includegraphics[width=0.33\textwidth]{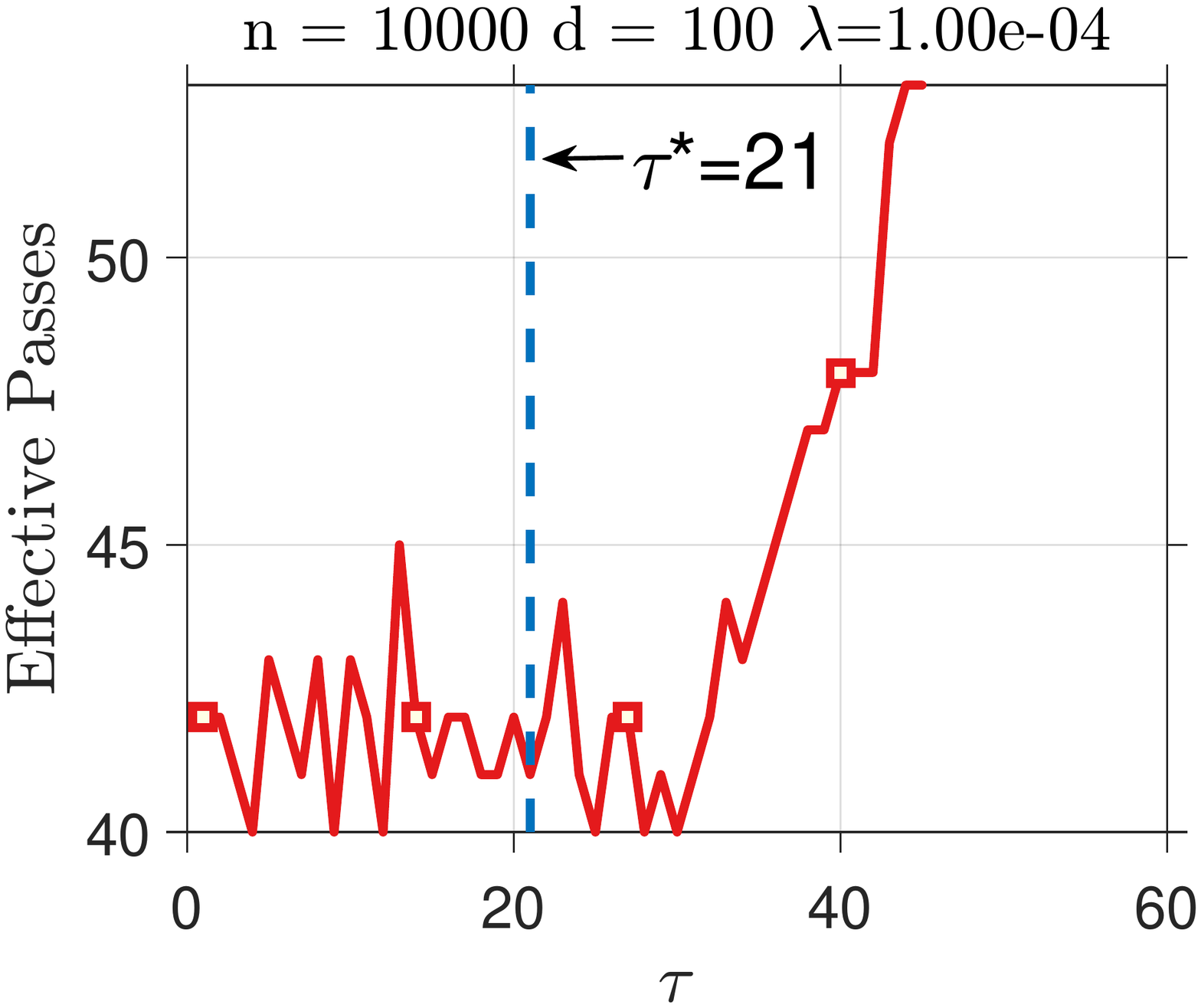}
\includegraphics[width=0.33\textwidth]{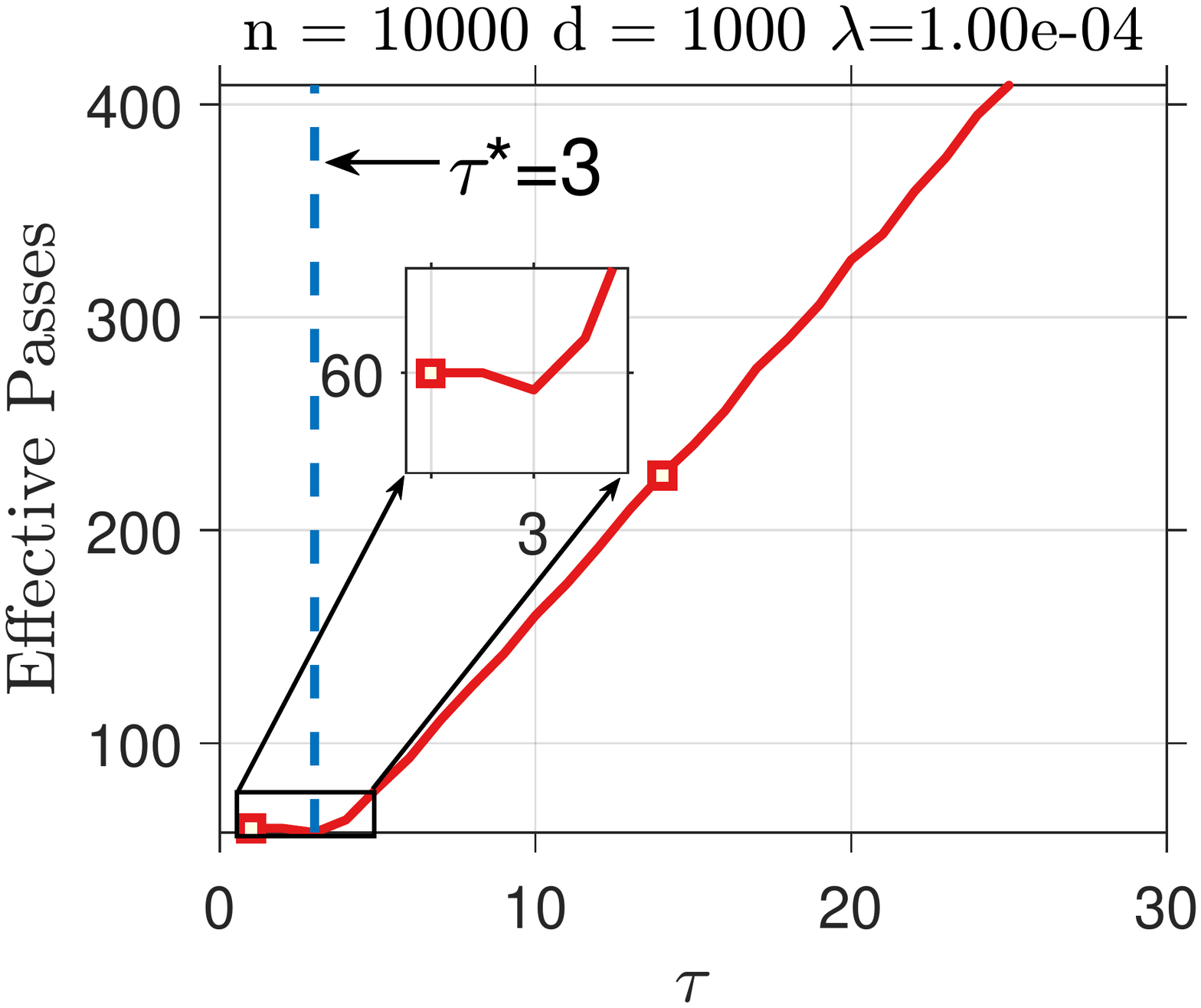}
\vspace{-1.75cm}
\end{subfigure}
\caption{Shows the number of effective data passes required for a given ($q^*,\tau$) pair to achieve $\epsilon$ accuracy. The experiments are conducted on synthetic datasets of varying $n$ and $d$. As it is noted from the plots that the theoretical $\tau^*$ matches very closely the best $\tau$  observed in practice.}
\label{sensitivity_of_tau}
\vspace{-0.5cm}
\end{figure*}

\begin{figure*}[h!]
\vspace{-0.75cm}
\begin{subfigure}[ht]{0.99\linewidth}
\includegraphics[width=0.33\textwidth]{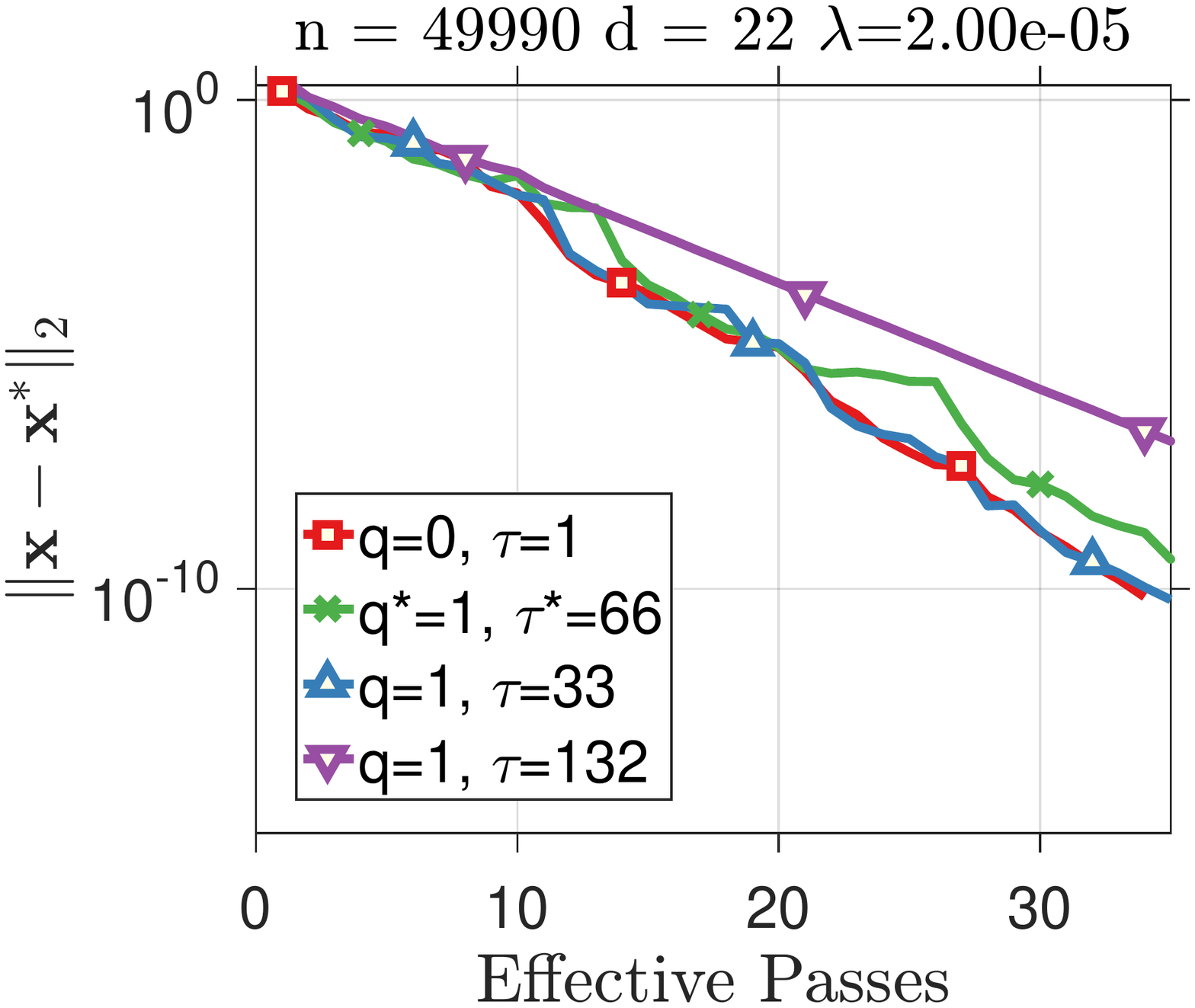}
\includegraphics[width=0.33\textwidth]{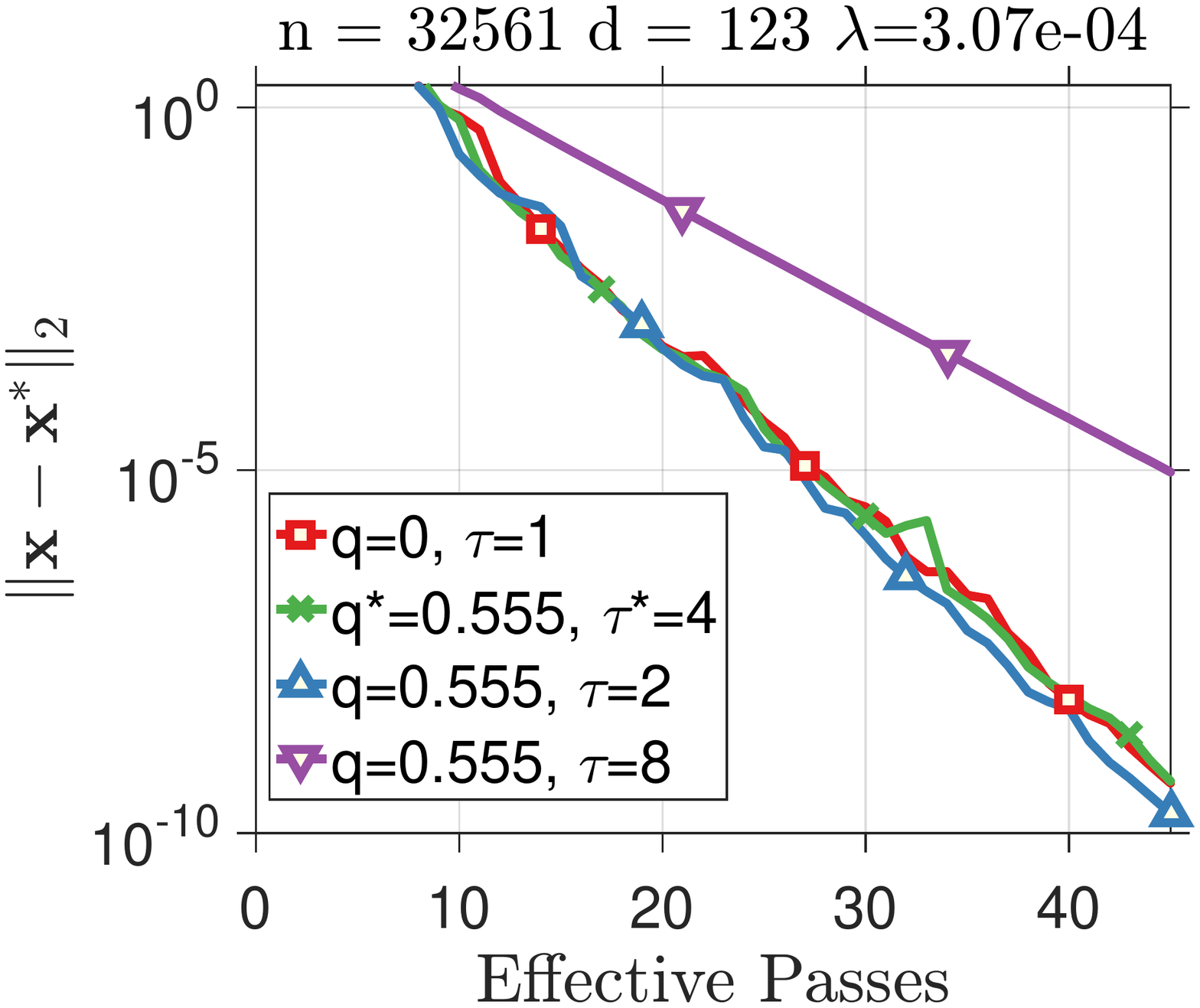}
\includegraphics[width=0.33\textwidth]{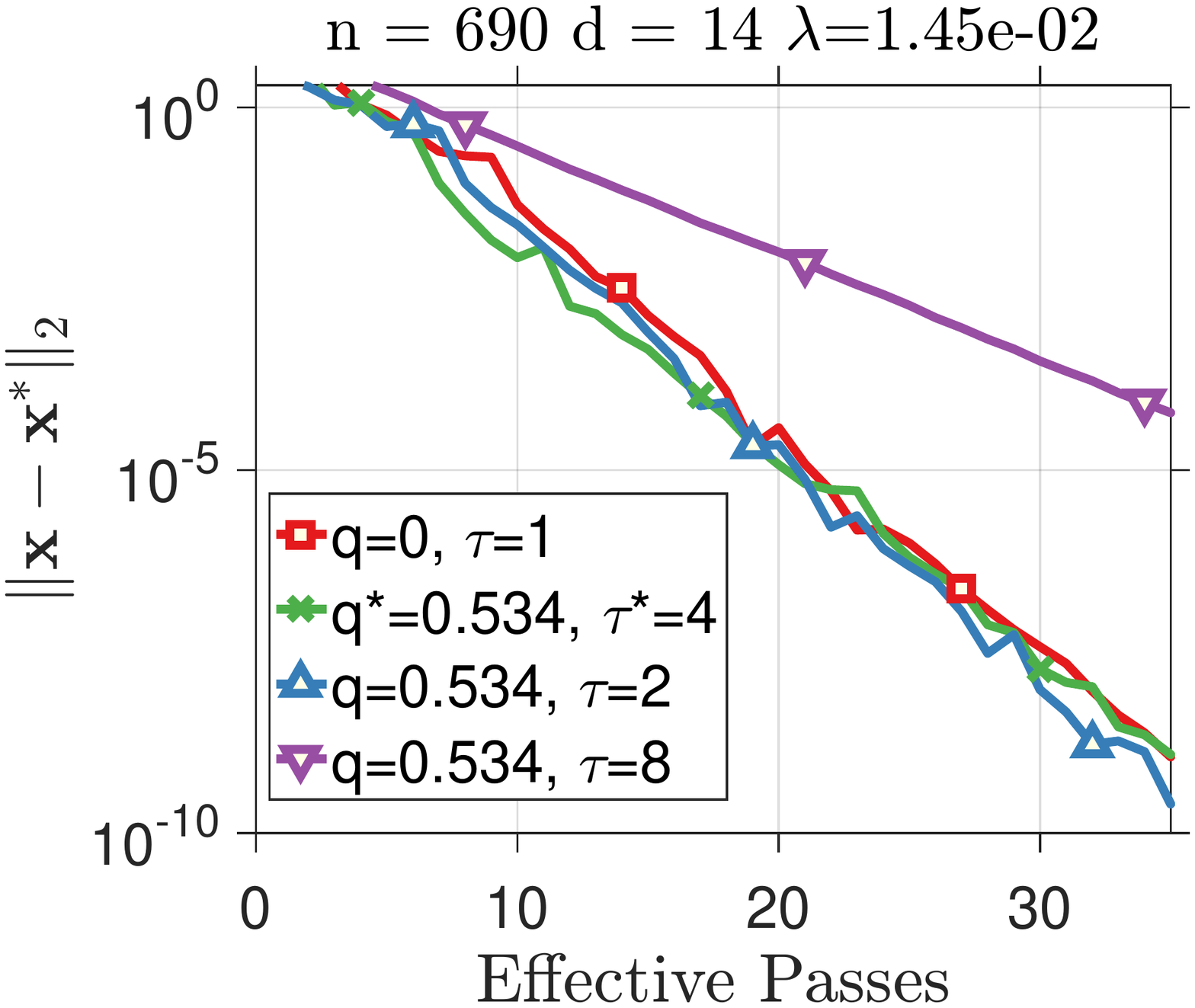}
\vspace{-1.5cm}
\end{subfigure}
\caption{Comparison of SAGD  with ($q^*,\tau^*$), ($q^*,\frac{\tau^*}{2}$) and ($q^*,2\tau^*$) and SAGA (which is SAGD with ($0,1$)). We use three  LIBSVM datasets: \texttt{ijcnn1}, \texttt{a9a} and \texttt{australian} (left to right).}
\vspace{-0.5cm}
\label{real_datasets_comparisons}
\end{figure*}

\textbf{Sensitivity to Optimal} $\tau^*$. In this experiment, we demonstrate that the theoretical results of the optimal minibatch size $\tau^*$ for a given $q^*$ match the one observed in practice. We conduct experiments on six different synthetic datasets that vary in both the number of samples $n$ and the problem dimension $d$. In all experiments, we compare the number of effective data passes required to achieve $\|x - x^* \|_2 \leq 10^{-10}$ accuracy over multiple choices of $\tau$. First as can be noted from Figure~\ref{sensitivity_of_tau}, increasing the minibatch size beyond a certain value will linearly increase the number of required data passes to achieve the same accuracy. Moreover, the theoretical optimal minibatch size $\tau^*$ is very close to the best $\tau$ observed in practice. In the few cases where the theoretical $\tau^*$ does not precisely match the best $\tau$ observed in practice, the difference in the number of required data passes to achieve $\epsilon$ accuracy is very small. For instance, for the experiment ($n=1000,d=10$) the $\tau^*$ minibatch required 34 data passes for achieving $\epsilon$ accuracy as compared to the 33 passes for the best $\tau$ in practice. The difference is also 1 data pass between $\tau^*$ and the best $\tau$ is for the experiments ($n=1000,d=100,1000$). As for the experiments ($n=10000,d=10,100$), the difference is 2 data passes. Lastly for (n=$10000$,$d=1000$), the theory precisely predicts $\tau^*=3$ in practice.

\textbf{Real Datasets.} In this subsection, we conduct similar types of experiments as the one described in the synthetic datasets but on real data. Three real datasets used in the experiments are from the LIBSVM database \cite{chang2011libsvm} namely, \texttt{ijcnn1}, \texttt{a9a} and \texttt{australian}. We compare SAGA against SAGD with ($q^*,\tau^*$) and two non optimal minibatch sizes ($q^*,\tau^*/2$) and ($q^*,2\tau^*$). On the \texttt{a9a} and \texttt{australian} datasets, we set $\lambda = \frac{10}{n}$ to showcase an example of the non trivial $q=0.555$ and $q=0.534$, respectively. As can be noted again from Figure~\ref{real_datasets_comparisons}, a direct linear speedup of SAGD can be achieved with a parallel implementation on real datasets too.

%% file: sections/conclusion.tex
\section{Conclusion}
We have demonstrated that the optimal interpolation between SAGA and gradient descent is not trivially SAGA. We thus presented formulas for the optimal probabilistic interpolation factor $q$ between SAGA and the optimal minibatch SAGA of size $\tau$. We demonstrated experimentally that the proposed SAGD is consistently better than SAGA on both synthetic and real dataset. We have also shown that the theory predicts what is observed in practice for optimal batch size $\tau^*$ of minibatch SAGA. Interestingly, experiments strongly suggest that a linear speedup is achieved with a parallel implementation.

%% file: sections/appendix.tex
\newpage
\onecolumn
\appendix

\section{Proof of Theorem~\ref{theo:t-n} }
\label{sec:newesttheo}

The proof will follow by simply computing the constants in Theorem~\ref{main_theorem} and then analyzing the resulting complexity.

When $\tau =n$  we have that 
\[\theta = \frac{n}{q(\tau-1)+1} = \frac{n}{q(n-1)+1}\] and consequently 
\begin{align}
\rho \overset{\eqref{eq:rhointheo}}{ =} \begin{cases}
\frac{n(1-q)}{(q(n-1)+1)^2}    & 0\leq  q \leq  \frac{1}{(n-1)^2}\\
\frac{n(1-q)}{(q(n-1)+1)^2}+ n \left(\frac{n^2q}{(q(n-1)+1)^2}-1\right)  \qquad \qquad & \qquad q \geq  \frac{1}{(n-1)^2}
\end{cases}\label{eq:rhotau-n1}
\end{align}
Furthermore from~\eqref{eq:cL1new} and $\tau =n$ we have
\begin{eqnarray}
     \cL_1 & = &  \frac{(1-q) L_{\max}+qn^2  \overline{L}}{(q(n-1)+1)^2}. \label{eq:cL1tau-n1}
\end{eqnarray}
Inserting~\eqref{eq:rhotau-n1} and~\eqref{eq:cL1tau-n1} into~\eqref{alpha_value} gives a stepsize of
\begin{equation}
    \alpha =  \min \left\{\frac{1}{4} \frac{(q(n-1)+1)^2}{(1-q) L_{\max}+qn^2  \overline{L}}, \frac{n(q(n-1)+1)}{\mu n^2 +4 \rho L_{\max}(q(n-1)+1)}  \right\} \label{eq:alpha1}
\end{equation}
and a total complexity~\eqref{total_compleixty_SAGD} of
\begin{eqnarray}
    \Omega^{q,\tau} &=& \left(q(n-1)+1\right)  \max \left\{\frac{4 }{\mu} \frac{(1-q) L_{\max}+qn^2  \overline{L}}{(q(n-1)+1)^2}, \frac{n}{q(n-1)+1} + \frac{4 \rho L_{\max} }{\mu n} \right\}  \log \left(\frac{1}{\epsilon}\right) \nonumber \\
    &= &  \max \left\{\underbrace{\frac{4 }{\mu} \frac{(1-q) L_{\max}+qn^2  \overline{L}}{q(n-1)+1}}_{X(q)}, \underbrace{n + \frac{4 \left(q(n-1)+1\right)  \rho L_{\max} }{\mu n}}_{Y(q)} \right\}  \log \left(\frac{1}{\epsilon}\right).
    \label{eq:total_compleixty_tau-n1}
\end{eqnarray}

Now we separate the analysis into two cases depending on the condition number $  \frac{4\overline{L}}{\mu}.$

Consider the \emph{well conditioned} case where $\frac{4\overline{L}}{\mu} \leq n-1$ and let $q= \frac{1}{(n-1)^2}$. From~\eqref{eq:rhotau-n1} we have that 
\[\rho = \frac{n(1-q)}{(q(n-1)+1)^2}  =  \frac{n(1-\frac{1}{(n-1)^2})}{(\frac{1}{(n-1)^2}(n-1)+1)^2} =
\frac{n\frac{ (n-1)^2-1}{(n-1)^2}}{(\frac{n}{n-1})^2} =
n-2, \] 
which  together with  $q= \frac{1}{(n-1)^2}$ when plugged into~\eqref{eq:alpha1} gives~\eqref{eq:alphawell}, and when plugged into~\eqref{eq:total_compleixty_tau-n1} gives
\begin{eqnarray}
    \Omega^{q,\tau} 
    &= &   \max{ \left\{ \frac{n-2}{n-1}\frac{4L_{max}}{\mu} +  \frac{n}{n-1} \frac{4 \overline{L}}{\mu}, n+ \frac{n-2}{n-1} \frac{4L_{\max}}{\mu} \right\} } \log \left(\frac{1}{\epsilon}\right) \nonumber\\
    & \overset{\frac{4\overline{L}}{\mu}\leq n-1}{=}& 
    n+ \frac{n-2}{n-1} \frac{4L_{\max}}{\mu}\log \left(\frac{1}{\epsilon}\right).
    \nonumber
\end{eqnarray}

Now consider the \emph{badly conditioned} case where $\frac{4\overline{L}}{\mu} \geq n-1$ and let
$q= \frac{\mu}{4n \overline{L}}$. It follows that $q\leq \frac{1}{(n-1)^2}$ since
\[  q = \frac{\mu}{4n \overline{L}} \overset{\frac{4\overline{L}}{\mu}\geq n-1}{\leq} \frac{1}{n(n-1)} \leq \frac{1}{(n-1)^2}. \]
Consequently 
\[\rho = \frac{n(1-q)}{(q(n-1)+1)^2} =  \frac{n(1-\frac{\mu}{4n \overline{L}})}{(\frac{\mu}{4n \overline{L}}(n-1)+1)^2}  = 
 \frac{4n^2 \overline{L}(4n \overline{L}-\mu)}{(\mu(n-1)+4n \overline{L})^2}. 
\]
Plugging in the above and  $q= \frac{\mu}{4n \overline{L}}$  into  $X(q)$ defined in~\eqref{eq:total_compleixty_tau-n1} gives
\begin{eqnarray}
  X(q) &= &  \frac{4 }{\mu} \frac{(1-\frac{\mu}{4n \overline{L}}) L_{\max}+\frac{\mu}{4n \overline{L}}n^2  \overline{L}}{\frac{\mu}{4n \overline{L}}(n-1)+1} \nonumber \\
  &= & \frac{4L_{\text{max}}}{\mu} \frac{4n\overline{L} - \mu}{n(4\overline{L} + \mu) - \mu} + \frac{4\overline{L}}{\mu} \frac{n^2 \mu}{n(4\overline{L}+\mu)-\mu}
\end{eqnarray}
and plugging into $Y(q)$ gives
\begin{eqnarray}
  Y(q) &= &   n + \frac{4 \left(\frac{\mu}{4n \overline{L}}(n-1)+1\right)  \rho L_{\max} }{\mu n}  \nonumber \\
  &= &   n + \frac{4L_{\text{max}}}{\mu} \frac{4n\overline{L} - \mu}{n(4\overline{L}+\mu)-\mu}.
\end{eqnarray}
Furthermore $Y(q) \geq X(q)$ which follows by comparing the two and noting that
\begin{eqnarray}
    n &\geq &  \frac{4\overline{L}}{\mu} \frac{n^2 \mu}{n(4\overline{L}+\mu)-\mu} \nonumber \quad \Leftrightarrow \\
     \frac{n(4\overline{L}+\mu)-\mu}{n \mu} &\geq &  \frac{4\overline{L}}{\mu}  \nonumber \quad \Leftrightarrow  \\
        \frac{4\overline{L}}{\mu} +  \frac{\mu(n-1)}{n \mu} &\geq &  \frac{4\overline{L}}{\mu}.  \nonumber \quad 
\end{eqnarray} 
Consequently, the total complexity~\eqref{eq:total_compleixty_tau-n1} is dominated by $Y(q)$  and given by
\begin{eqnarray}
    \Omega^{q,\tau} 
    &= &  \left(n + \frac{4L_{\text{max}}}{\mu} \frac{4n\overline{L} - \mu}{n(4\overline{L}+\mu)-\mu} \right)\log \left(\frac{1}{\epsilon}\right).
    \nonumber
\end{eqnarray}
Taking into consideration that $Y(q) \geq X(q)$ and that the stepsize is given by 
\begin{eqnarray*}
\alpha &=&  \frac{q(n-1)+1}{\mu}\min \left\{\frac{1}{X(q)}, \frac{1}{Y(q)} \right\} \\
&=&  \frac{n(4\overline{L}+\mu)-\mu}{4n \overline{L}}  \frac{1}{Y(q)} \\
&= &
  \frac{1}{4n \overline{L}} \frac{(n(4\overline{L}+\mu)-\mu)^2}{\mu n(n(4\overline{L}+\mu)-\mu)+4L_{\max}(4n\overline{L} - \mu)}  \\
  &= &   \frac{1}{4 \overline{L}} \frac{\left((4\overline{L}+\mu)-\frac{\mu}{n}\right)^2}{\mu n((4\overline{L}+\frac{\mu}{n})-\mu)+4L_{\max}(4\overline{L} - \frac{\mu}{n})}
\end{eqnarray*}
Finally calculating the constant in Lyaponov constant in~\eqref{eq:lyapgen}
we have
\[  \frac{\alpha}{2 \cL_2} \overset{\eqref{eq:cL2calc}}{=}  =  \frac{\theta}{n L_{\max}}\frac{\alpha}{2}  =\frac{1}{ (q(n-1)+1)L_{\max}}\frac{\alpha}{2}. \qed \]

\section{Proof of Theorem~\ref{main_theorem} }

Our proof of Theorem~\ref{main_theorem} is based on the more  the general convergence Theorem 3.6 in \cite{2018arXiv180502632G}. In the setting~\cite{2018arXiv180502632G}, the authors allow for different Euclidean norms based on a positive definite matrix $\mW \in \mathbb{R}^{d\times d}$. Here we only consider the standard Euclidean norm (thus $\mW = \mI$ ).
Let
\begin{equation}
    \cG \quad \eqdef \quad \{C \; : \; C \subset \{1,\ldots, n \}, \; |C| = \tau \},
    \end{equation}
and consider the sampling given in~\eqref{eq:one-or-tau}.
That is we either sample a coordinate $j$  with probability $\left.(1-q)\right/n$ or we 
sample a minibatch $C \in \cG$ with probability $\left. q\right/\binom{n}{\tau}$. Note that there are exactly $\binom{n}{\tau}$ sets in $\cG$. Furthermore, if we fix a coordinate $j$, there are exactly  $c_1 \eqdef \binom{n-1}{\tau -1}$ sets in $\cG$ that contain $j$, in other words
\begin{equation}\label{eq:c1def}
    c_1 \eqdef \sum_{C \in \mathcal{G} \,: \, j \in C}  1 \;=\; \binom{n-1}{\tau-1},
\end{equation}
independent on our choice of $j$.
We also define a projection matrix $\Pi_S$ based on this sampling in~\eqref{eq:one-or-tau}.
\begin{equation}
    \begin{aligned}
    \mS = 
    \begin{cases}
    \mI_C & \text{with prob} ~~ \frac{q}{\binom{n}{\tau}}, \quad C \in \cG \\
    e_j   & \text{with prob} ~~ \frac{1-q}{n}
    \end{cases}
    \end{aligned}
        \qquad \quad
    \begin{aligned}
    \Pi_S = 
    \begin{cases}
    \mI_C \mI_C^\top & \text{with prob} ~~ \frac{q}{{\binom{n}{\tau}}} , \quad C \in \cG \\
    e_j e_j^\top & \text{with prob} ~~ \frac{1-q}{n}
    \end{cases}
    \end{aligned}\label{eq:one-or-tau}
\end{equation}
where $\mI_C$ is a column submatrix of the identity matrix and $e_j$ is the all zero vector except for the $\text{j}^{\text{th}}$ coordinate with a value of one. In this setting, we can now fully analyze our proposed algorithm as a direct implication of Theorem 3.6 in \cite{2018arXiv180502632G}  whose statement we repeat here for ease of reference.

\begin{theorem}[Convergence of JacSketch for General Sketches] \label{theo:convgen}  Let $f$ be $\mu$--strongly convex and each $f_i$ be convex and $L_i$--smooth. Let  $\mS$ be defined by~\eqref{eq:one-or-tau} and let $\theta_{\mS}$ be the associated bias-correcting random variable that satisfies
\begin{equation}\label{eq:thetadef}
    \mathbb{E} \left[\theta_S \Pi_S e \right]e  = e,
\end{equation}
where $e \in \mathbb{R}^n$ is the vector of all ones.
Let $\cL_1$ and $\cL_2$ be the \emph{expected smoothness} constants that satisfy
\begin{eqnarray}\label{eq:expsmoothcL1}
    \frac{1}{n^2} \mathbb{E} \left[\|\mR \theta \Pi_S e \|_2^2\right] &\leq& 2 \cL_1 \left(f(x) - f(x^*)\right) \\ \label{eq:expsmoothcL2}
     \mathbb{E} \left[\|\mR \Pi_S \|_2^2\right] &\leq& 2 \cL_2 \left(f(x) - f(x^*)\right),
\end{eqnarray}
where  $\mR = \nabla \mF(x) - \nabla \mF(x^*)$. Let 
 \begin{equation} \label{eq:kappdef} \kappa  = \lambda_{\text{min}} \left(\mathbb{E}[\Pi_S]\right), 
 \end{equation}
 be the stochastic condition number and assume that $\kappa >0.$
  Let the sketch residual $\rho$ be defined as 
  \begin{equation}\label{eq:rhodef}
     \rho \eqdef \lambda_{\max}\left( \mathbb{E} \left[\theta_{\mS}^2\Pi_S e e^\top \Pi_S \right]- ee^\top \right) \geq 0.
  \end{equation}
  
Choose any $x^0\in \mathbb{R}^d$ and $\mJ^0\in \mathbb{R}^{d\times n}$. Let $\{x^k,\mJ^k\}_{k\geq 0}$ be the random iterates produced by Jacsketch Algorithm in \cite{2018arXiv180502632G}. Consider the Lyapunov function
\begin{equation}\label{eq:lyapgen}
\Psi^k \eqdef \|x^k-x^*\|_2^2 +  \frac{\alpha}{2 \cL_2} \|\mJ^{k} -\nabla F(x^*)\|^2.
\end{equation}
 If the stepsize satisfies
	\begin{equation}
	0 \leq \alpha \leq \min\left\{ \frac{1}{4 \cL_1 }, \, \frac{\kappa}{4 \cL_2 \rho/n^2   +\mu  }\right\},
	\label{eq:alphaboundXX}
	\end{equation}
	then 
	\begin{equation} \label{eq:conv2XX}
	\E{\Psi^{k}} \quad \leq \quad (1-\mu \alpha)^k \cdot \Psi^0,
	\end{equation}
If we choose $\alpha$ to be equal to the upper bound in \eqref{eq:alphaboundXX}, then
	\begin{equation}\label{eq:itercomplexgen}
	k\geq \max \left\{ \frac{4\cL_1}{\mu}, \; \frac{1}{\kappa} + \frac{4\rho \cL_2}{\kappa \mu n^2}  \right\} \log\left (\frac{1}{\epsilon} \right) \quad \Rightarrow \quad \E{\Psi^{k}}  \leq \epsilon \Psi^0.\end{equation}
\end{theorem}

Now all we need to do to apply this theorem is to compute the constants $\theta_S,\kappa,\cL_1,\cL_2,\rho$ under this new sampling $\mS$.

 In our upcoming computations, we will often make use of the following type of equality
\begin{equation}\label{eq:doublecount}
    \sum_{C \in \cG} \sum_{i \in C} a_{i,C} =    \sum_{i=1}^n \sum_{C \in \cG \; : \; i \in C} a_{i,C},
\end{equation}
which holds for any $a_{i,C} \in \mathbb{R}$ and can be proved using a \emph{double counting} argument.

\newpage
\underline{\textbf{The Bias Correcting Term:}} $\theta_S$

It turns out that a constant bias correcting variable $\theta_S \equiv \theta$ is sufficient to guarantee~\eqref{eq:thetadef} since
\begin{eqnarray}
     \displaystyle \mathbb{E} \left[\theta_S \Pi_S e \right]e 
    &= &\sum_i^n \theta e_i e_i^\top e \frac{1-q}{n} 
    + \sum_{C \in \cG}\theta \, \mI_C \mI_C^\top e \frac{q}{\binom{n}{\tau}} \\
    &=& \theta \left(\frac{1-q}{n} \mI + \frac{q}{\binom{n}{\tau}}\mI \binom{n-1}{\tau -1} \right) e = e,\label{eq:asd89j8js}
\end{eqnarray}
where we used in the second equality that 
\begin{equation}\label{eq:sumIcIc}
    \sum_{C \in \cG} I_C I_C^\top = \binom{n-1}{\tau-1} I.
\end{equation}
From~\eqref{eq:asd89j8js} and the definition~\eqref{eq:thetadef} we have that 
\begin{equation}\label{eq:theta}
    \theta = \frac{n}{q(\tau - 1) + 1}.
    \end{equation}
\underline{\textbf{The Stochastic Condition Number:}} $\kappa$

To calculate $\kappa$ we first compute
\begin{eqnarray}
    \mathbb{E}[\Pi_S] &=& \sum_i^n \frac{1-q}{n} e_i e_i^\top +  \sum_{C \in \cG}  \mI_C \mI_C^\top \frac{q}{\binom{n}{\tau}}\nonumber \\
    &\overset{\eqref{eq:sumIcIc}}{=}& \left(\frac{1-q}{n}+ \frac{q}{\binom{n}{\tau}} \binom{n-1}{\tau-1} \right) \mI \nonumber \\
    & =&    \frac{q(\tau - 1) + 1}{n} \quad \overset{\eqref{eq:theta} }{=}\quad  \frac{1}{\theta}. \label{eq:EPiS}
\end{eqnarray}
Consequently
\begin{equation} 
    \kappa  \overset{\eqref{eq:kappdef} }{=} \lambda_{\min} \left(\mathbb{E}[\Pi_S]\right) \overset{\eqref{eq:EPiS}}{=}
    \frac{1}{\theta}. \label{eq:kappaapp}
\end{equation}

\underline{\textbf{Expected Smoothness Constant of the Stochastic Gradient:}} $\cL_1$

Here we deduce the formula~\eqref{eq:cL1new} from the definition of expected smoothness in~\eqref{eq:expsmoothcL1}.

Since each $f_i(x)$ is $L_i$--smooth and convex, then for $C \subseteq [n]$, $f_C(x) = \frac{1}{|C|} \sum_{i \in C} f_i(x)$ is $L_C$ smooth where $L_C = \frac{1}{|C|}\sum_{i \in C} L_i$. 
Therefore, the following is equivalent to the smoothness assumption (See 2.1.7 in \cite{NesterovBook}):
\begin{align}
    \|\nabla f_C(x) - \nabla f_C(x^*) \|_2^2 \leq 2 L_C \left(f_C(x) - f_C(x^*)\right) - \langle \nabla f_C(x^*), x-x^* \rangle,
    \label{nestrovs_smoothness}
\end{align}
and by an analogous argument and summing up we have
\begin{equation}\label{eq:nabfismooth}
      \sum_{i=1}^n \|\nabla f_i(x) - \nabla f_i(x^*) \|_2^2 \leq 2 nL_{\max} \left(f(x) - f(x^*)\right).
\end{equation}
It now follows that
\begin{align*}
    &\frac{1}{n^2} \mathbb{E}\left[\|\mR \theta \Pi_S e \|_2^2\right] = \frac{\theta^2}{n^2} \mathbb{E} \left[\|\mR \Pi_S e \|_2^2\right] = \frac{\theta^2}{n^2}\sum_{C \in \cG}\|\mR \mI_C \mI_C^\top e \|_2^2 \frac{q}{\binom{n}{\tau}} + \frac{\theta^2}{n^2} \sum_i^n \|\mR e_i \|_2^2 \frac{1-q}{n} \\
    & = \frac{\theta^2}{n^2}\left(\frac{q}{\binom{n}{\tau}}\sum_{C \in \cG} \|\mR e_C\|_2^2 + \frac{1-q}{n}\sum_i^n \|\mR e_i\|_2^2\right) \\
    & \overset{|C| = \tau}{=} \frac{\theta^2}{n^2}\left(\tau^2 \frac{q}{\binom{n}{\tau}} \sum_{C \in \cG }\|\nabla f_C(x) - \nabla f_C(x^*)\|_2^2  + \frac{1-q}{n} \sum_i^n \|\nabla f_i(x) - \nabla f_i (x^*) \|_2^2\right) \\
    &\overset{\eqref{eq:nabfismooth}}{\leq} \frac{\theta^2}{n^2} \left(2 \tau^2  \frac{q}{\binom{n}{\tau}} \sum_{C \in \cG }L_C \left((f_C(x) - f_C(x^*) - \langle \nabla f_C(x^*),x-x^* \rangle\right)\right.) \left.+  2n\frac{1-q}{n}  L_{\text{max}} (f(x) - f(x^*) \right) \\
    & \overset{\eqref{eq:doublecount}}{=} \frac{\theta^2}{n^2}\left(2\tau \frac{q}{\binom{n}{\tau}}  \sum_{i=1}^n \sum_{C \in \cG \; : \; i \in C}L_C   \left((f_i(x) - f_i(x^*) - \langle \nabla f_i(x^*),x-x^* \rangle\right) \right.\left. + 2(1-q) L_{\text{max}} \left(f(x) - f(x^*)\right)\right) \\
    &\leq \frac{\theta^2}{n^2} \left(\frac{2q \tau n}{\binom{n}{\tau}} \left(f(x) - f(x^*)\right)\max_{i=1,\ldots, n}\sum_{C \in \cG \; : \; i \in C}L_C + 2(1-q)L_{\text{max}}(f(x) - f(x^*)\right) \\
    &= 2\frac{\theta^2}{n^2} \left(\frac{q\tau n}{\binom{n}{\tau}} \left(\max_{i=1,\ldots, n}\sum_{C \in \cG \; : \; i \in C}L_C\right) + (1-q) L_{\text{max}} \right) (f(x) - f(x^*)).
\end{align*}

Note that $e_C$ is the vector of all zeros except for the indices selected by the sample C have a value of one. Note that since the sampling is $\tau$-uniform therefore $|C| = \tau$; therefore, the forth equality follows from the fact that $\frac{1}{|C|} \nabla \mF(x) e_C = \nabla f_C(x)$. The second part of the first inequality follows from first order optimality conditions where $\sum_i \nabla f_i (x^*) = 0$.

It now follows from the definition~\eqref{eq:expsmoothcL1} and the above that
\begin{equation}\label{eq:cL1res}
    \cL_1 = \frac{\theta^2}{n} \frac{q \tau}{\binom{n}{\tau}}\left(\max_{i=1,\ldots, n}\sum_{C \in \cG \; : \; i \in C}L_C\right) + \frac{\theta^2 (1-q)}{n^2} L_{\text{max}}.
    \end{equation}

We can bound the max term in~\eqref{eq:cL1res}  under the assumption that $L_C = \frac{1}{\tau} \sum_{i \in C}L_i,$ as follows. Using another double counting argument we have that
\begin{eqnarray}
    \sum_{C \in \cG \; : \; i \in C}L_C & =&  
    \frac{1}{\tau}  \sum_{C \in \cG \; : \; i \in C} \;\sum_{j \in C} L_j \nonumber \\
    & \overset{\eqref{eq:doublecount}}{=}& 
     \frac{1}{\tau} \sum_{j =1}^n \;\sum_{C \in \cG \; : \; i \in C, \; j \in C}  L_j  \nonumber\\
    &= &    \frac{1}{\tau} \sum_{j \neq i}\;
    \sum_{C \in \cG \; : \; i \in C, \; j \in C} L_j+\frac{1}{\tau} 
    \sum_{C \in \cG \; : \; i \in C} L_i \nonumber\\
    & = &  \frac{1}{\tau} \sum_{j \neq i} \binom{n-2}{\tau -2} L_j+\frac{1}{\tau} 
    \binom{n-1}{\tau -1}L_i   \nonumber\\
    & =&  \frac{n}{\tau} \binom{n-2}{\tau -2} \overline{L}+\frac{1}{\tau} 
    \left(\binom{n-1}{\tau -1}-\binom{n-2}{\tau -2}\right)L_i \nonumber \\
    &= &  \binom{n-2}{\tau -2}\left( \frac{n}{\tau}  \overline{L}+\frac{1}{\tau} 
    \frac{n-\tau}{\tau-1}L_i\right), \label{eq:cLpaaa}
\end{eqnarray}
where $\overline{L} \eqdef \tfrac{1}{n} \sum_{i=1}^n L_i.$
Inserting the above into~\eqref{eq:cL1res} gives
\begin{eqnarray}
     \cL_1 &\overset{\eqref{eq:cLpaaa}}{=}& \frac{\theta^2}{n} 
     \frac{q \tau}{\binom{n}{\tau}}\binom{n-2}{\tau -2}\left( \frac{n}{\tau}  \overline{L}+\frac{1}{\tau} 
    \frac{n-\tau}{\tau-1}L_{\max}\right) + \frac{\theta^2 (1-q)}{n^2} L_{\max} \nonumber \\
    & =&  q \theta^2 
     \frac{\tau(\tau-1)}{n(n-1)}\left(  \overline{L}+
    \frac{n-\tau}{n(\tau-1)}L_{\max}\right) + \frac{\theta^2 (1-q)}{n^2} L_{\max}  \nonumber \\
    & =&   \frac{\theta^2}{n^2}\left( q \left(\frac{\tau(n-\tau)}{n-1}-1\right) +1 \right) L_{\max} + q \theta^2 
     \frac{\tau(\tau-1)}{n(n-1)}\overline{L}.\label{eq:newcL1}
\end{eqnarray}

\underline{\textbf{Expected Smoothness Constant of the Stochastic Jacobian:}} $\cL_2$

It follows from (87) in \cite{2018arXiv180502632G} that
\begin{equation}\label{eq:oaksk93k9}
    \mathbb{E} \left[\|\mR \Pi_S \|_F^2 \right] = \text{trace}\left(\mR^\top \mR \mE\left[\Pi_S \Pi_S^\top\right]\right) \leq 2n \lambda_{\max}\left(\mD_L^{\frac{1}{2}}\mE[\Pi_S]\mD_L^{\frac{1}{2}}\right) \left(f(x) - f(x^*)\right) 
\end{equation}
where $\mD_L = \text{Diag}\left(L_1,\dots,L_n\right)$. Moreover from~\eqref{eq:EPiS} we have that  $\mathbb{E}[\Pi_S] = \frac{1}{\theta} \mI$. Thus comparing~\eqref{eq:expsmoothcL1} and~\eqref{eq:oaksk93k9} we have that
\begin{equation} \label{eq:cL2calc}\cL_2 = \frac{n}{\theta} \max_i L_i = \frac{n L_{\max}}{\theta}.\end{equation}

\underline{\textbf{The Sketch Residual:}} $\rho$

To compute the sketch residual we first note that
\begin{align*}
    &\Theta \eqdef \theta^2 \mathbb{E} \left[\Pi_S e e^\top \Pi_S \right]- e e^\top \\
    &= \theta^2 \left(\sum_{C \in \cG} \frac{q}{\binom{n}{\tau}}e_C e_C^\top+ \sum_i^n \frac{1-q}{n} e_i e_i^\top\right) - ee^\top \\
    & = \theta^2 \left(\frac{q}{\binom{n}{\tau}}c_2 \left(ee^\top - \mI \right) + \frac{q}{\binom{n}{\tau}} c_1 \mI + \frac{1-q}{n} \mathbf{I}\right) -ee^\top
    \\
    & =\mI \underbrace{\theta^2 \left(\frac{1-q}{n} + \frac{qc_1}{\binom{n}{\tau}} - \frac{q c_2}{\binom{n}{\tau}}\right)}_{\alpha} + ee^\top \underbrace{\left(\frac{q c_2}{\binom{n}{\tau}} \theta^2 - 1\right)}_{\beta}.
    \end{align*}
Note that $c_1 = \binom{n}{\tau} \frac{\tau}{n} \forall \tau \ge 1$ and $c_2 = \binom{n}{\tau} \frac{\tau}{n} \frac{\tau-1}{n-1} \forall \tau >1$. Note that $c_2 = 0$ for $\tau = 1$. Since $\Theta$ is a circulant matrix with a first column being $[\alpha + \beta; \beta \mathbf{e}^{n-1}]$, where $\mathbf{e}^{n-1}$ is the vector of all ones of size $(n-1)$, a well known result for calculating eigenvalues of circulant matrices~\cite{varga1954}  states that
the eigenvalues of $\Theta$ are simply the 1-D Discrete Fourier transform of the first column, that is $\left(\text{DFT}[\alpha + \beta; \beta \mathbf{e}^{n-1}]\right) = [\alpha +n\beta, \alpha e_{n-1}]$.  Since $\Theta$ only has two distinct eigenvalues we have that
\[\rho \overset{\eqref{eq:rhodef} }{=} \lambda_{\max} \left(\Theta\right) = \max \left(\text{DFT}[\alpha + \beta; \beta \mathbf{e}^{n-1}]\right)  = \max \left(\alpha + n \beta,\alpha\right).\] 
Consequently
\begin{align*}
    \rho = 
    \begin{cases}
    \alpha + n \beta & \beta \ge 0,\\
    \alpha & \beta \leq 0.
    \end{cases}
\end{align*}
Inserting these values for $\alpha$ and $\beta$ gives
\begin{align*}
    \rho = 
    \begin{cases}
    \theta^2\left(\frac{1-q}{n} + \frac{q(c_1-c_2)}{\binom{n}{\tau}}\right) + n \left(\frac{\theta^2 q c_2}{\binom{n}{\tau}} -1\right)  & \left(\frac{\theta^2 q c_2}{\binom{n}{\tau}} -1\right) \ge 0,\\
    \theta^2\left(\frac{1-q}{n} + \frac{q(c_1-c_2)}{\binom{n}{\tau}}\right) & \left(\frac{\theta^2 q c_2}{\binom{n}{\tau}} -1\right) \leq 0.
    \end{cases}
\end{align*}
After further simplifications we have that
\begin{align}
\rho = \begin{cases}
\theta^2 \left(\frac{1-q}{n} + q \frac{\tau}{n} \frac{n-\tau}{n-1}\right) \qquad \qquad \qquad \qquad \qquad q\theta^2 \le \frac{n}{\tau}\frac{n-1}{\tau-1} \\
\theta^2 \left(\frac{1-q}{n} + q \frac{\tau}{n}\frac{n-\tau}{n-1}\right) + n \left(\theta^2 q \frac{\tau}{n} \frac{\tau-1}{n-1}-1\right)  ~ \quad q\theta^2 \geq \frac{n}{\tau}\frac{n-1}{\tau-1}
\end{cases}\label{eq:rhoconstfin}
\end{align}

By substituting all the previously computed constants~\eqref{eq:theta}, \eqref{eq:kappaapp}, \eqref{eq:cL1onlyLmax}, \eqref{eq:cL2calc} and \eqref{eq:rhoconstfin} for the sketch $\mS$ in Theorem~\ref{theo:convgen}, the proof of Theorem 2 is complete.\qed

\section{Proof of Lemma 1}
First note that when $L_i = L_{\max}$ we have that~\eqref{eq:cL1new} becomes
\begin{equation}\label{eq:cL1onlyLmax}
    \cL_1 =  \frac{ q (\tau^2-1) +1 }{(q(\tau-1)+1)^2}L_{\max}.
\end{equation}
The rest of the proof is straightforward.

\section{Proof of Lemma 2}
Using~\eqref{eq:cL1onlyLmax} we have 
\begin{align*}
    g_1 &= \frac{4 \cL_1}{\mu}\left(q(\tau-1)+1\right) 
    \overset{\eqref{eq:cL1onlyLmax}}{=} \frac{4 L_{\max}}{\mu} \left(\frac{q\tau^2}{(q(\tau-1)+1)^2} + \frac{1-q}{(q(\tau-1)+1)^2}\right) \left(q(\tau-1) + 1\right) \\
     & = \frac{4L_{\max}}{\mu} \frac{q(\tau^2-1)+1}{q(\tau-1) + 1},
\end{align*}
and thus the derivative with respect to $q$ is:
\begin{align*}
    \frac{\partial g_1}{\partial q} &= \frac{4L_{\max}}{\mu} \left(\frac{\tau^2 - 1}{q(\tau-1) + 1}- \frac{(\tau-1 )(1+q(\tau^2-1))}{\left(1+ q(\tau-1)\right)^2} \right) \\
    &= \frac{4L_{\max}}{\mu} \left(\frac{(\tau^2-1)(q(\tau-1)+1) - (\tau-1)(q(\tau^2-1)+1)}{(q(\tau-1)+1)^2}\right)
\end{align*}
Then it is sufficient to show that if $(\tau^2-1)(q(\tau-1)+1) - (\tau-1)(q(\tau^2-1)+1) \ge 0 $ then $g_1$ is monotonically. Note that

\begin{align*}
    &(\tau^2-1)(q(\tau-1)+1) - (\tau-1)(q(\tau^2-1)+1) \\
    & = \tau^2 -1 - (\tau-1)  = \tau(\tau - 1) \ge 0 \forall ~~~~ \tau \ge 1 \\
\end{align*}
Therefore $g_1$ is monotonically increasing in $q$.

\section{Proof of Lemma 3}

We will denote $f \sim_{\mathcal{D}} g$ to indicate that $f$ and $g$ have the same monotonicity on the domain $\mathcal{D}$, i.e., (if $\frac{\partial f}{\partial q}  \ge 0$ then $\frac{\partial g}{\partial q} \ge 0$ on the domain $\mathcal{D}$ and vice verse.)
Moreover, $f\simeq g$ indicates that both functions have the same concavity direction, i.e, (if $\frac{\partial^2 f}{\partial q^2} \ge 0$ then $\frac{\partial^2 g}{\partial q^2} \ge 0$ on the domain $\mathcal{D}$ and vice verse). For instance, consider the domain of $q$ to be the real numbers, that is $\mathcal{D} = \mathbb{R}$ and note that:
\begin{eqnarray*}
g_2 &=& \left(\theta + \frac{4 \rho L_{\max}}{\mu n} \right) \left(q (\tau - 1) +1 \right)\\
&\overset{\eqref{eq:thetadef} }{=}& \left( n + \frac{ 4 \rho L_{\max}}{\mu n} \left(q (\tau-1) +1 \right)   \right)  \sim  \rho \left( q(\tau-1)+1 \right).
\end{eqnarray*}
First we will start by showing that $g_2$ is a decreasing function in $q$ in the interval $q \in [0, q_{-}] \bigcup [q_{+}, 1]$. Therefore, consider the domain of $q$ where $\mathcal{D} =  [0, q_{-}] \bigcup [q_{+}, 1]$. Then $\rho = \theta^2 \left(\frac{1-q}{n} + q \frac{\tau}{n} \frac{n-\tau}{n-1}\right)$ and thus:

\begin{align*}
g_2 &\sim \frac{n^2}{(q(\tau-1)+1)^2} \left(\frac{1-q}{n} + q \frac{\tau}{n} \frac{n-\tau}{n-1} \right)\left(q(\tau-1)+1\right)  = \frac{n^2}{(q(\tau-1)+1)} \left(\frac{1-q}{n} + q \frac{\tau}{n} \frac{n-\tau}{n-1} \right) \\
&\sim \frac{1-q}{q(\tau - 1) +1} + \frac{ q \tau (n-\tau)}{(n-1)(q (\tau - 1) +1)} = \frac{n-1-q n+q+q\tau n - q \tau^2}{(n-1)(q (\tau-1)+1)} \\
&\sim \frac{q(\tau n - \tau^2 + 1 - n) + n -1}{q (\tau-1)+1} =
\frac{ q (\tau-1) (n-\tau-1) + n-1}{q(\tau - 1) + 1} \\
&=\frac{(n-\tau - 1)( q (\tau-1) +1) + \tau } {q (\tau - 1) + 1} = n- \tau - 1 + \frac{\tau}{q (\tau - 1) + 1} 
\end{align*}
which is decreasing in $q$, i.e., $\frac{\partial}{\partial q} \left( n- \tau - 1 + \frac{\tau}{q (\tau - 1) + 1}\right) \leq 0$). Since $g_2 \sim  \left( n- \tau - 1 + \frac{\tau}{q (\tau - 1) + 1}\right)$, then $\frac{\partial g_2}{\partial q} \leq 0$; thus, $g_2$ is a decreasing function in $q$ in the domain $\mathcal{D} = [0, q_{-}] \bigcup [q_{+}, 1]$.

Now, with a similar argument but considering the second derivative with the domain of $q$ as $\mathcal{D} =  [q_{-}, q_{+}]$. Note that in this domain $\rho =\theta^2 \left(\frac{1-q}{n} + q \frac{\tau}{n} \frac{n-\tau}{n-1}\right) + n \left( \theta^2 q \frac{\tau (\tau-1)}{n (n-1)} -1 \right)$
\begin{align*}
g_2 & = \frac{n (1-q)}{q (\tau - 1) +1} + \frac{n q \tau (n-\tau)}{(n-1)(q(\tau - 1)+1)} 
+ n (q(\tau-1)+1) \left(\frac{n q \tau (\tau-1)}{(n-1)(q (\tau - 1)+1)^2} - 1 \right) \\
&\simeq \frac{n (1-q)}{q( \tau - 1)+1} + \frac{n q \tau (n-\tau)}{(n-1)(q(\tau - 1)+1)} + n^2  \frac{q \tau (\tau-1)}{(n-1)(q (\tau - 1)+1)} \\
&= \frac{ n^2 - n - n^2 q + n q + n^2 q \tau - n q \tau^2 + n^2 q \tau^2 - n^2 q \tau }{(n-1)(q(\tau-1)+1)} \\
&\simeq \frac{ q ( -n^2 + n - n \tau^2 + n^2 \tau^2) + n^2 - n}{q(\tau-1)+1} = \frac{q (\tau^2-1)(n^2-n) + n^2 - n}{q(\tau-1)+1} \\
&\simeq \frac{q (\tau^2-1) + 1}{q(\tau-1)+1} = \frac{ (\tau+1)(q(\tau-1)+1) -\tau}{q(\tau-1)+1} = \tau + 1 - \frac{\tau}{q(\tau-1)+1} \\
&\simeq -\frac{\tau}{q(\tau-1)+1}
\end{align*}
Note that since $ -\frac{\tau}{q(\tau-1)+1}$ is a concave function in $q$ since $\frac{\partial^2}{\partial^2 q} \left( -\frac{\tau}{q(\tau-1)+1}\right) \leq 0$ and since $g_2 \simeq  -\frac{\tau}{q(\tau-1)+1}$ on the domain $\mathcal{D} = [q_-,q_+]$; thus, $\frac{\partial^2 g_2}{\partial^2 q} \leq 0$ on the domain $\mathcal{D} = [q_-,q_+]$.
\newpage

\section{Proof of Lemma 4}
We first start by showing that $q_-$ is a valid probability. First note that we are only interested in the case where $\tau > 4 \frac{n-1}{n}$, or equivalently $\tau \geq 4$ in which $q_-$ is real. Therefore,

\begin{align*}
    \frac{\partial q_-}{\partial \tau} = \frac{n - \frac{n\sqrt{n\tau}}{2\sqrt{4(1-n)+n\tau}}- \frac{n\sqrt{4(1-n)+n\tau}}{2\sqrt{n\tau}}}{2(n-1)(\tau-1)}- \frac{2 + n(\tau-2) - \sqrt{n\tau}\sqrt{4(1-n)+n\tau)}}{2(n-1)(\tau-1)^2}.
\end{align*}

Let $D = \sqrt{4(1-n) + n\tau}$, then we need to show that $ \frac{\partial q_-}{\partial \tau} \leq 0 \forall \tau$ and thus $q_-$ is a decreasing function in $\tau$. Thus we need to establish the following:
\begin{align*}
    &n(\tau-1) - \frac{n(\tau-1)\sqrt{n\tau}}{2D} - \frac{n(\tau-1)D}{2 \sqrt{n\tau}} -2 -n(\tau-2) + \sqrt{n\tau}D \leq 0  \\
    &\iff n + \sqrt{n\tau}D \leq \frac{n(\tau-1)\sqrt{n\tau}}{2D} + \frac{n(\tau-1)D}{2\sqrt{n\tau}} + 2  \\
    &\iff 2Dn\sqrt{n\tau} + 2D^2 (n\tau) \leq n(\tau-1)(n\tau) + D^2n(\tau-1) + 4D\sqrt{n\tau}  \\
    &\iff 2Dn\left[\sqrt{n\tau} + D\tau \right] \leq 2Dn \left[2\sqrt{\frac{\tau}{n}} + \frac{1}{2}D(\tau-1) \right] + n(\tau-1)(n\tau)  \\
    &\iff n(\tau-1) n\tau \ge 2Dn \left[\sqrt{n\tau} + D\tau - 2 \sqrt{\frac{\tau}{n}} - \frac{1}{2}D(\tau - 1) \right]  \\
    &\iff n(\tau-1) n\tau \ge 2Dn \left[\sqrt{n\tau}  - 2 \sqrt{\frac{\tau}{n}} + \frac{1}{2}D(\tau + 1) \right] \\
    &\iff n(\tau-1) n\tau \ge D^2n(\tau + 1)  + 2Dn \sqrt{n\tau}  - 4 D\sqrt{n \tau}  \\
    &\iff n^2 \tau^2 - n^2 \tau \ge (4 - 4n + n\tau)(\tau+1)n + 2Dn \sqrt{n\tau}  - 4 D\sqrt{n \tau} \\
    &\iff n^2 \tau^2 - n^2 \tau \ge (4n\tau + 4n - 3n^2\tau - 4n^2 + n^2\tau^2) + 2Dn \sqrt{n\tau}  - 4 D\sqrt{n \tau} \iff \\
    &\iff - n^2 \tau \ge (4n\tau + 4n - 3n^2\tau - 4n^2) + 2Dn \sqrt{n\tau}  - 4 D\sqrt{n \tau}  \\
    &\iff 2n^2 \tau + 4n^2\ge 4n(\tau + 1) + 2D\sqrt{n \tau} \left(n-2\right) \\
    &\iff n^2 \tau + 2n^2\ge2n(\tau + 1) + D\sqrt{n \tau} \left(n-2\right)  \\
    &\iff n^2 (\tau + 2) \ge 2n(\tau + 1) + D\sqrt{n \tau} \left(n-2\right) \\
\end{align*}

Note that since:
\begin{align*}
    D = \sqrt{4(1-n)+n\tau} \leq \sqrt{n\tau}
\end{align*}
Then:
\begin{align*}
   2n(\tau+1) + D\sqrt{n\tau}(n-2) \leq 2n(\tau+1) + n\tau(n-2)  = n(n\tau +2) \leq n^2(\tau+2)
\end{align*}

Thus the last inequality always holds which shows that $q_{-}$ is a decreasing function of $\tau$. Furthermore, $\lim_{\tau \to 4} q_- = \frac{n+1-2\sqrt{n}}{3(n-1)}$ and that $\lim_{\tau \to n} q_- = \frac{1}{(n-1)^2}$. Therefore, since $\frac{1}{(n-1)^2} \leq q_- \leq \frac{n+1-2\sqrt{n}}{3(n-1)}$, $q_-$ is a valid probability. A similar argument can be used to show that $q_+$ is increasing in $\tau$ and that it is a valid probability.
\newpage

\section{Proof of Lemma 5}
To find the intersection points between $g_1$ and $g_2$, we consider the two cases of $\rho$ separately. For a simplification, note that the intersection point between $g_1$ and $g_2$ is the same as for $\frac{g_1}{q(\tau-1)+1}$ and $\frac{g_2}{q(\tau-1)+1}$.

\textbf{Case I:} $q\theta^2 \leq \frac{n}{\tau} \frac{n-1}{\tau-1}$

\begin{align*}
    \frac{g_2- g_1}{q(\tau-1)+1} &= \theta \left(1+ \frac{4\theta L_{\max}}{\mu n} \left(\frac{1-q}{n} + q \frac{\tau}{n}\frac{n-\tau}{n-1} \right)\right) - \frac{4\theta L_{\max}}{\mu n}\frac{q(\tau^2-1)+1}{q(\tau-1)+1} = 0 \\
    &\frac{n}{q(\tau-1) +1} + \frac{4L_{\max}}{\mu} \frac{1}{q(\tau-1)+1} \left(\frac{nq\tau(1-\tau)}{(n-1)(q(\tau-1)+1)}\right) = 0 \\
    &\frac{4L_{\max}}{\mu} \frac{q \tau (\tau-1)}{(n-1)(q(\tau-1)+1)} = 1 \\
    &q_{i1} = \frac{1-n}{(\tau-1)(n-1) - \tau(\tau-1) \frac{4L_{\max}}{\mu}} = \frac{n-1}{(\tau-1)(\tau \frac{4L_{\max}}{\mu}+1 -n)}
\end{align*}

Now that we have the intersection point, we can next find the ranges of $\tau$ in which $q_{i1}$ occur. In particular, the range of $\tau$ in which the following is satisfied $q_{i1} \theta^2 \leq \frac{n}{\tau} \frac{n-1}{\tau-1}$. Note that for $q_{i1}$, we have $\theta = \frac{n(-n+1+\tau \frac{4L_{\max}}{\mu})}{\tau \frac{4L_{\max}}{\mu}}$. Therefore, it is easy now to show that for 


\textbf{(a)}: if $n \leq \frac{4L_{\max}}{\mu}$ then $q_{i1} \ge 0$. However, for $q_{i1} \leq 1$ then we have:

\begin{align*}
    &n-1 \leq \left(\tau-1\right) \left(\tau \frac{4L_{\max}}{\mu} + 1 - n\right) \\
    & \tau \left(\tau \frac{4L_{\max}}{\mu} - \frac{4L_{\max}}{\mu} -n + 1\right) \ge 0
\end{align*}

Therefore for $0 \leq q_{i1} \leq 1$ where $n \leq \frac{4L_{\max}}{\mu}$ then:

\begin{align*}
    \tau \ge \tau_{\text{min}} = \frac{n}{\frac{4L_{\max}}{\mu}} +1 - \frac{\mu}{4L_{\max}}
\end{align*}

\textbf{(b)}: Similarly , if $n > \frac{4L_{\max}}{\mu}$ then for $0 \leq q_{i1} \leq 1$ we we need:

\begin{align*}
    \tau \leq \tau_{\text{max}} = \left(\frac{n(n-1)\mu}{(n - \frac{4L_{\max}}{\mu})4L_{\max}}\right)
\end{align*}

Thus the range of
\begin{align*}
    \tau \in \left[\max\left(4, \tau_{\text{min}}\right), \min\left(\tau_{\text{max}},n\right)\mathbbm{1}_{n > \frac{4L_{\max}}{\mu}} + n\mathbbm{1}_{n \leq \frac{4L_{\max}}{\mu}}\right]
\end{align*}

\newpage
\textbf{Case II:} $q\theta^2 \ge \frac{n}{\tau} \frac{n-1}{\tau-1}$

\begin{align*}
    \frac{g_2- g_1}{q(\tau-1)+1}  &= \theta \left(1+ \frac{4}{\mu n^2} \frac{nL_{\text{max}}}{\theta}\left(\theta^2 \left(\frac{1-q}{n} + q\frac{\tau}{n}\frac{n -\tau}{n-1}\right)+n\left(\theta^2q \frac{\tau}{n}\frac{\tau-1}{n-1}-1\right)\right)\right)\\
    &\Rightarrow -  \frac{4\theta L_{\text{max}}}{\mu n}\frac{q(\tau^2-1)+1}{q(\tau-1)+1} = 0\\
    &\Rightarrow \theta+ \frac{4}{\mu n^2} nL_{\text{max}}\left(\theta^2 \left(\frac{1-q}{n} + q\frac{\tau}{n}\frac{n -\tau}{n-1}\right)+n\left(\theta^2q \frac{\tau}{n}\frac{\tau-1}{n-1}-1\right)\right)\\
    &-  \frac{4\theta L_{\text{max}}}{\mu n}\frac{q(\tau^2-1)+1}{q(\tau-1)+1} = 0\\
    &\Rightarrow \theta+ \frac{4\theta L_{\text{max}}}{\mu n} \left(\theta \left(\frac{1-q}{n} + q\frac{\tau}{n}\frac{n -\tau}{n-1}\right)+n\left(\theta q \frac{\tau}{n}\frac{\tau-1}{n-1}-\frac{1}{\theta}\right)-\frac{q(\tau^2-1)+1}{q(\tau-1)+1}\right) = 0 \\
    &\Rightarrow \theta+ \frac{4\theta L_{\text{max}}}{\mu n} \left(\frac{nq\tau(1-\tau)}{(n-1)(q(\tau-1)+1)} +n\left(\theta q \frac{\tau}{n}\frac{\tau-1}{n-1}-\frac{1}{\theta}\right)\right) = 0 \\
    &\Rightarrow 1 = \frac{4 L_{\text{max}}}{\mu n} \left(\frac{nq\tau(\tau - 1)}{(n-1)(q(\tau-1)+1)} -n\left(\theta q \frac{\tau}{n}\frac{\tau-1}{n-1}-\frac{1}{\theta}\right)\right) = 0 \\
    &\Rightarrow 1 = \frac{4 L_{\text{max}}}{\mu n} \left(\frac{nq\tau(\tau - 1)}{(n-1)(q(\tau-1)+1)} -n\left(\frac{q\tau}{q(\tau-1)+1} \frac{\tau-1}{n-1}-\frac{1}{\theta}\right)\right) = 0 \\
    &\Rightarrow 1 = \frac{4 L_{\text{max}}}{\mu n} \left(\frac{nq\tau(\tau - 1)}{(n-1)(q(\tau-1)+1)} -\left(\frac{qn\tau}{q(\tau-1)+1} \frac{\tau-1}{n-1}-(q(\tau-1)+1)\right)\right) = 0 \\
    &\Rightarrow 1 = \frac{4 L_{\text{max}}}{\mu n} \left(q(\tau-1)+1\right) = 0  ~~~~~~~ \Rightarrow q_{i2} = \frac{1 - \frac{4L_{\text{max}}}{\mu n}}{\frac{4L_{\text{max}}}{\mu n}(\tau-1)} = \frac{n - \frac{4L_{\max}}{\mu}}{\frac{4L_{\max}}{\mu}\left(\tau - 1\right)}
\end{align*}

and that $\theta = \frac{4L_{\max}}{\mu}$. For $q_{i2} \theta^2 \ge \frac{n(n-1)}{\tau(\tau-1)}$ then 

\begin{align*}
    &\frac{4L_{\max}}{\mu} \tau \left(n - \frac{4L_{\max}}{\mu}\right) \ge n(n-1) \\
    & \tau \ge \frac{n(n-1)}{(n - \frac{4L_{\max}}{\mu})\frac{4L_{\max}}{\mu}} = \tau_{\text{max}}
\end{align*}

Then the range of $\tau$ that assures $0 \leq q_{i2} \leq 1$ are as follows:
\begin{align*}
    \tau \in \left[\tau_{\text{max}},n\right]
\end{align*}

\newpage
\section{Experiments with logistic loss}
\begin{figure}[h]
\begin{subfigure}[ht]{0.99\linewidth}
\includegraphics[width=0.33\textwidth]{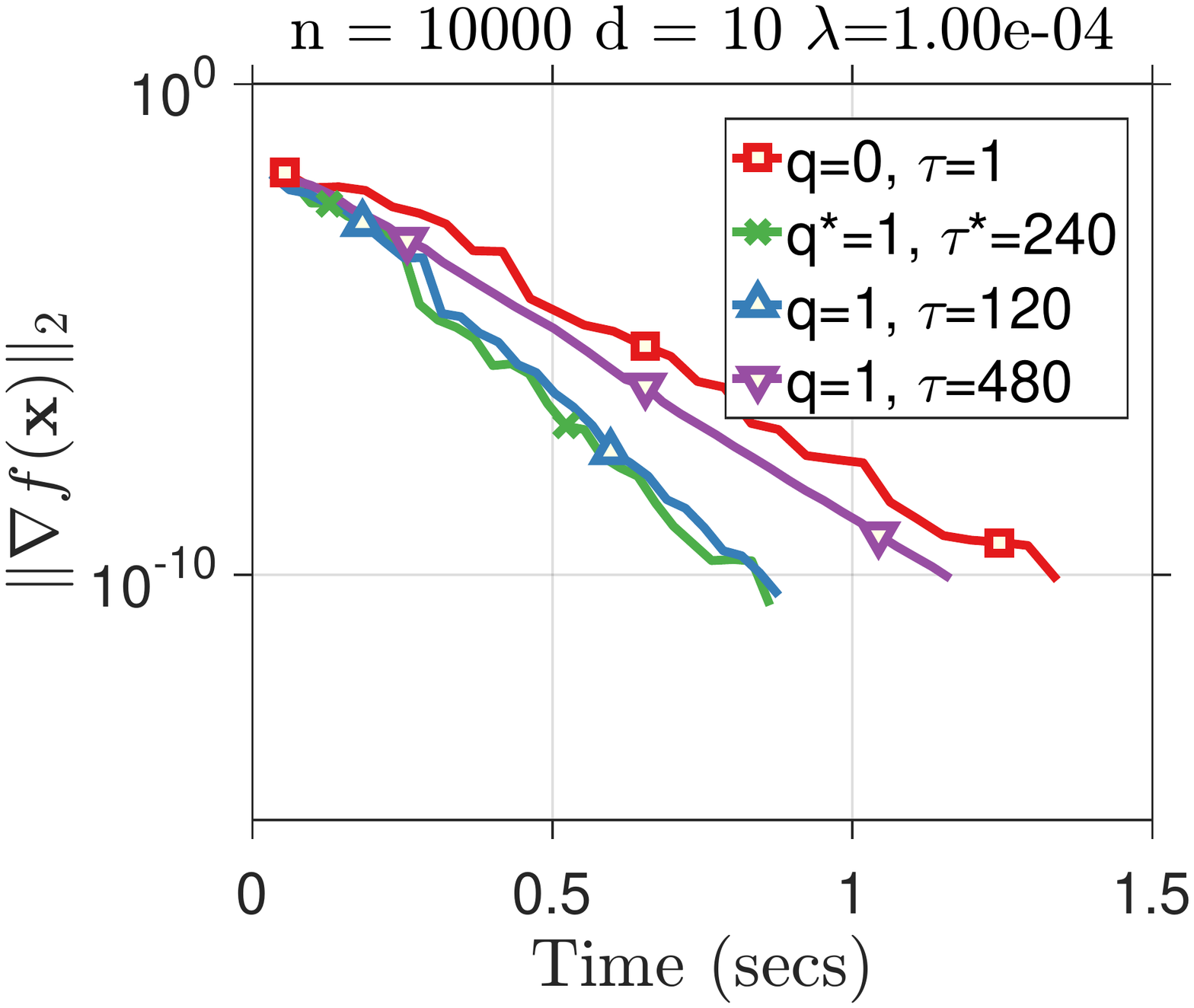}
\includegraphics[width=0.33\textwidth]{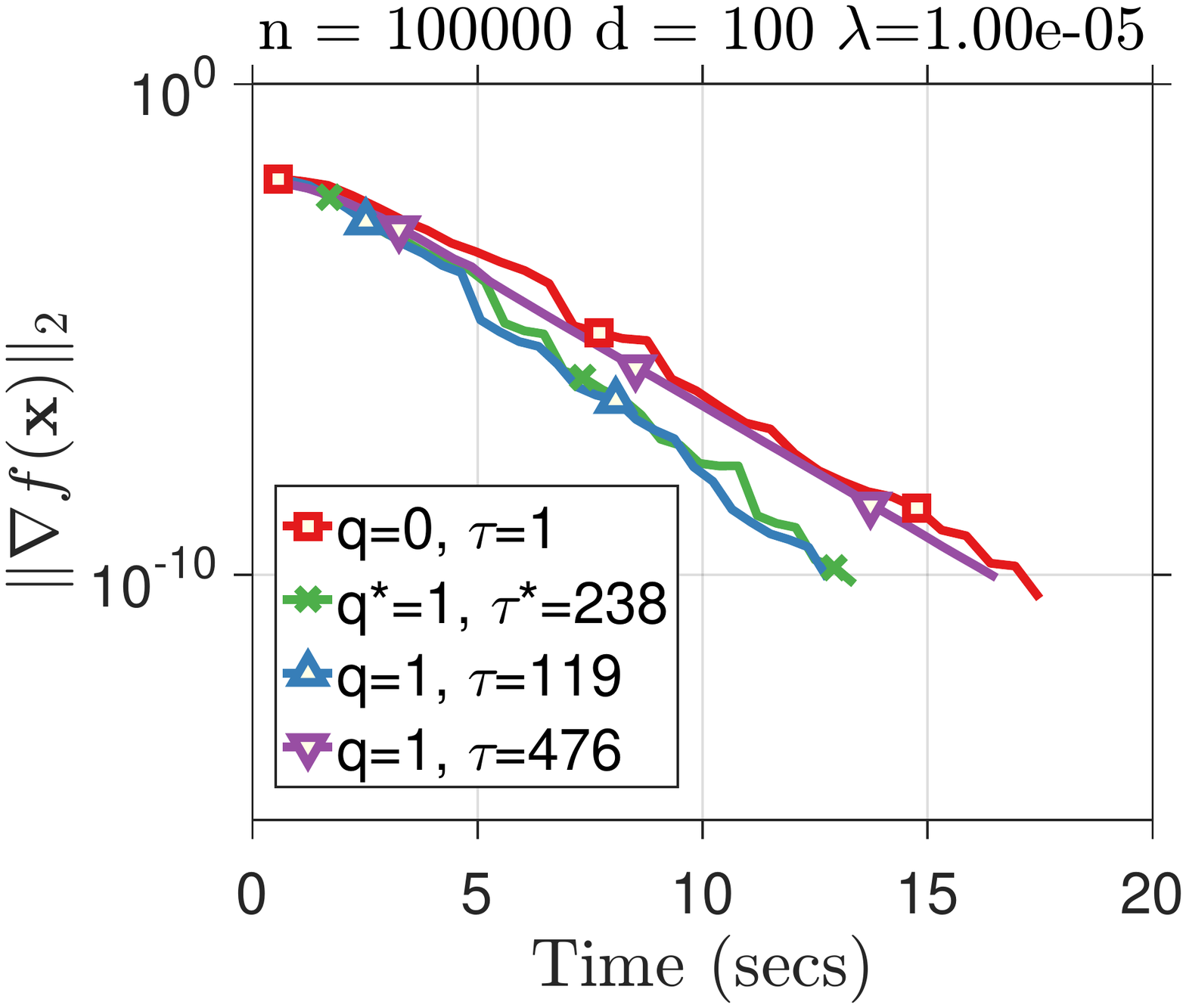}
\includegraphics[width=0.33\textwidth]{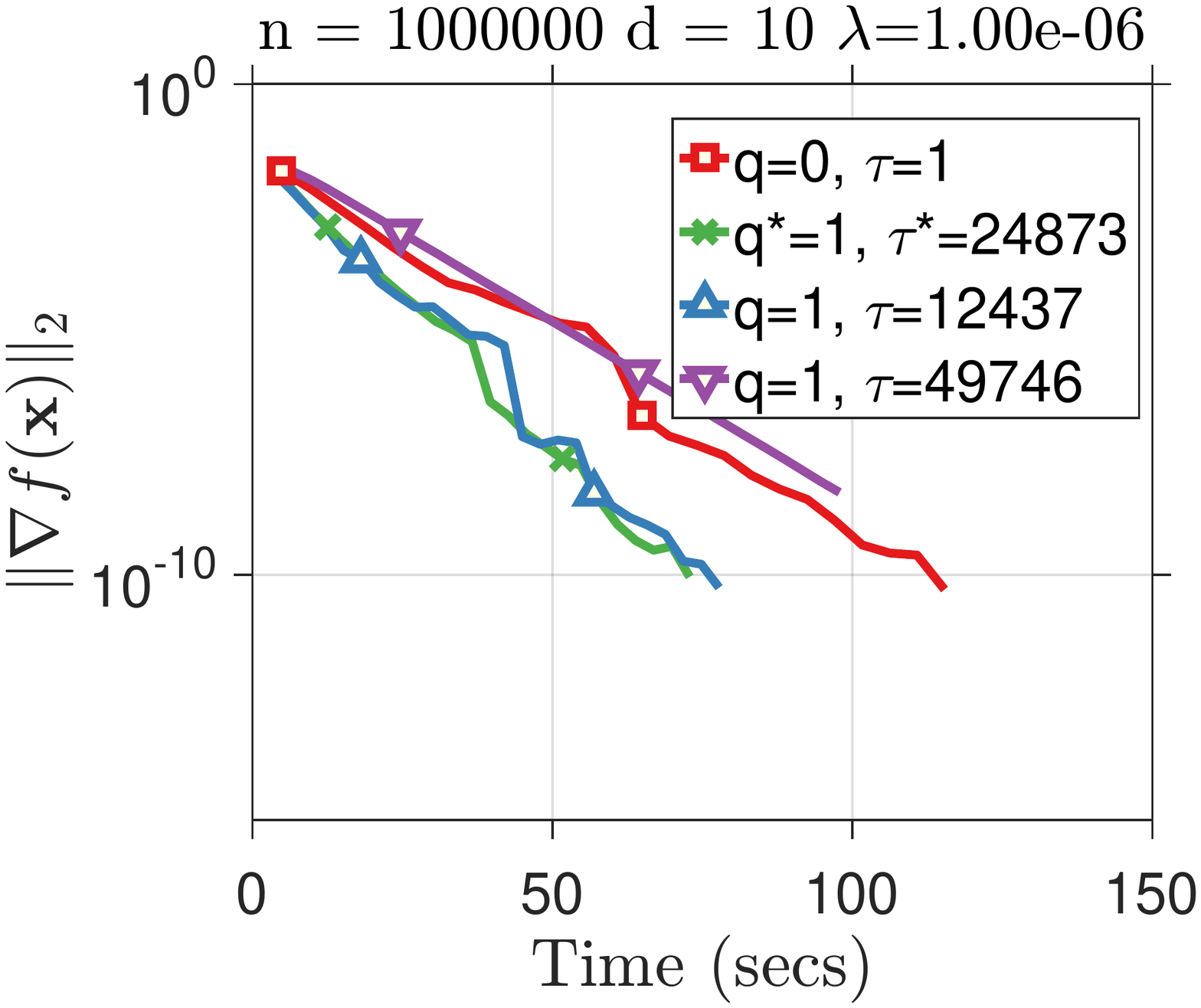}
\vspace{-2.5cm}
\end{subfigure}
\begin{subfigure}[ht]{0.99\linewidth}
\includegraphics[width=0.33\textwidth]{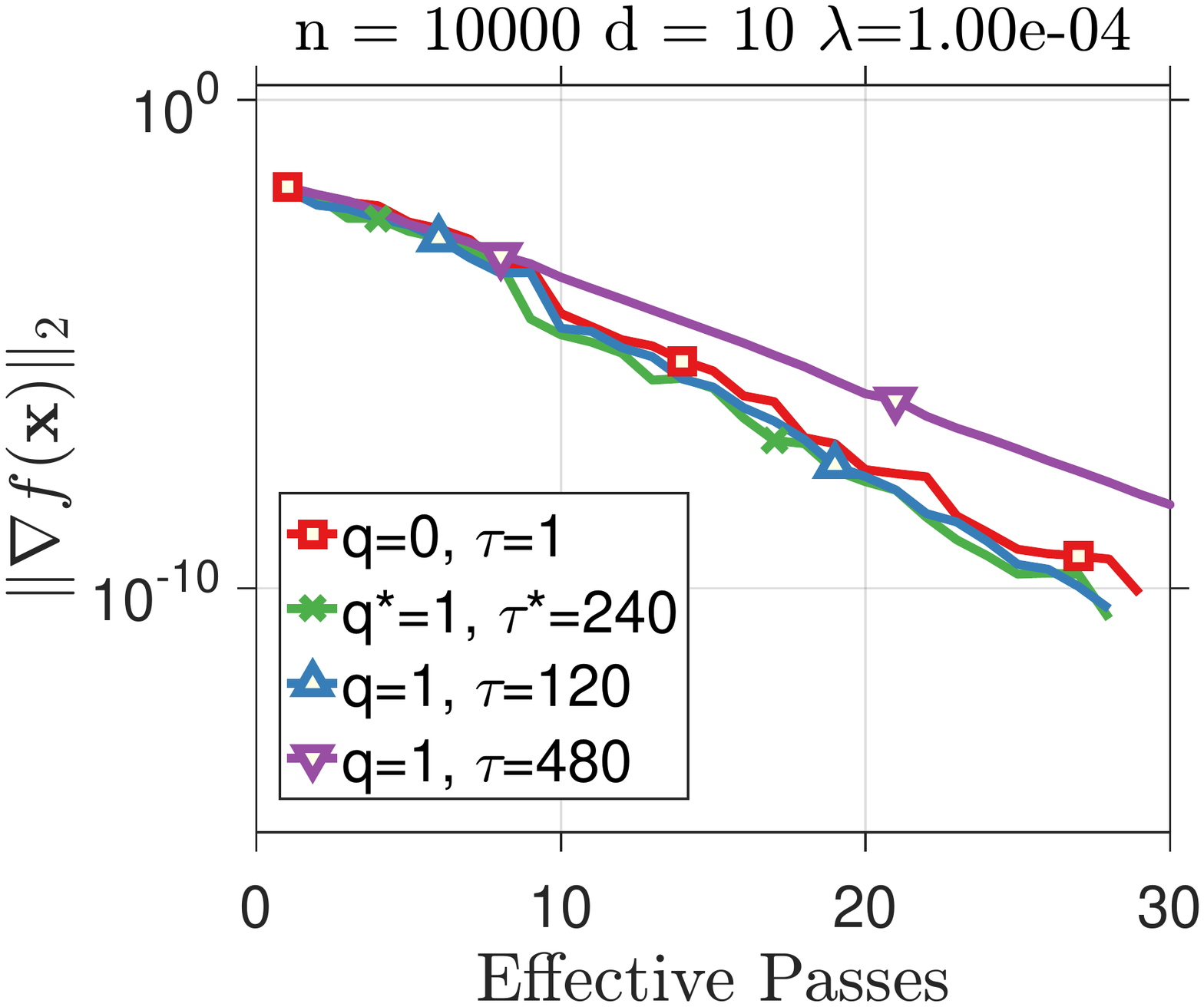}
\includegraphics[width=0.33\textwidth]{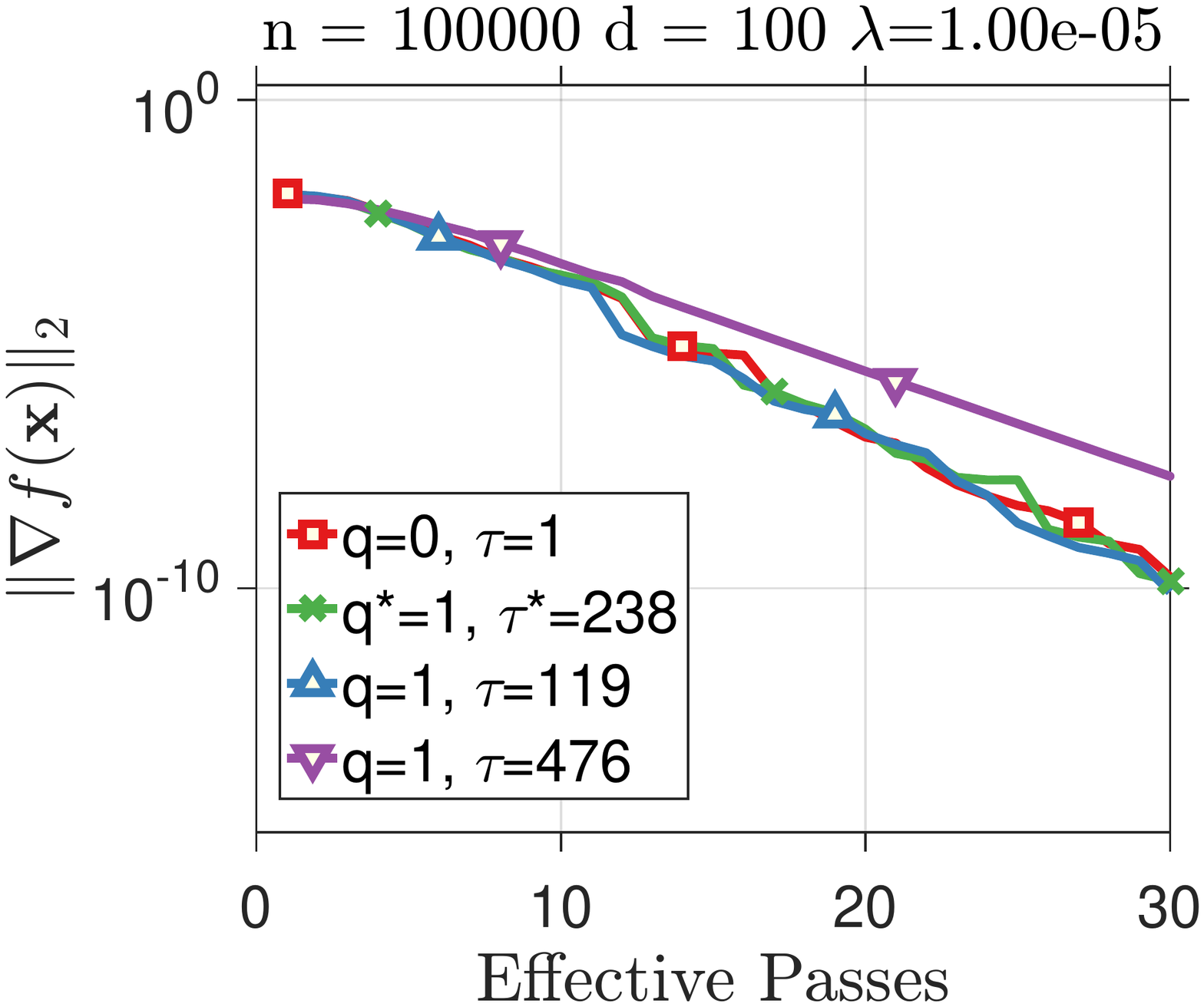}
\includegraphics[width=0.33\textwidth]{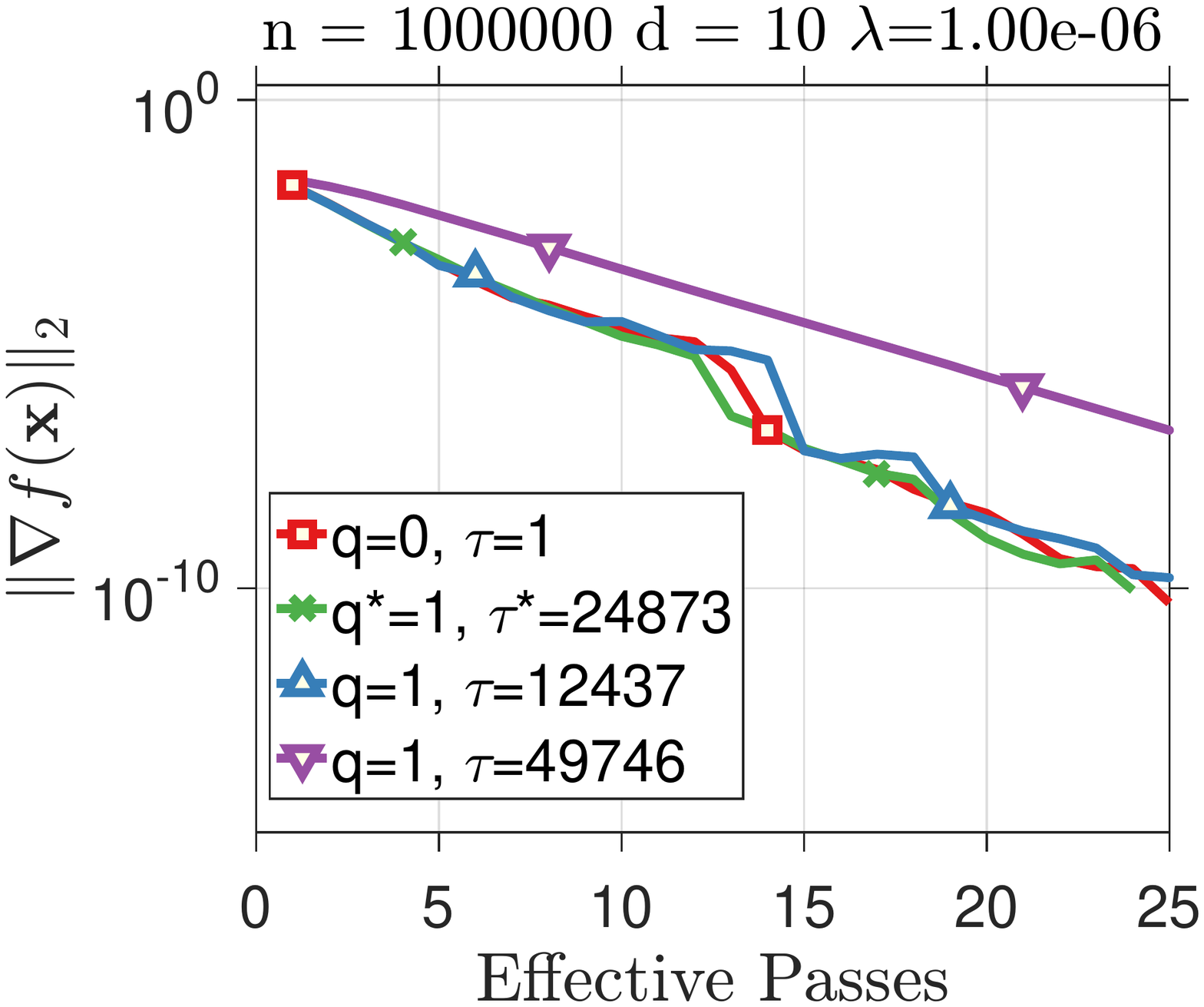}
\vspace{-0.5cm}
\end{subfigure}
\caption{Shows a comparison between SAGA and our method with optimal ($q^*,\tau^*$) pair. We also compare them both against ($q^*,\frac{\tau^*}{2}$) and ($q^*,2\tau^*$). The first row shows a comparison in time while the second row shows the same experiment as compared in the total number of effective passes.}
\label{first_time_epocs_comp_logistic}
\end{figure}

\begin{figure}[t]
\begin{subfigure}[ht]{0.99\linewidth}
\includegraphics[width=0.33\textwidth]{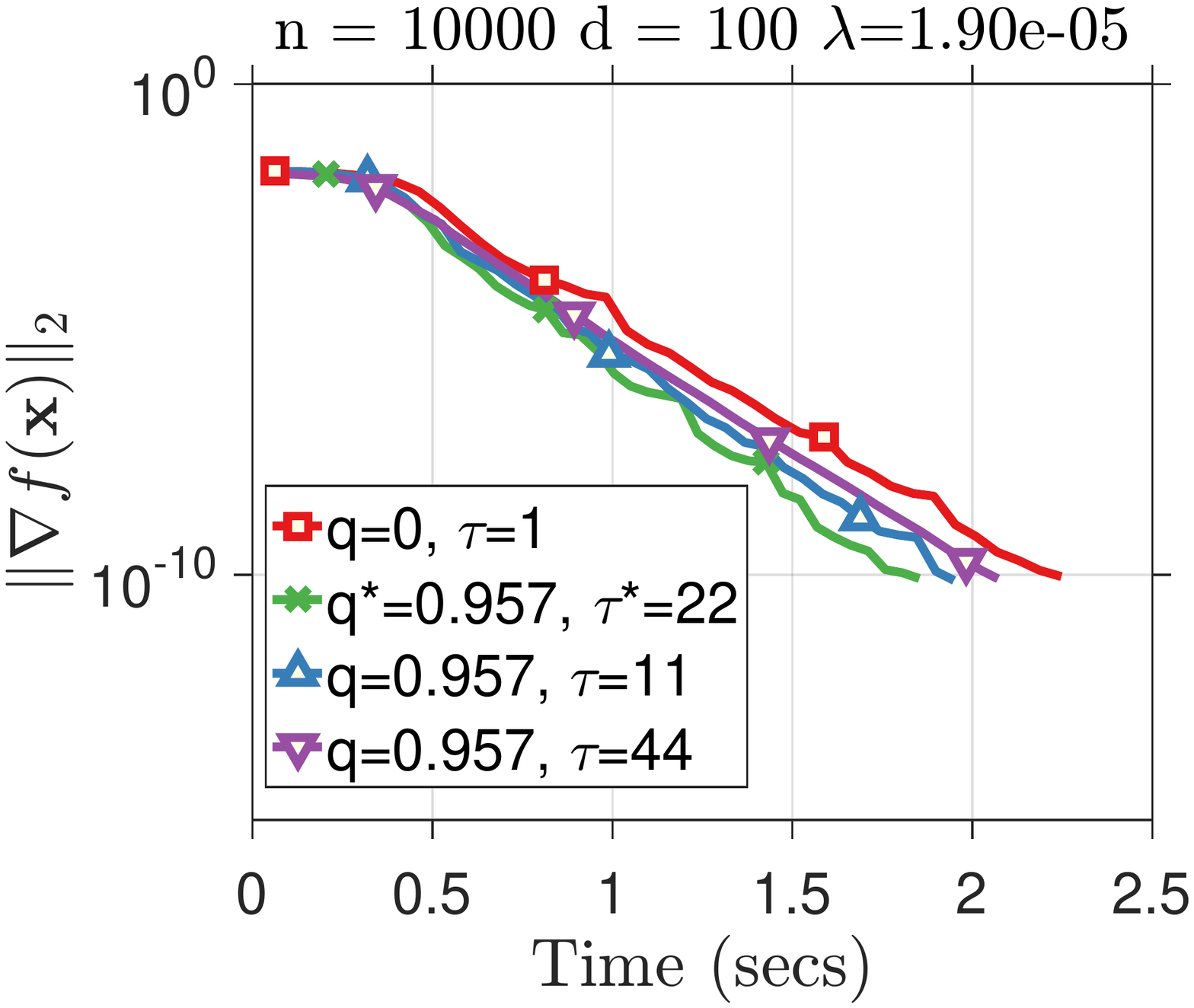}
\includegraphics[width=0.33\textwidth]{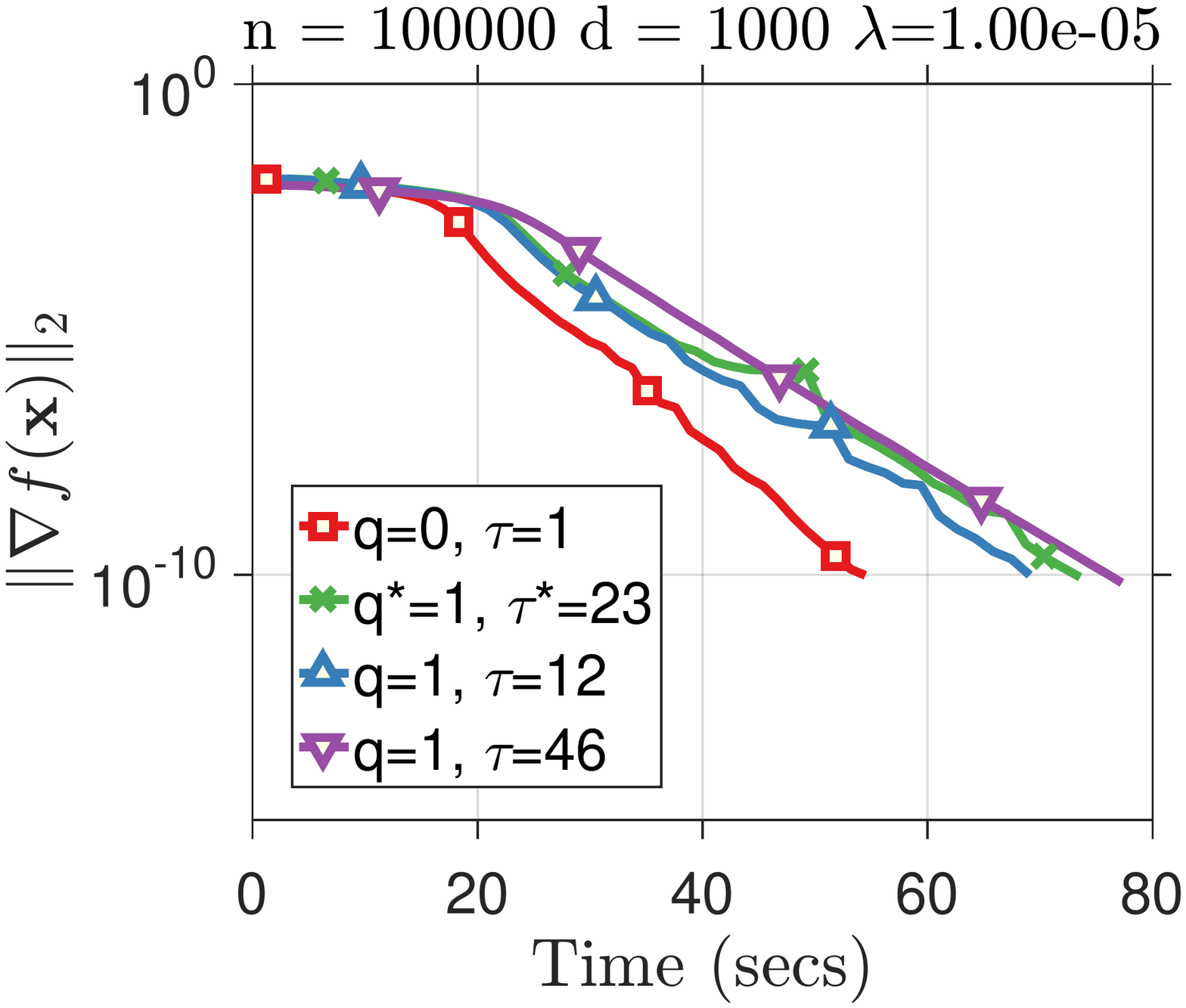}
\includegraphics[width=0.33\textwidth]{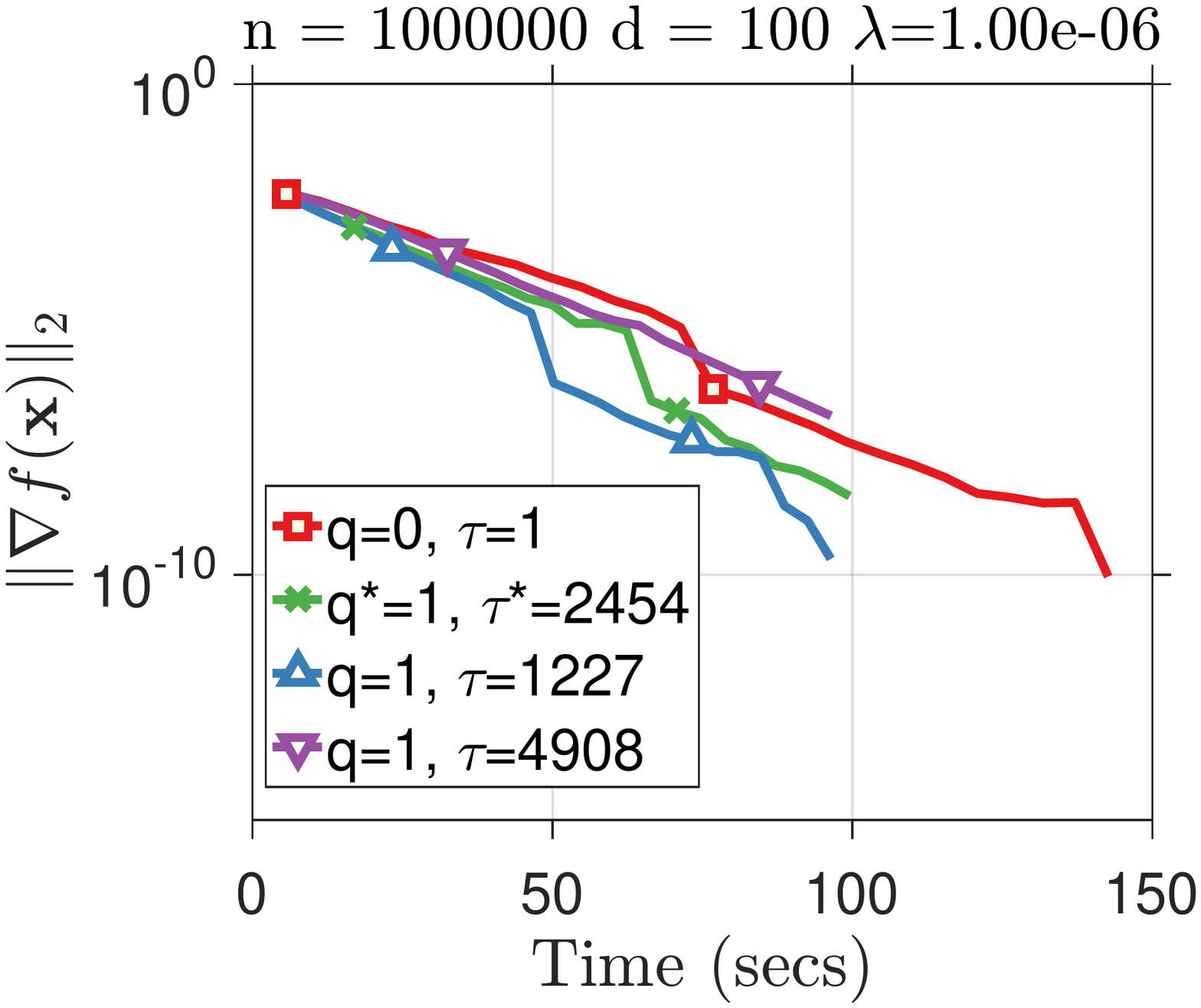}
\vspace{-2.5cm}
\end{subfigure}
\begin{subfigure}[t]{0.99\linewidth}
\includegraphics[width=0.33\textwidth]{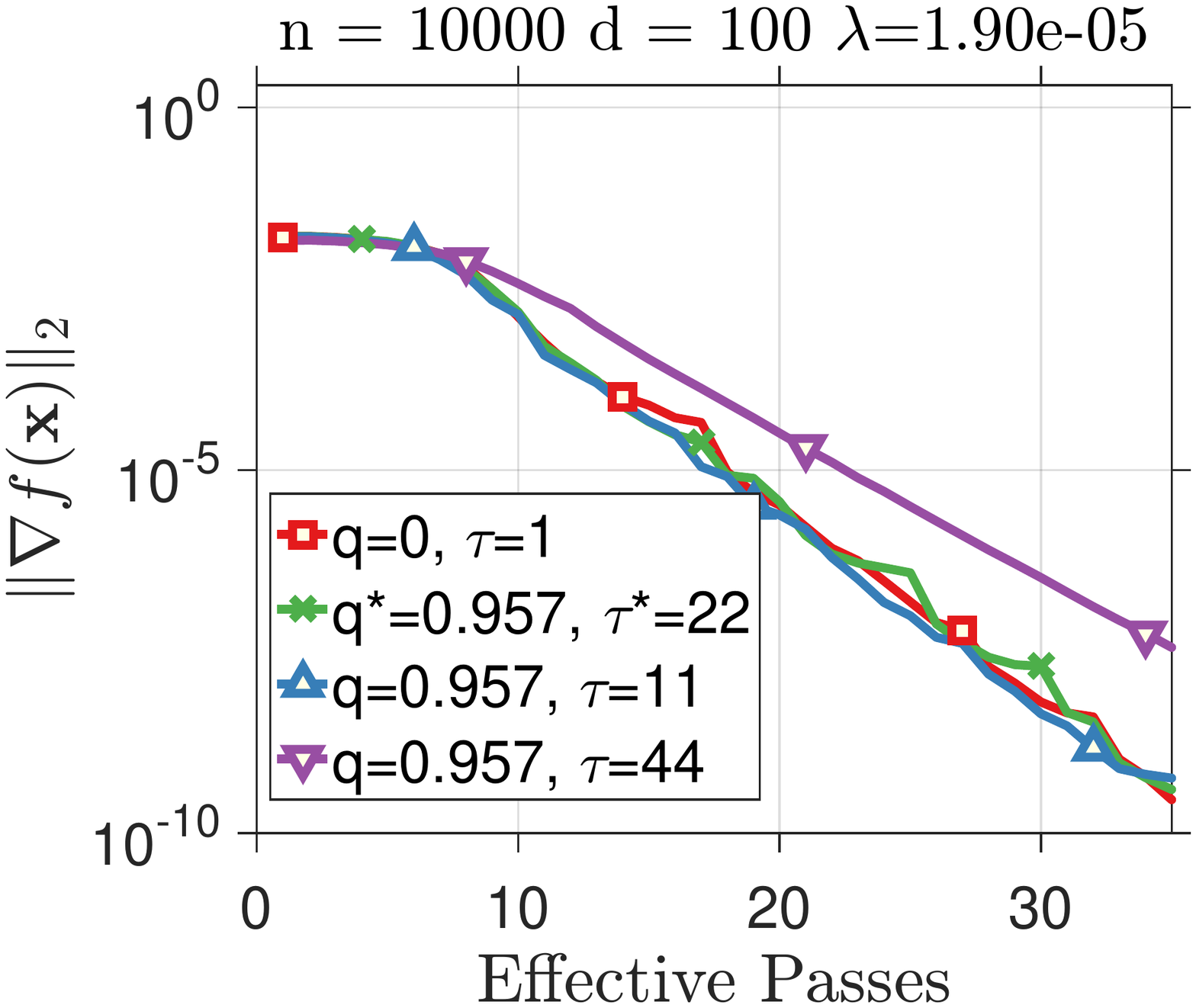}
\includegraphics[width=0.33\textwidth]{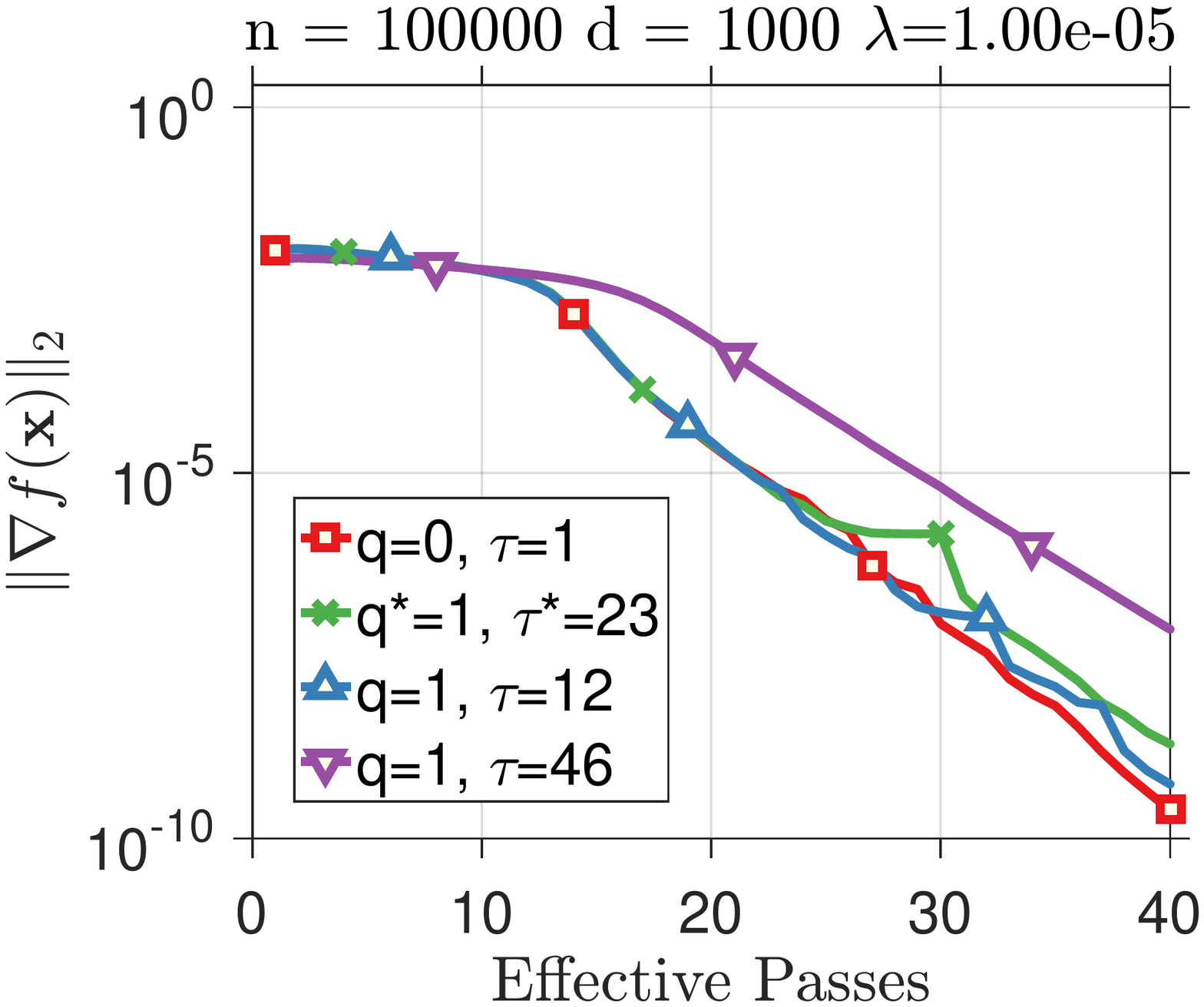}
\includegraphics[width=0.33\textwidth]{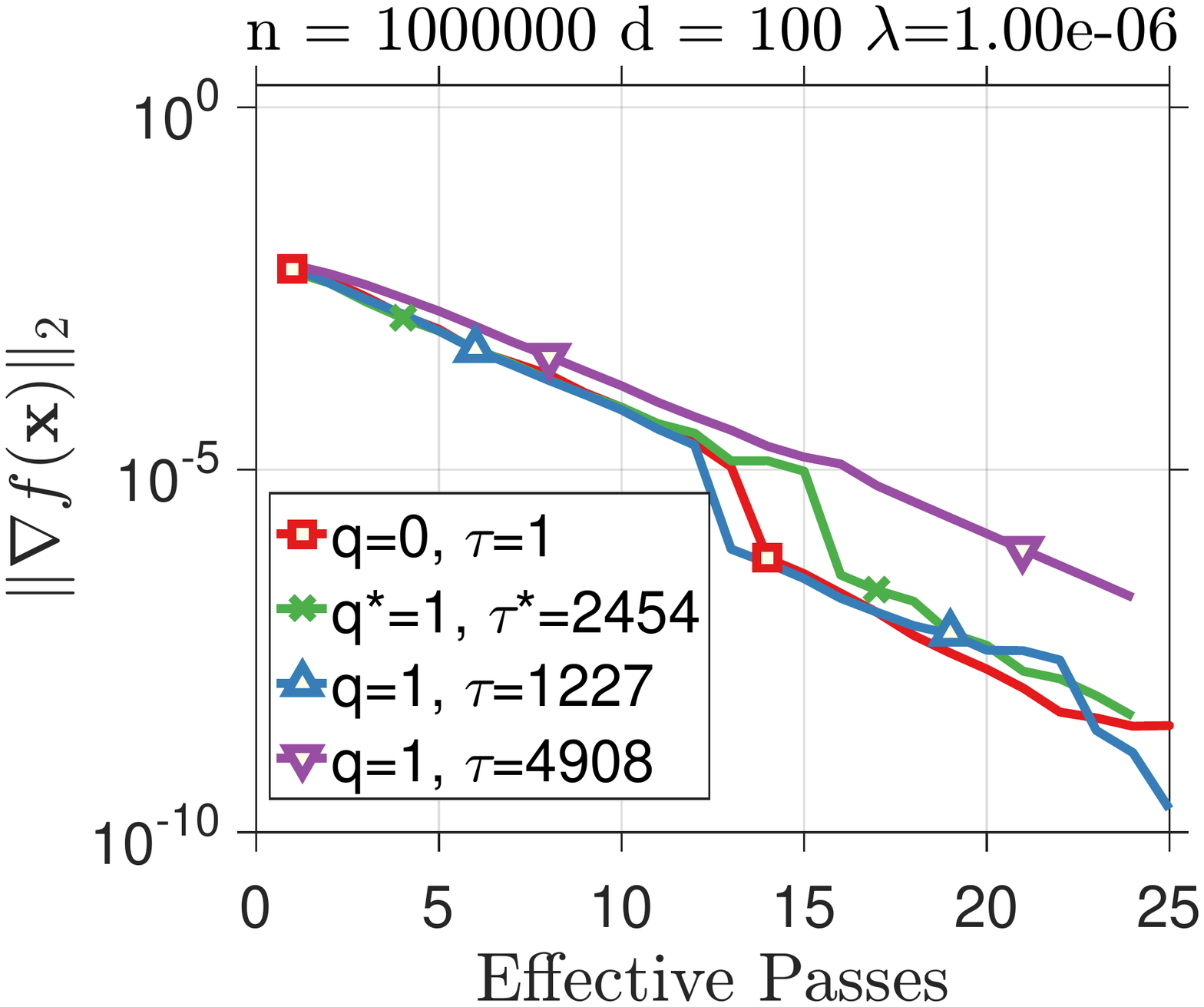}
\end{subfigure}
\vspace{-0.5cm}
\caption{Shows a comparison between SAGA and our method with optimal ($q^*,\tau^*$) pair. We also compare them both against ($q^*,\frac{\tau^*}{2}$) and ($q^*,2\tau^*$). The first row shows a comparison in time while the second row shows the same experiment as compared in the total number of effective passes.}
\vspace{-0.1cm}
\label{second_time_epocs_comp_logistic}
\end{figure}

\begin{figure}[t]
\begin{subfigure}[t]{0.99\linewidth}
\includegraphics[width=0.33\textwidth]{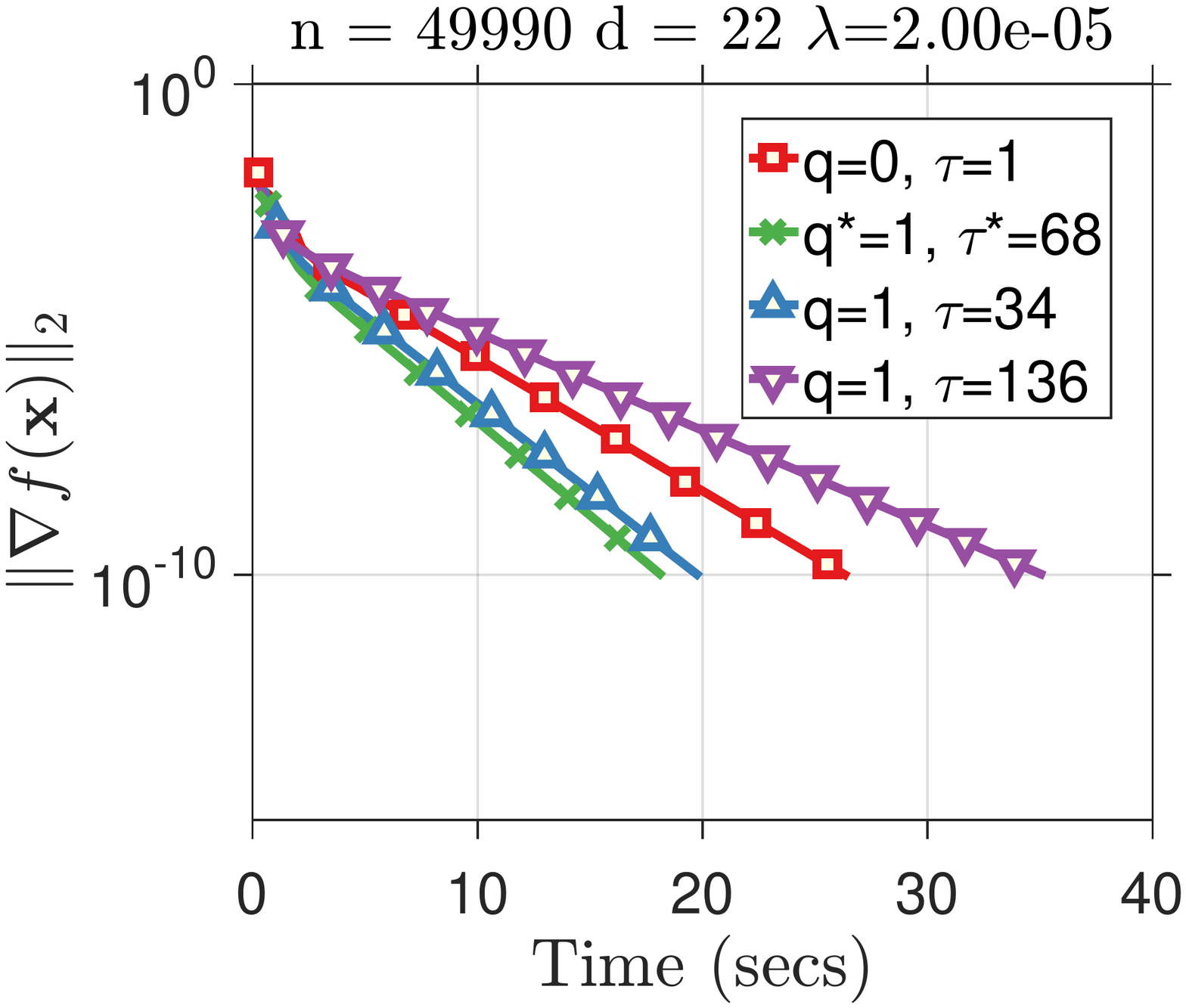}
\includegraphics[width=0.33\textwidth]{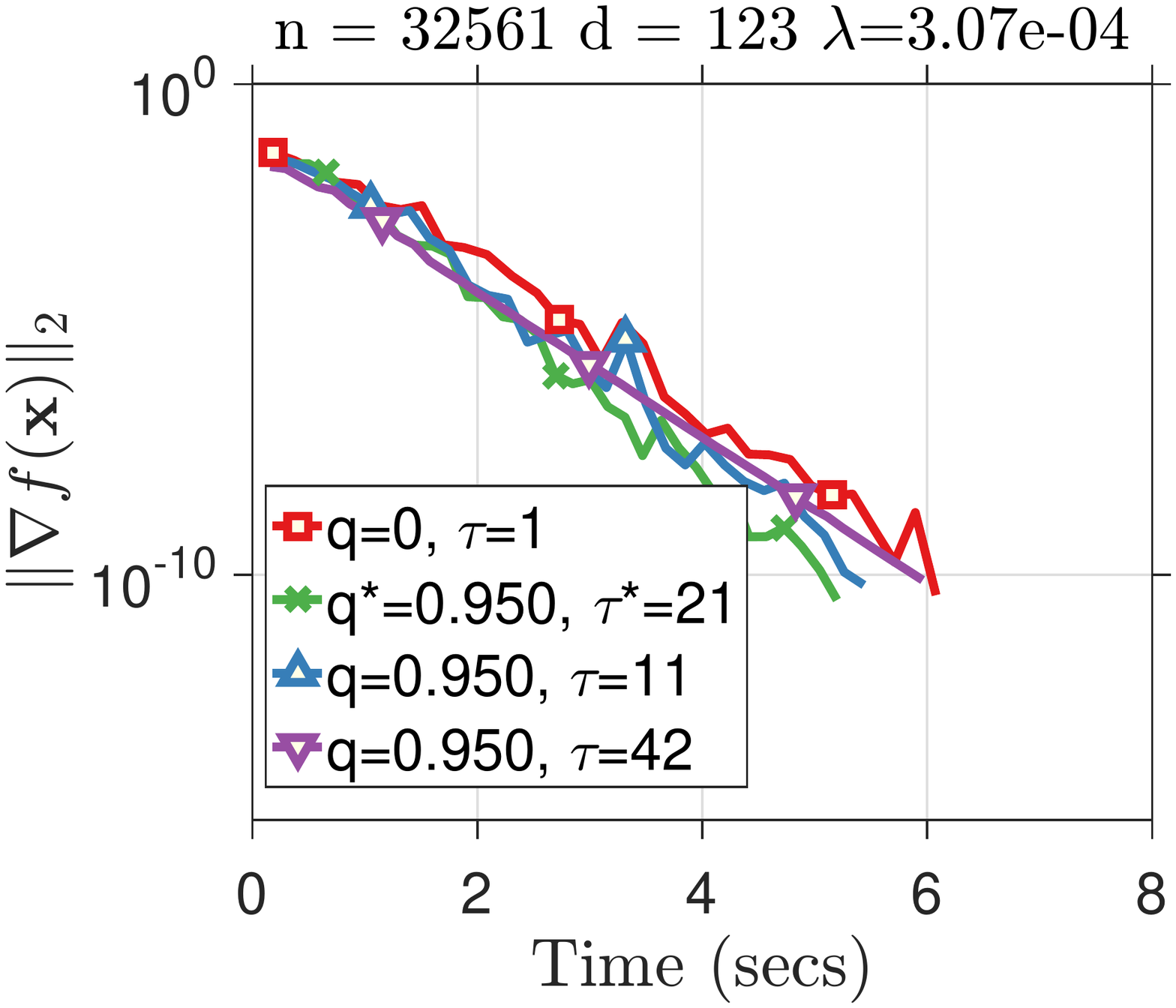}
\includegraphics[width=0.33\textwidth]{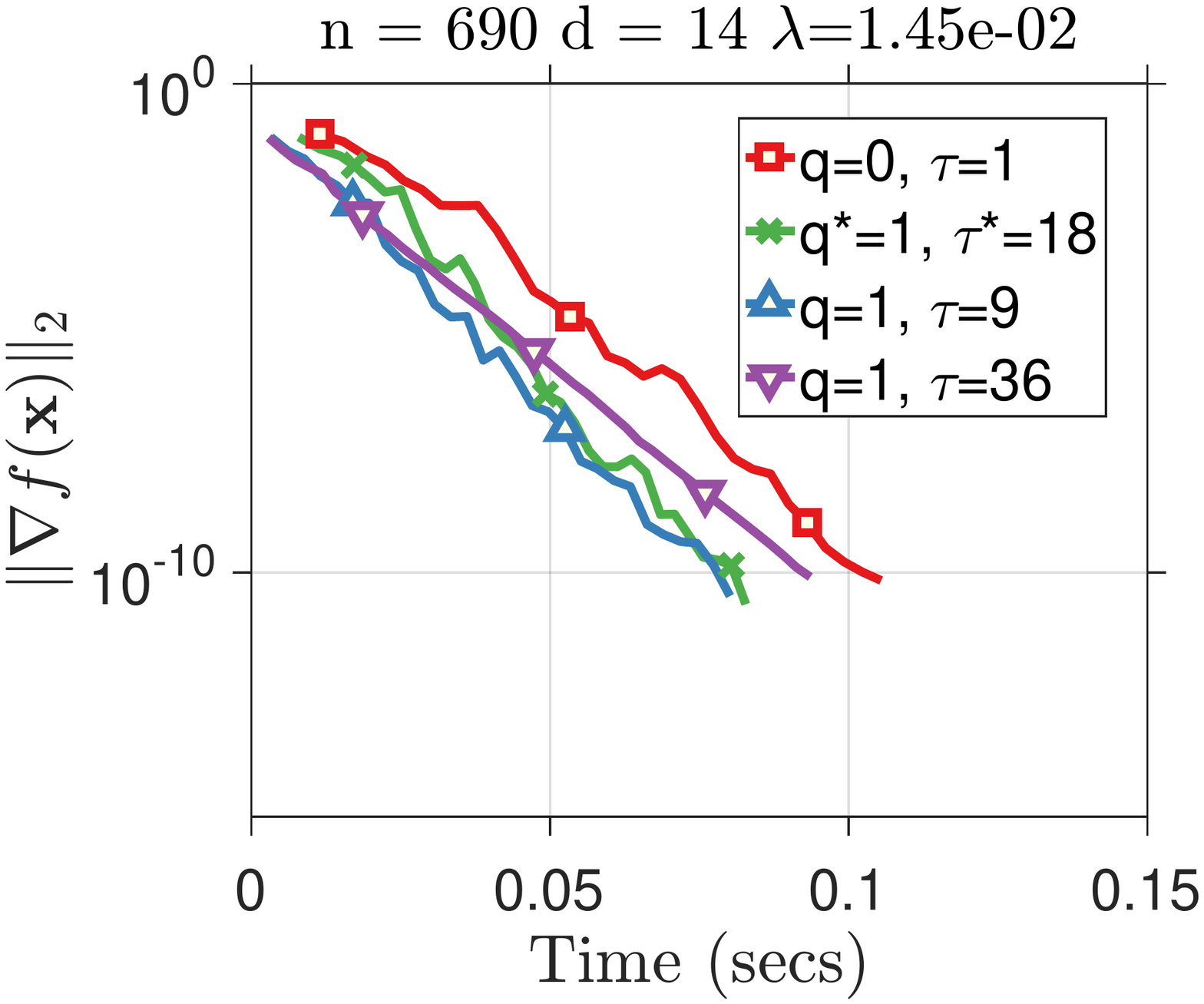}
\end{subfigure}
\begin{subfigure}[t]{0.99\linewidth}
\includegraphics[width=0.33\textwidth]{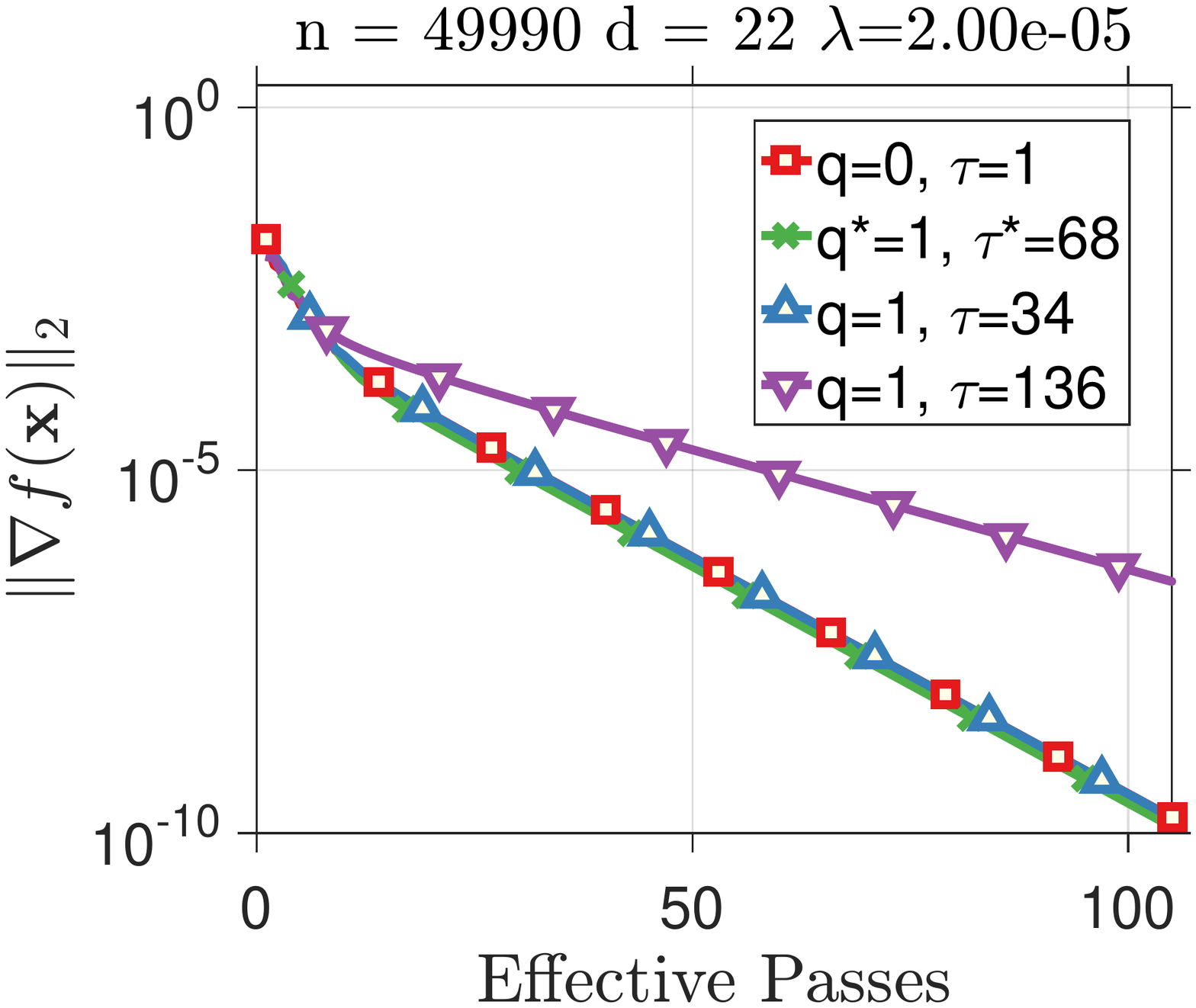}
\includegraphics[width=0.33\textwidth]{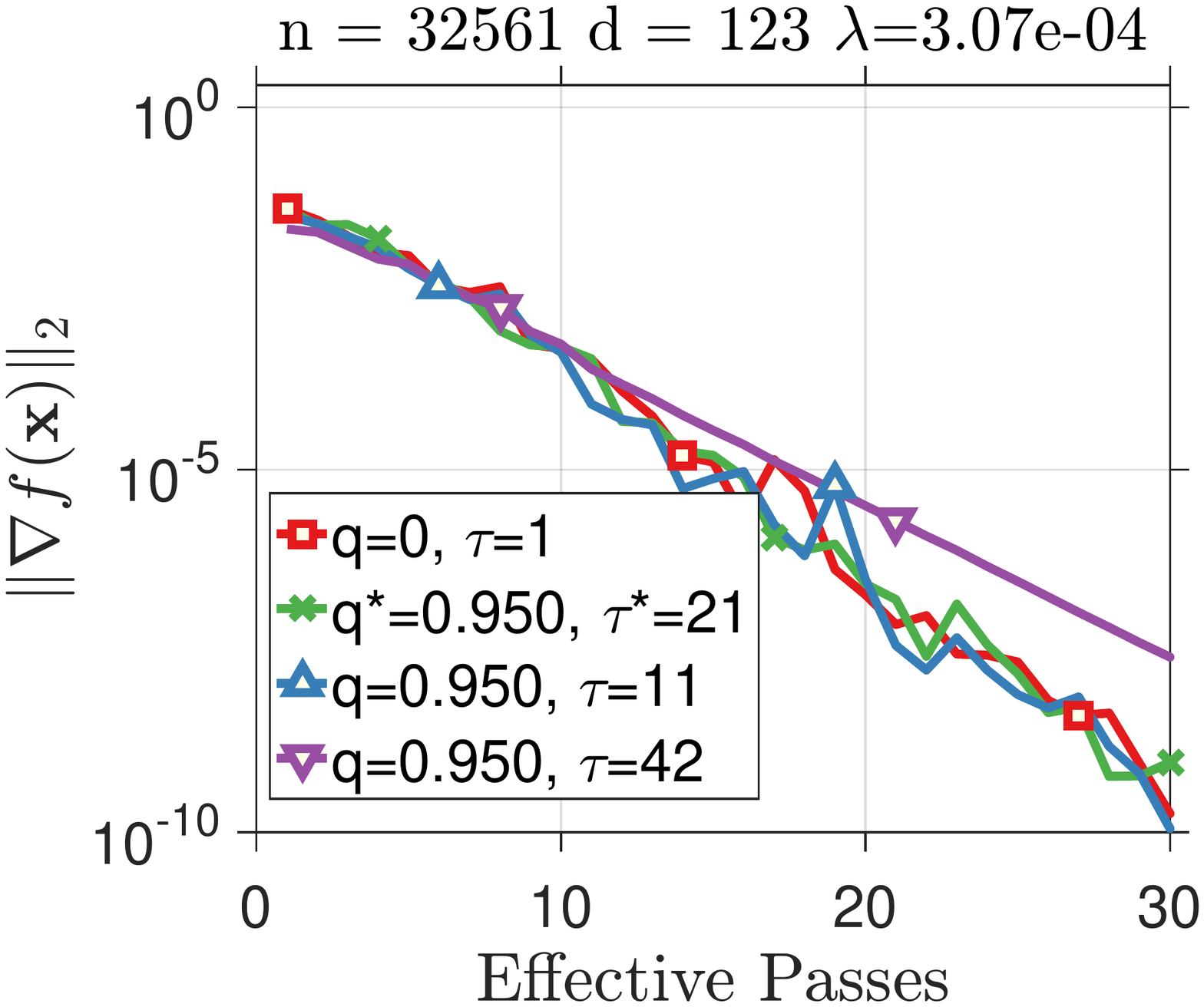}
\includegraphics[width=0.33\textwidth]{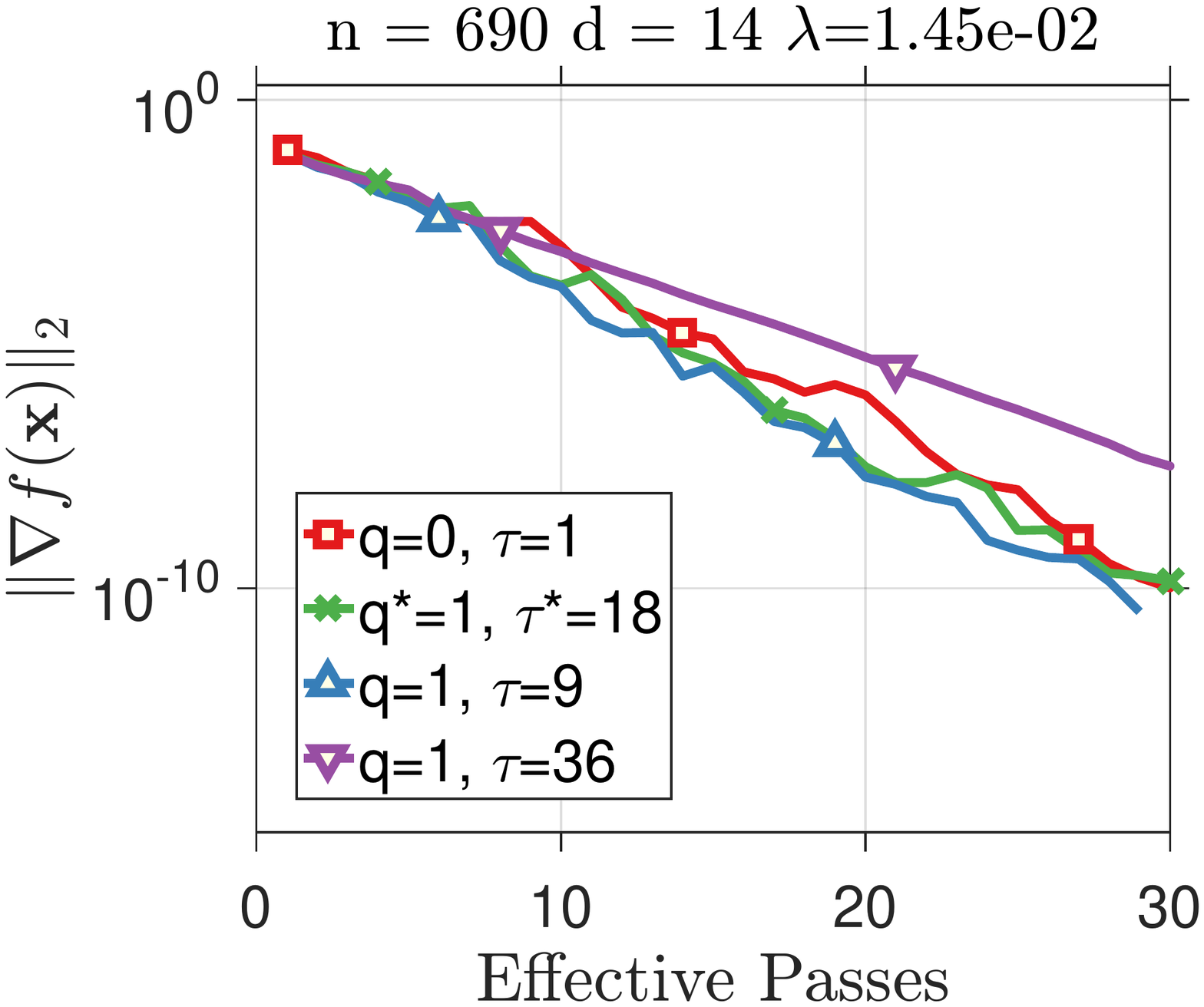}
\end{subfigure}
\caption{Shows a similar comparison but on the real datasets ijcnn1, a9a and australian.}
\label{real_dataset_logitic}
\end{figure}

In this section we conduct some further experiments on the logistic loss:
\begin{align*}
    f(x) = \frac{1}{2n} \sum_{i=1}^n \log\left(1+ \exp(-y_i a_i^\top x)\right) + \frac{\lambda}{2} \|x\|_2^2
\end{align*}
The experiments are conducted on both synthetic and real datasets similarly to section~\ref{experiments}. Figure~\ref{first_time_epocs_comp_logistic} demonstrates time comparisons and epochs comparisons among SAGA, SAGD with ($q^*,\tau^*$), SAGD with ($q^*, \frac{\tau^*}{2}$) and ($q^*,2\tau^*$) on three different datasets. SAGD with optimal ($q^*,\tau^*$) pair achieve the best performance and suggest a linear speed up in minibatch size for a parallel implementation. The same conclusion can be drawn from the the second synthetic datasets~\ref{second_time_epocs_comp_logistic} that has a larger dimensional problem. Moreover we conduct some further experiments on real datasets~\ref{real_dataset_logitic}, namely ijcnn1, a9a and australian, demonstrating the effectivness of SAGD.


%% file: main.bbl
\begin{thebibliography}{10}

\bibitem{shai_book}
S.~Shalev-Shwartz and S.~Ben-David, {\em Understanding machine learning: from
  theory to algorithms}.
\newblock Cambridge University Press, 2014.

\bibitem{SchmidtBADCS15}
M.~W. Schmidt, R.~Babanezhad, M.~O. Ahmed, A.~Defazio, A.~Clifton, and
  A.~Sarkar, ``Non-uniform stochastic average gradient method for training
  conditional random fields,'' in {\em Proceedings of the Eighteenth
  International Conference on Artificial Intelligence and Statistics, {AISTATS}
  2015}, 2015.

\bibitem{RobinsMonro1951}
H.~Robbins and S.~Monro, ``A stochastic approximation method,'' {\em Annals of
  Mathematical Statistics}, vol.~22, no.~3, pp.~400--407, 1951.

\bibitem{Nemirovski-Juditsky-Lan-Shapiro-2009}
A.~Nemirovski, A.~Juditsky, G.~Lan, and A.~Shapiro, ``Robust stochastic
  approximation approach to stochastic programming,'' {\em SIAM Journal on
  Optimization}, vol.~19, no.~4, pp.~1574--1609, 2009.

\bibitem{schmidt2017minimizing}
M.~Schmidt, N.~Le~Roux, and F.~Bach, ``Minimizing finite sums with the
  stochastic average gradient,'' {\em Mathematical Programming}, vol.~162,
  no.~1-2, pp.~83--112, 2017.

\bibitem{shalev2013stochastic}
S.~Shalev-Shwartz and T.~Zhang, ``Stochastic dual coordinate ascent methods for
  regularized loss minimization,'' {\em Journal of Machine Learning Research},
  vol.~14, no.~Feb, pp.~567--599, 2013.

\bibitem{UCDC}
P.~Richt\'{a}rik and M.~Tak\'{a}\v{c}, ``Iteration complexity of randomized
  block-coordinate descent methods for minimizing a composite function,'' {\em
  Mathematical Programming}, vol.~144, pp.~1--38, 2014.

\bibitem{johnson2013accelerating}
R.~Johnson and T.~Zhang, ``Accelerating stochastic gradient descent using
  predictive variance reduction,'' in {\em Advances in Neural Information
  Processing Systems}, pp.~315--323, 2013.

\bibitem{S2GD}
J.~Kone\v{c}n\'{y} and P.~Richt\'{a}rik, ``{S2GD}: {S}emi-stochastic gradient
  descent methods,'' {\em Frontiers in Applied Mathematics and Statistics},
  pp.~1--14, 2017.

\bibitem{defazio2014saga}
A.~Defazio, F.~Bach, and S.~Lacoste-Julien, ``{SAGA}: A fast incremental
  gradient method with support for non-strongly convex composite objectives,''
  in {\em Advances in neural information processing systems}, pp.~1646--1654,
  2014.

\bibitem{pegasos2}
M.~Tak{\'a}\v{c}, A.~Bijral, P.~Richt{\'a}rik, and N.~Srebro, ``Mini-batch
  primal and dual methods for {SVM}s,'' in {\em Proceedings of the 30th
  International Conference on Machine Learning}, 2013.

\bibitem{PCDM}
P.~Richt\'{a}rik and M.~Tak\'{a}\v{c}, ``Parallel coordinate descent methods
  for big data optimization,'' {\em Mathematical Programming}, vol.~156,
  no.~1-2, pp.~433--484, 2016.

\bibitem{PCD_complexity}
R.~Tappenden, M.~Tak{\'a}{\v{c}}, and P.~Richt{\'a}rik, ``On the complexity of
  parallel coordinate descent,'' {\em Optimization Methods and Software},
  vol.~32, no.~2, pp.~372--395, 2018.

\bibitem{mS2GD}
J.~Kone\v{c}n\'{y}, J.~Lu, P.~Richt\'{a}rik, and M.~Tak\'{a}\v{c}, ``Mini-batch
  semi-stochastic gradient descent in the proximal setting,'' {\em IEEE Journal
  of Selected Topics in Signal Processing}, 2016.

\bibitem{ASDCA}
S.~Shalev-Shwartz and T.~Zhang, ``Accelerated mini-batch stochastic dual
  coordinate ascent,'' in {\em Advances in Neural Information Processing
  Systems 26}, pp.~378--385, 2013.

\bibitem{dmSDCA}
M.~Tak\'{a}\v{c}, P.~Richt\'{a}rik, and N.~Srebro, ``Distributed mini-batch
  {SDCA},'' {\em to appear in Journal of Machine Learning Research
  (arXiv:1507.08322)}, pp.~1--15, 2015.

\bibitem{csiba2016importance}
D.~Csiba and P.~Richt{\'a}rik, ``Importance sampling for minibatches,'' {\em to
  appear in Journal of Machine Learning Research (arXiv:1602.02283)}, 2016.

\bibitem{hofmann2015variance}
T.~Hofmann, A.~Lucchi, S.~Lacoste-Julien, and B.~McWilliams, ``Variance reduced
  stochastic gradient descent with neighbors,'' in {\em Advances in Neural
  Information Processing Systems}, pp.~2305--2313, 2015.

\bibitem{richtarik2016optimal}
P.~Richt{\'a}rik and M.~Tak{\'a}{\v{c}}, ``On optimal probabilities in
  stochastic coordinate descent methods,'' {\em Optimization Letters}, vol.~10,
  no.~6, pp.~1233--1243, 2016.

\bibitem{QUARTZ}
Z.~Qu, P.~Richt\'{a}rik, and T.~Zhang, ``Quartz: Randomized dual coordinate
  ascent with arbitrary sampling,'' in {\em Advances in Neural Information
  Processing Systems 28}, pp.~865--873, 2015.

\bibitem{ALPHA}
Z.~Qu and P.~Richt{\'a}rik, ``Coordinate descent with arbitrary sampling {I}:
  algorithms and complexity,'' {\em Optimization Methods and Software},
  vol.~31, no.~5, pp.~829--857, 2016.

\bibitem{allen2016even}
Z.~Allen-Zhu, Z.~Qu, P.~Richt{\'a}rik, and Y.~Yuan, ``Even faster accelerated
  coordinate descent using non-uniform sampling,'' in {\em International
  Conference on Machine Learning}, pp.~1110--1119, 2016.

\bibitem{APPROX}
O.~Fercoq and P.~Richt\'{a}rik, ``Accelerated, parallel and proximal coordinate
  descent,'' {\em SIAM Journal on Optimization}, vol.~25, pp.~1997--2023, 2015.

\bibitem{Hydra2}
O.~Fercoq, Z.~Qu, P.~Richt\'{a}rik, and M.~Tak{\'a}{\v{c}}, ``Fast distributed
  coordinate descent for minimizing non-strongly convex losses,'' {\em IEEE
  International Workshop on Machine Learning for Signal Processing}, 2014.

\bibitem{lin2014accelerated}
Q.~Lin, Z.~Lu, and L.~Xiao, ``An accelerated proximal coordinate gradient
  method,'' in {\em Advances in Neural Information Processing Systems},
  pp.~3059--3067, 2014.

\bibitem{allen2017katyusha}
Z.~Allen-Zhu, ``Katyusha: The first direct acceleration of stochastic gradient
  methods,'' in {\em Proceedings of the 49th Annual ACM SIGACT Symposium on
  Theory of Computing}, pp.~1200--1205, ACM, 2017.

\bibitem{Schmidt2009}
M.~W. Schmidt, E.~V.~D. Berg, M.~P. Friedlander, and K.~Murphy, ``Optimizing
  costly functions with simple constraints: A limited-memory projected
  quasi-newton algorithm,'' {\em International Conference on Artificial
  Intelligence and Statistics}, vol.~5, pp.~456--463, 2009.

\bibitem{Schmidt2011a}
M.~Schmidt, D.~Kim, and S.~Sra, ``Projected {N}ewton-type methods in machine
  learning,'' in {\em Optimization for Machine Learning} (S.~Sra, S.~Nowozin,
  and S.~J. Wright, eds.), pp.~1--25, The MIT press, 2011.

\bibitem{Erdogdu2015nips}
M.~A. Erdogdu and A.~Montanari, ``Convergence rates of subampled {N}ewton
  methods,'' in {\em Advances in Neural Information Processing Systems 28}
  (C.~Cortes, N.~D. Lawrence, D.~D. Lee, M.~Sugiyama, and R.~Garnett, eds.),
  pp.~3052--3060, Curran Associates, Inc., 2015.

\bibitem{GowerGold2016}
R.~M. Gower, D.~Goldfarb, and P.~Richt\'{a}rik, ``Stochastic block {BFGS}:
  Squeezing more curvature out of data,'' in {\em Proceedings of the 33rd
  International Conference on Machine Learning}, 2016.

\bibitem{moritz2016linearly}
P.~Moritz, R.~Nishihara, and M.~I. Jordan, ``A linearly-convergent stochastic
  {L-BFGS} algorithm,'' in {\em International Conference on Artificial
  Intelligence and Statistics}, pp.~249--258, 2016.

\bibitem{GRB-HessianSVRG}
R.~M. Gower, N.~Le~Roux, and F.~Bach, ``Tracking the gradients using the
  {H}essian: {A} new look at variance reducing stochastic methods,'' {\em
  Proceedings of the 21st International Conference on Artificial Intelligence
  and Statistics}, vol.~84, pp.~707--715, 2018.

\bibitem{acceleratedBFGSrules}
R.~M. Gower, F.~Hanzely, F.~Richt\'{a}rik, and S.~Stich, ``Accelerated
  stochastic matrix inversion: general theory and speeding up {BFGS} rules for
  faster second-order optimization,'' {\em arXiv:1801.05490}, 2018.

\bibitem{2018arXiv180502632G}
R.~M. Gower, P.~Richt{\'a}rik, and F.~Bach, ``Stochastic quasi-gradient
  methods: variance reduction via {J}acobian sketching,'' {\em
  ArXiv:1805.02632}, 2018.

\bibitem{SDNA}
Z.~Qu, P.~Richt\'{a}rik, M.~Tak\'{a}\v{c}, and O.~Fercoq, ``{SDNA:} stochastic
  dual {N}ewton ascent for empirical risk minimization,'' in {\em Proceedings
  of The 33rd International Conference on Machine Learning}, pp.~1823--1832,
  2016.

\bibitem{SCSG}
L.~Lei and M.~I. Jordan, ``Less than a single pass: Stochastically controlled
  stochastic gradient,'' in {\em PMLR: Proceedings of Machine Learning Research
  (AISTATS 2017)}, vol.~54, 2017.

\bibitem{chang2011libsvm}
C.-C. Chang and C.-J. Lin, ``{LIBSVM}: a library for support vector machines,''
  {\em ACM Transactions on Intelligent Systems and Technology (TIST)}, vol.~2,
  no.~3, p.~27, 2011.

\bibitem{NesterovBook}
Y.~Nesterov, {\em Introductory Lectures on Convex Optimization: A Basic Course
  (Applied Optimization)}.
\newblock Kluwer Academic Publishers, 2004.

\bibitem{varga1954}
R.~S. Varga, ``Eigenvalues of circulant matrices,'' {\em Pacific J. Math.},
  vol.~1, pp.~151--160., 1954.

\end{thebibliography}
